\documentclass[reqno]{amsart}
\usepackage{amssymb,amsmath,enumerate,mathrsfs}

\numberwithin{equation}{section}

\newtheorem{theorem}[equation]{Theorem}
\newtheorem{proposition}[equation]{Proposition}
\newtheorem{lemma}[equation]{Lemma}
\newtheorem{corollary}[equation]{Corollary}

\theoremstyle{definition}
\newtheorem{definition}[equation]{Definition}
\newtheorem{notation}[equation]{Notation}
\newtheorem{remark}[equation]{Remark}

\def\Vsp{\rule[-1ex]{0ex}{3.5ex}}
\DeclareMathOperator{\sym}{ \sigma\!\!\!\sigma}

\def\C{\mathbb C}

\def\Dom{\mathcal D}
\def\K{\mathcal K}
\def\L{\mathscr L}
\def\N{\mathbb N}
\def\R{\mathbb R}
\def\S{\mathscr S}
\def\Sing{\mathcal E}

\def\csym{\,{}^c\!\sym}
\def\csymb{\,{}^c\!\sym_{\partial}}
\def\cT{\,{}^c T}
\def\cpi{\,{}^c\hspace{-1.5pt}\pi}
\def\bsym{\,{}^b\!\sym}
\def\bsymb{\,{}^b\!\sym_{\partial}}
\def\bT{\,{}^b T}
\def\bpi{\,{}^b\hspace{-1.5pt}\pi}

\def\e{\mathrm e}
\def\embed{\hookrightarrow}
\def\eps{\varepsilon}
\def\id{I}
\def\m{\mathfrak m}
\def\open#1{{\mathring{#1}}}
\def\s{\mathrm s}
\def\st{;\;}
\def\set#1{\{#1\}}
\def\vp{\varphi}
\def\minus{\backslash}

\DeclareMathOperator{\Ind}{ind}
\DeclareMathOperator{\Diff}{Diff}

\DeclareMathOperator{\spec}{spec}

\DeclareMathOperator{\Hom}{Hom}

\begin{document}
\title{Resolvents of elliptic boundary problems on conic manifolds}
\author{Thomas Krainer}
\address{Institut f\"ur Mathematik\\ Universit\"at Potsdam\\ Postfach 60 15 53\\
D-14415 Potsdam, Germany}
\email{krainer@math.uni-potsdam.de}
\begin{abstract}
We prove the existence of sectors of minimal growth for realizations of boundary value problems
on conic manifolds under natural ellipticity conditions. Special attention is devoted to
the clarification of the analytic structure of the resolvent.
\end{abstract}

\subjclass[2000]{Primary: 35J40; Secondary: 58J50, 58J05, 35J70, 35K35}
\keywords{Resolvents, boundary value problems, manifolds with conical singularities, spectral theory,
elliptic and parabolic equations, parametrices}

\maketitle


\section{Introduction}
\label{sec-Introduction}

The present article is a continuation of the investigation of resolvents of elliptic operators on
conic manifolds from \cite{GKM1,GKM2} to the case of manifolds with boundary and
realizations of operators under boundary conditions. Our principal focus of interest are resolvents of
boundary value problems satisfying a parameter-dependent ellipticity condition that resembles
the Shapiro-Lopatinsky condition.

While the study of operators on a conic manifold without boundary is mainly motivated by questions
from spectral theory and geometric analysis, the analysis of boundary value problems has in addition
a wide range of applications in partial differential equations.
In particular, results about the structure and growth of resolvents of operators with respect to the
spectral parameter have immediate consequences as regards the existence, uniqueness, and maximal
regularity of solutions to parabolic linear and semilinear equations.
Hence this paper naturally belongs to the study of partial differential equations in nonsmooth
domains, a subject which due to its importance in models from applications has recently attracted
increased attention (see \cite{KaSchu03}, \cite{KozMazRos}, \cite{MimiTay}, \cite{MiNistor},
\cite{MiTayVas}, \cite{Nistor}, to mention only a few).
Our results, in particular, give a fairly complete picture about the existence of sectors of
minimal growth for $L^2$-based realizations of general elliptic boundary value problems in domains
with cone-like singularities on the boundary. They seem to be new also when specialized to second
order equations under Dirichlet, Neumann, or oblique derivative boundary conditions.

We begin this paper with a brief discussion on manifolds with boundary and conical singularities.
We recall the definitions and some properties of totally characteristic and cone differential
operators, and give a short description of their symbols. We consider boundary value problems for
both of these operator classes and give the suitable definition of parameter-dependent ellipticity (called
$b$- and $c$-ellipticity). Our primary focus, however, are cone operators.

In case of boundary value problems for totally characteristic operators, there is just one realization of
an elliptic operator, and the $b$-ellipticity with parameter in a sector $\Lambda \subset \C$ already
implies that $\Lambda$ is a sector of minimal growth for this realization.
The situation of boundary value problems for cone operators is completely different. There
is a variety of domains for every boundary condition, cf. \cite{CoSeiSchr04}, \cite{GiMe01}, \cite{Le97}.
Each of these domains can be characterized by the asymptotic behavior of its elements near the singularities,
and $c$-ellipticity with parameter in a sector $\Lambda$ is not sufficient
to insure that $\Lambda$ is a sector of minimal growth for a given domain ($c$-ellipticity just entails
Fredholm solvability). Similar to the boundaryless case in \cite{GKM1,GKM2}, it turns out that an
additional spectral condition needs to be required for a model boundary value problem on the model cone
with an associated domain. The model boundary value problem is obtained by ``freezing the coefficients''
at the singularities, and the additional spectral condition can be regarded as a kind of
Shapiro-Lopatinsky condition, but associated with the ``singular boundary'', coming from blowing up the
conical points.
The arbitrariness in the choice of a domain is the source of substantial difficulties which cannot be
overcome merely by considering spaces with weights, a technique that is widely used in the literature since
Kondratiev's seminal paper \cite{K1}.

In the Sections \ref{sec-Realizations}--\ref{sec-AssociatedDomains} we discuss the
domains of realizations of parameter-dependent elliptic boundary value problems for cone operators,
the associated model boundary value problem, and the link between these two.

An important component of this paper is the parametrix construction that we perform in
Section \ref{sec-Parametrix}. Technically, our approach makes use of Boutet de Monvel's calculus
away from the singularities (cf. \cite{GrubbBuch}, \cite{SzWiley98}), several aspects of pseudodifferential
boundary value problems on conic manifolds without parameters (cf. \cite{SchroSchu}), and, near
the singularities, we employ in addition some techniques from the edge-calculus as studied in the
monograph \cite{KaSchu03}.

Finally, in Section \ref{sec-Resolvents}, we use our parametrix to describe the pseudodifferential
structure of the resolvent and prove in Theorem \ref{Maintheorem} the existence of sectors of minimal
growth for realizations of boundary value problems for cone operators under natural
ellipticity conditions, the main result of this note:

\medskip

\noindent
{\bfseries Theorem:}\/ {\itshape Let $A$ be a cone operator of order $m>0$, and $T$ be a vector of
boundary conditions for $A$. Assume that the associated boundary value problem is $c$-elliptic
with parameter in the closed sector $\Lambda \subset \C$, see Section \ref{sec-Realizations}.

Let $\Dom_{\min}(A_T) \subset \Dom(A_T) \subset \Dom_{\max}(A_T)$ be any domain of $A$ under
the boundary condition $Tu = 0$ in a weighted $L_b^2$-space on the manifold, and let $\Dom_{\wedge,\min}(A_{\wedge,T_{\wedge}}) \subset
\Dom_{\wedge}(A_{\wedge,T_{\wedge}}) \subset \Dom_{\wedge,\max}(A_{\wedge,T_{\wedge}})$ be the
associated domain to $\Dom(A_T)$ for the model operator $A_{\wedge}$ under the boundary
condition $T_{\wedge}u = 0$ on the model cone as described
in the Sections \ref{sec-Modeloperator} and \ref{sec-AssociatedDomains}.

If $\Lambda$ is a sector of minimal growth for $A_{\wedge}$ with domain $\Dom_{\wedge}(A_{\wedge,T_{\wedge}})$,
then $\Lambda$ is a sector of minimal growth for the operator $A$ with domain $\Dom(A_T)$
(see Definition~\ref{Sectorminimalgrowthdef} below).
Moreover, the resolvent of $A$ with this domain can be written in the form
$$
\bigl(A_{\Dom}-\lambda\bigr)^{-1} = {\mathcal B}_T(\lambda) + \bigl(A_{\Dom}-\lambda\bigr)^{-1}\Pi_T(\lambda),
$$
where ${\mathcal B}_T(\lambda)$ is a parameter-dependent parametrix of $A-\lambda$ taking values in the
minimal domain $\Dom_{\min}(A_T)$, and $\Pi_T(\lambda)$ is a regularizing projection operator onto a
(finite dimensional) complement of the range of $A-\lambda$ on $\Dom_{\min}(A_T)$ (see Section~\ref{sec-Parametrix}).
}

\medskip

The idea for proving this theorem is the construction of an invertible abstract reference problem
(Theorem \ref{MainTheoremParametrix}) that is used to reduce
the resolvent constructions to the finite dimensional contribution beyond the minimal domain.
In this sense, we perform a reduction to the ``singular boundary'', and the resulting operator
family plays technically a similar role as, e.g., the Dirichlet-to-Neumann map in regular boundary value
problems. However, a canonical reference domain for the operator $A$ with boundary condition $Tu = 0$,
which could be regarded as a ``Dirichlet extension'' with respect to the singular boundary, does not
exist.

Compared to the boundaryless case in \cite{GKM1,GKM2}, managing the boundary conditions requires special
care and forces setting up a much more elaborated machinery. By consistently working with full operator
matrices (with domains that are associated with the inhomogeneous boundary value problem), we are able to
transfer some methods from the boundaryless situation.
This also makes it possible to concentrate on the essential singular part of the problem when coming
to the resolvent constructions in Section \ref{sec-Resolvents}.

\medskip

\emph{Acknowledgement}\/: This article emerged from a collaboration with Juan Gil (Penn State Altoona) and Gerardo Mendoza
(Temple University, Philadelphia). The author wishes to express his gratitude for many invaluable
discussions on the subject of this paper.

\medskip

Throughout this article, let
$$
\Lambda = \{z=re^{i\theta} \st r \geq 0, \; |\theta-\theta_0| \leq \eps\}
$$
with $\theta_0 \in \R$ and $\eps > 0$ be a closed sector in $\C$.

\begin{definition}\label{Sectorminimalgrowthdef}
Let $H$ be a Hilbert space and $A : \Dom(A) \subset H \to H$ be unbounded, densely defined, and closed, i.e. $\Dom(A)$ is
complete in the graph norm
\begin{equation}\label{Graphnormdef}
\|u\|_A = \|u\|_H + \|Au\|_H.
\end{equation}
We call $\Lambda$ a sector of minimal growth for the operator $A$ (with domain $\Dom(A)$) if
$$
A-\lambda : \Dom(A) \to H
$$
is bijective for large $|\lambda| > 0$ in $\Lambda$, and if the following equivalent norm
estimates for the resolvent $(A-\lambda)^{-1} : H \to \Dom(A)$ are satisfied:
\begin{enumerate}[i)]
\item $\|(A-\lambda)^{-1}\|_{\L(H)} = O(|\lambda|^{-1})$ as $|\lambda| \to \infty$.
\item $\|(A-\lambda)^{-1}\|_{\L(H,\Dom(A))} = O(1)$ as $|\lambda| \to \infty$.
\end{enumerate}
\end{definition}


\section{Manifolds with boundary and conical singularities}
\label{sec-Manifolds}

\begin{definition}\label{conicmanifold}
A compact conic manifold with boundary is a second countable compact
Hausdorff topological space $M_{\textup{sing}}$ such that there exists a finite
subset $S$ with the following properties:
\begin{enumerate}[i)]
\item $M_{\textup{sing}}\setminus S$ is a smooth manifold with boundary.
\item Every $s \in S$ has a neighborhood $U(s) \subseteq M_{\textup{sing}}$ which
is homeomorphic to a neighborhood $\tilde{U}$ of
$$
\tilde{s} = (\{0\}\times\overline{Y})/(\{0\}\times\overline{Y})
$$
in $(\overline{\R}_+\times\overline{Y}) / (\{0\}\times\overline{Y})$,
where $\overline{Y}$ is a compact smooth manifold with boundary, and the
homeomorphism restricts to a diffeomorphism $U(s) \minus \{s\}
\cong \tilde{U} \minus \{\tilde{s}\}$.
\end{enumerate}
The set $S$ is the singular set in $M_{\textup{sing}}$, the elements of $S$ are
called conical points.
\end{definition}

As in the boundaryless case, analysis on conic manifolds with boundary is performed
away from the conical points. Consequently, by eventually passing to
$M_{\textup{sing}} / S$, we can and will assume henceforth that $M_{\textup{sing}}$ has
only one conical point $s$. Note that the manifold $\overline{Y}$ in Definition \ref{conicmanifold}
is not assumed to be connected.
 
It is evident from the definition that $N_{\textup{sing}}$, the boundary of $M_{\textup{sing}}$, is
a compact conic manifold without boundary and conical point $s$.

\begin{definition}\label{conestructure}
A cone structure on $M_{\textup{sing}}$ is a maximal atlas consisting of a
differential structure for the smooth manifold with boundary
$M_{\textup{sing}} \setminus \{s\}$, as well as coordinate neighborhoods of the conical
point $s$ of the form $U(s)$ from Definition \ref{conicmanifold},
where away from $s$ the coordinate changes are $C^{\infty}$-diffeomorphisms,
and the coordinate change
$$
\R_+\times\overline{Y} \supseteq \tilde{U} \minus\{\tilde{s}\} \cong \tilde{V}\minus\{\tilde{s}\} \subseteq \R_+\times\overline{Y}
$$
of any two charts near the conical point $s$ extends to a $C^{\infty}$-mapping
$$
((-\infty,0]\times\overline{Y}) \cup (\tilde{U} \minus\{\tilde{s}\}) \longrightarrow
((-\infty,0]\times\overline{Y}) \cup (\tilde{V} \minus\{\tilde{s}\})
$$
between these open subsets of $\R\times\overline{Y}$. Note, in particular, that by
continuity the cocycle property of coordinate changes near $s$ is preserved up to the origin,
i.e., up to $\{0\}\times\overline{Y}$.

Any cone structure on $M_{\textup{sing}}$ gives rise to a unique cone structure on the boundary $N_{\textup{sing}}$. We will always assume that a cone structure on
$M_{\textup{sing}}$ is fixed, and the boundary will be given the induced cone structure.
\end{definition}

Let $\overline{M}$ be the compact space obtained from $M_{\textup{sing}}$ by
blowing up the conical point $s$ to $\overline{Y}$. Note that $\overline{M}$
and the canonical embedding $\overline{Y} \embed \overline{M}$ are invariants
of the cone structure, and each local chart of $M_{\textup{sing}}$ near
the conical point $s$ gives rise to a collar neighborhood of $\overline{Y}$
in $\overline{M}$.

We have a canonical identification $\overline{M}/\overline{Y} \cong M_{\textup{sing}}$
as compact conic manifolds with boundary. The double
$$
2\overline{M} = \overline{M} \cup_{\overline{Y}} \overline{M}
$$
is a compact smooth manifold with boundary, where the $C^{\infty}$-structure
is inherited from collar neighborhoods of $\overline{Y}$ in $\overline{M}$.
Evidently, $\overline{N} = \partial(2\overline{M})\cap\overline{M}$ is the blow-up
of $N_{\textup{sing}}$, a compact manifold with boundary $\partial\overline{Y}$.

Let us fix a defining function $x \in C^{\infty}(2\overline{M})$ for $\overline{Y}$ with
$x > 0$ on $\overline{M} \minus \overline{Y}$.

\begin{definition}\label{coneoperators}
By $\Diff_b^m(\overline{M})$ we denote the restrictions to $\overline{M}$ of the $m$-th order
differential operators on $2\overline{M}$ which are totally characteristic with
respect to $\overline{Y}$. Thus $\Diff_b^\ast(\overline{M})$ is the enveloping algebra
generated by the restrictions to $\overline{M}$ of the vector fields on $2\overline{M}$
which are tangent to $\overline{Y}$ and $C^{\infty}(2\overline{M})$.

Observe that $\overline{N} \minus \partial\overline{Y}$, the regular part of the boundary of
$\overline{M}$, is not necessarily characteristic for the elements of $\Diff_b^\ast(\overline{M})$.

Correspondingly, let $\Diff_b^m(\overline{M};E,F)$ be the space of $m$-th order totally
characteristic differential operators acting in sections of the bundles $E$ and $F$.
Note that we consider here complex vector bundles that are restrictions of smooth bundles on
$2\overline{M}$ to $\overline{M}$.

The operators $A \in x^{-m}\Diff_b^m(\overline{M};E,F)$ are the cone operators of
order $m$. If $(x,y)$ are local coordinates near $p \in \overline{Y} \subset \overline{M}$ with
$x \geq 0$ on $\overline{M}$ and $x = 0$ on $\overline{Y}$, then $A$ takes the form
$$
A = x^{-m}\sum\limits_{k+|\alpha|\leq m}a_{k,\alpha}(x,y)D_y^{\alpha}(xD_x)^k
$$
with coefficients $a_{k,\alpha}(x,y)$ that are smooth up to $x = 0$.
\end{definition}

Totally characteristic operators $A \in \Diff_b^m(\overline{M};E,F)$, $m \in \N_0$, have
an invariant $b$-principal symbol on the compressed cotangent bundle
$\bT^*\overline{M}$, see \cite{RBM1, RBM2}.
Recall that $\bT\overline{M}$ is the bundle on $\overline{M}$ whose smooth sections are the
restrictions of the vector fields on $2\overline{M}$ to $\overline{M}$ which are tangent to
$\overline{Y}$. The compressed cotangent bundle $\bT^*\overline{M}$ is the dual of
$\bT\overline{M}$.

In \cite{GKM1} the $c$-cotangent bundle was introduced, and it was proved that cone
operators have invariantly defined principal symbols there. Consequently, with an
operator $A \in x^{-m}\Diff_b^m(\overline{M};E,F)$, $m \in \N_0$, we associate its $c$-principal
symbol $\csym(A)$ on $\cT^*\overline{M}$, a section of the bundle
$\Hom(\cpi^*E,\cpi^*F)$. Here $\cpi : \cT^*\overline{M} \to \overline{M}$ is the canonical
projection map.

The definition of the $c$-cotangent bundle $\cT^*\overline{M}$ is similar to the $b$-construction,
and its space of smooth sections are the restrictions of $1$-forms from $2\overline{M}$ to $\overline{M}$
which are, over $\overline{Y}$, sections of the conormal bundle of $\overline{Y}$ in $T^*(2\overline{M})$.

There is a second principal symbol associated with an operator $A$, the $b$- or $c$-principal
boundary symbol $\bsymb(A)$ or $\csymb(A)$, respectively.

\begin{definition}\label{boundarysymbol}
Let $A \in x^{-m}\Diff_b^m(\overline{M};E,F)$, and let $U(p) \subset 2\overline{M}$ be a small
neighborhood of the point $p \in \partial(2\overline{M})\cap\overline{M}$ on the
boundary of $\overline{M}$, and consider local coordinates $(z',z_n)$ centered at $p$ with
$z_n \geq 0$ and $z_n = 0$ on $\partial(2\overline{M})$.
Let $\csym(A)(z',z_n;\zeta',\zeta_n)$ be the local representation of the $c$-principal symbol in
these coordinates, a $(N_- \times N_+)$-matrix function with $N_- = \dim E$ and $N_+ = \dim F$.

The operator family
\begin{equation}\label{localboundarysymbol}
\csymb(A)(z';\zeta') = \csym(A)(z',0;\zeta',D_{z_n}) : \S(\overline{\R}_+)\otimes\C^{N_-} \to
\S(\overline{\R}_+)\otimes\C^{N_+}
\end{equation}
is then a local representation of the $c$-principal boundary symbol $\csymb(A)$ of $A$.

It is more tedious than hard to see that the $c$-principal boundary symbol is invariantly
defined on $\cT^*\overline{N}$, a section of the bundle $\Hom({}^c\S_+ \otimes \cpi^*E,{}^c\S_+ \otimes \cpi^*F)$.
Here ${}^c\S_+$ is a bundle on $\cT^*\overline{N}$ with fiber $\S(\overline{\R}_+)$ which
comes up canonically when changing to a different local representation \eqref{localboundarysymbol}
near $p$.

Analogously, with a totally characteristic operator $A \in \Diff_b^m(\overline{M};E,F)$ we
associate the $b$-principal boundary symbol $\bsymb(A)$, a section of the bundle
$\Hom({}^b\S_+ \otimes \bpi^*E,{}^b\S_+ \otimes \bpi^*F)$ over $\bT^*\overline{N}$.
\end{definition}

\begin{definition}\label{ellipticity}
Let $0 \in \Omega \subset \C$ be a conical subset. The operator family $A - \lambda \in
\Diff_b^m(\overline{M};E)$, $\lambda \in \Omega$, is called $b$-elliptic with parameter, if
$\spec(\bsym(A)(z,\zeta)) \cap \Omega = \emptyset$ for $(z,\zeta) \in \bT^*\overline{M}\minus 0$.

Analogously, for a cone operator $A \in x^{-m} \Diff_b^m(\overline{M};E)$ the family $A - \lambda$
is called $c$-elliptic with parameter $\lambda \in \Omega$ iff $\spec(\csym(A)(z,\zeta)) \cap \Omega = \emptyset$
for $(z,\zeta) \in \cT^*\overline{M}\minus 0$.
\end{definition}

Obviously, $A - \lambda \in x^{-m} \Diff_b^m(\overline{M};E)$ is $c$-elliptic with parameter if
and only if $x^mA - \lambda \in \Diff_b^m(\overline{M};E)$ is $b$-elliptic with parameter, and,
if $\Omega = \{0\}$, $b$- and $c$-ellipticity with parameter reduces to ordinary
$b$- and $c$-ellipticity.


\section{Boundary problems for totally characteristic operators}
\label{sec-Totally}

In this section we consider boundary value problems for totally characteristic operators
$A \in \Diff_b^m(\overline{M};E)$, $m \in \N$.

We assume henceforth that $A-\lambda$ is $b$-elliptic with parameter in the sector $\Lambda \subset \C$.

\begin{proposition}
The $b$-principal boundary symbol
$$
\bsymb(A)(z',\zeta') - \lambda : {}^b\S_+ \otimes \bpi^*E \to {}^b\S_+ \otimes \bpi^*E
$$
is surjective and has finite dimensional kernel for all
$(z',\zeta',\lambda) \in \bigl(\bT^*\overline{N}\times\Lambda\bigr)\minus 0$.

Consequently, the kernels form a vector bundle $\K$ on
$\bigl(\bT^*\overline{N}\times\Lambda\bigr)\minus 0$.
\end{proposition}

For a sufficiently smooth section $u$ of a bundle $F$ on $\overline{M}\minus\overline{Y}$
we denote by $\gamma u$ its restriction to the regular part $\overline{N}\minus\partial\overline{Y}$
of the boundary, which gives rise to the restriction operator $\gamma$.

Now let $B_j \in \Diff_b^{m_j}(\overline{M};E,F_j)$, $m_j < m$ , be totally characteristic,
$j = 1,\ldots,K$. We consider the family of boundary value problems
\begin{equation}\label{BVPTotally}
\begin{aligned}
(A - \lambda)u &= f \quad \text{in } \open{\overline{M}} = \overline{M}\minus\overline{Y}, \\
Tu &= g \quad \text{on } \open{\overline{N}} = \overline{N}\minus\partial\overline{Y}
\end{aligned}
\end{equation}
for the operator $A$, where $T = (\gamma B_1,\ldots,\gamma B_K)^t$.

\begin{definition}\label{bShapiroLopatinsky}
The boundary value problem \eqref{BVPTotally} is called $b$-elliptic with parameter $\lambda \in \Lambda$ if
$$
\begin{pmatrix}
({}^b\gamma_0\otimes\id_{\bpi^*F_1})\bsymb(B_1)(z',\zeta') \\
\vdots \\
({}^b\gamma_0\otimes\id_{\bpi^*F_K})\bsymb(B_K)(z',\zeta')
\end{pmatrix} : \K_{(z',\zeta',\lambda)} \to \bigoplus\limits_{j=1}^K \bpi^*F_j
$$
is bijective for all $(z',\zeta',\lambda) \in \bigl(\bT^*\overline{N}\times\Lambda\bigr)\minus 0$.
Here ${}^b\gamma_0 : {}^b\S_+ \to \C$ is the canonical evaluation map at zero.

Note that this notion of $b$-ellipticity is the appropriate version of the Shapiro-Lopatinsky
condition for families of totally characteristic boundary problems.
\end{definition}

Let $\m > 0$ be a $b$-density on $2\overline{M}$, i.e. $x\m$ is a smooth everywhere
positive density. Let $L_b^2(\overline{M};E)$ be the $L^2$-space of sections of the bundle
$E$ on $\overline{M}$ with respect to $\m$ and a Hermitian inner product on $E$.

For $s \in \N_0$ let $H^s_b(\overline{M};E)$ be the Sobolev space of all $L_b^2$-sections
$u$ such that $Cu \in L_b^2$ for all totally characteristic operators $C$ of order $\leq s$, and let
$H^s_{b,0}(\overline{M};E)$ be the closure of all $C_0^{\infty}$-sections of $E$ in
$H^s_b(\overline{M};E)$. Note that the $C_0^{\infty}$-sections here are supported away from
the boundary $\overline{N}\cup\overline{Y}$.

For $s \in -\N$ let $H^s_b(\overline{M};E)$ be the dual of $H^{-s}_{b,0}(\overline{M};E)$,
and analogously let $H^s_{b,0}(\overline{M};E)$ be the dual of $H^{-s}_b(\overline{M};E)$ with
respect to the sesquilinear pairing induced by the $L^2_b$-inner product. For arbitrary real $s$
we define $H^s_{b,0}(\overline{M};E)$ and $H^s_b(\overline{M};E)$ by interpolation.
Analogously to the boundaryless case, we also consider weighted spaces $x^{\mu}H^s_b$ for
arbitrary weights $\mu \in \R$.

For every $s > -1/2$ we consider the boundary value problem \eqref{BVPTotally} as a family
of continuous operators
\begin{equation}\label{BVPTotSobsp}
\begin{pmatrix}
A - \lambda \\
T
\end{pmatrix} : x^{\mu}H_b^{s+m}(\overline{M};E) \to \begin{array}{c}
x^{\mu}H_b^{s}(\overline{M};E) \\ \oplus \\ \bigoplus\limits_{j=1}^K x^{\mu}H_b^{s+m-m_j-1/2}(\overline{N};F_j)
\end{array}.
\end{equation}
Correspondingly, we consider the realization of $A$ under the boundary condition $Tu = 0$, i.e.
the unbounded operator
$$
A_T : \Dom^s(A_T) \subset x^{\mu}H_b^s(\overline{M};E) \to x^{\mu}H_b^s(\overline{M};E)
$$
with domain $\Dom^s(A_T) = \{u \in x^{\mu}H_b^{s+m}(\overline{M};E) \st Tu = 0 \}$ that
acts like $A$. We sometimes also say that this is the $H^{s+m}$-realization (or just
$H$-realization) of $A$ in order to emphasize that we assume apriori smoothness of order $s+m$ in
$\open{\overline{M}}$ as it is custom in boundary value problems.

It is a part of Theorem \ref{BVPTotallyThm} below that the boundary condition
$$
T : x^{\mu}H_b^{s+m}(\overline{M};E) \to \bigoplus\limits_{j=1}^K x^{\mu}H_b^{s+m-m_j-1/2}(\overline{N};F_j)
$$
is surjective for $s > -1/2$. In particular, the operator \eqref{BVPTotSobsp} is invertible for
some $\lambda \in \Lambda$ if and only if $\lambda \notin \spec(A_T)$.

\begin{theorem}\label{BVPTotallyThm}
Let \eqref{BVPTotally} be $b$-elliptic with parameter $\lambda \in \Lambda$, and let
$\mu \in \R$. There exists $R=R(\mu) > 0$ such that for $\lambda \in \Lambda$,
$|\lambda| \geq R$, the operator \eqref{BVPTotSobsp} is invertible for all
$s > -1/2$.

Moreover, we have $\|(A_T - \lambda)^{-1}\| = O(|\lambda|^{-1})$ as $|\lambda| \to \infty$ for
the $\L(x^\mu L_b^2)$-norm of the resolvent of $A$ with domain
$\Dom(A_T) = \Dom^0(A_T) \subset x^{\mu}L_b^2(\overline{M};E)$, and the norm of the resolvent
in $\L(x^{\mu} H_b^s)$ of realizations in Sobolev spaces of higher regularity is polynomially
bounded as $|\lambda| \to \infty$.
\end{theorem}

As $A$ with domain $\Dom(A_T)$ is closed, we thus obtain that $b$-ellipticity with parameter implies
that $\Lambda$ is a sector of minimal growth for $A$.

The proof of Theorem \ref{BVPTotallyThm} follows by constructing a parametrix of
\eqref{BVPTotally} within a Boutet de Monvel's calculus of parameter-dependent
boundary value problems of totally characteristic pseudodifferential operators (a suitable
modification of the arguments given in Section \ref{sec-Parametrix} up to Proposition
\ref{Parametrix2} will do).

As our interest in this article lies in resolvents and spectral properties of boundary value
problems for cone operators we will not pursue this here.

Despite of the many similarities between cone operators and totally characteristic operators,
the spectral theory for cone operators is much more complicated than the spectral theory of
totally characteristic operators. This is underscored by a comparison of
Theorem \ref{BVPTotallyThm} and Theorem \ref{Maintheorem}.

Assuming parameter-dependent ellipticity, for every weight $\mu \in \R$ there is only one
$H$-realization of a totally characteristic operator $A$ under the boundary condition
$Tu = 0$, and this realization is well-behaved for large parameter values as
Theorem \ref{BVPTotallyThm} shows. In contrast, for every weight $\mu \in \R$ there are many
$H$-realizations of cone operators, and the spectrum of every such realization could be $\C$
(see also \cite{GKM1, GKM2} for a discussion of the boundaryless case).

\medskip

For later purposes, we recall the notion of the conormal symbol associated with a totally characteristic
boundary value problem on $\overline{M}$:

If $u$ is a smooth section of $E$ on $2\overline{M}$ that vanishes on $\overline{Y}$, then also
$Pu$ vanishes on $\overline{Y}$ for every $P \in \Diff_b^*(\overline{M};E,F)$. Consequently, if
$v$ is a section of $E$ on $\overline{Y}$ and $u$ is any extension of $v$, then $(Pu)|_{\overline{Y}}$
does not depend on the choice of the extension. Thus, associated with $P$ there is a differential operator
$$
\hat{P}(0) : C^{\infty}(\overline{Y};E) \to C^{\infty}(\overline{Y};F)
$$
of the same order as $P$. Since $\C \ni \sigma \mapsto x^{-i\sigma}Px^{i\sigma}$ is a family of
totally characteristic operators, we so obtain the operator valued polynomial
$$
\C \ni \sigma \mapsto \hat{P}(\sigma) \in \Diff^*(\overline{Y};E,F),
$$
the conormal symbol associated with $P$. If $P$ is part of a boundary condition $\gamma P$, we associate
with this condition its conormal symbol
$$
\C \ni \sigma \mapsto \hat{\gamma}\hat{P}(\sigma),
$$
where $\hat{\gamma}v$ denotes the restriction of the section $v$ on $\overline{Y}$ to the boundary
$\partial\overline{Y}$. In this way we obtain for each $\lambda \in \Lambda$ the conormal symbol of
the boundary value problem \eqref{BVPTotally}, a family of boundary value problems
\begin{equation}\label{ConormalsymbolTotally}
h(\sigma,\lambda) : C^{\infty}(\overline{Y};E) \to
\begin{array}{c}
C^{\infty}(\overline{Y};E) \\
\oplus \\
\bigoplus\limits_{j=1}^KC^{\infty}(\partial\overline{Y};F_j)
\end{array}
\end{equation}
on $\overline{Y}$ depending holomorphically on $\sigma \in \C$ and $\lambda \in \Lambda$.

Provided that \eqref{BVPTotally} is $b$-elliptic with parameter $\lambda \in \Lambda$, the conormal
symbol \eqref{ConormalsymbolTotally} is for every $\sigma \in \C$ and $\lambda \in \Lambda$ an
elliptic boundary value problem on $\overline{Y}$, which is even elliptic with parameter
$(\sigma,\lambda) \in \{\Im(\sigma) = \alpha\}\times\Lambda$ for every fixed $\alpha \in \R$.


\section{Realizations of boundary problems for cone operators}
\label{sec-Realizations}

Let $A \in x^{-m}\Diff_b^m(\overline{M};E)$, $m \in \N$, be a cone operator such that
$A-\lambda$ is $c$-elliptic with parameter in the sector $\Lambda \subset \C$.

Analogously to the case of totally characteristic operators, we then know that
the $c$-principal boundary symbol
$$
\csymb(A)(z',\zeta') - \lambda : {}^c\S_+ \otimes \cpi^*E \to {}^c\S_+ \otimes \cpi^*E
$$
is surjective and has finite dimensional kernel for all
$(z',\zeta',\lambda) \in \bigl(\cT^*\overline{N}\times\Lambda\bigr)\minus 0$.
Let $\K$ be the bundle of kernels on $\bigl(\cT^*\overline{N}\times\Lambda\bigr)\minus 0$.

Let $B_j \in x^{-m_j}\Diff_b^{m_j}(\overline{M};E,F_j)$, $m_j < m$ , be cone operators,
$j = 1,\ldots,K$, and consider the family of boundary value problems
\begin{equation}\label{BVPCone}
\begin{aligned}
(A - \lambda)u &= f \quad \text{in } \open{\overline{M}} = \overline{M}\minus\overline{Y}, \\
Tu &= g \quad \text{on } \open{\overline{N}} = \overline{N}\minus\partial\overline{Y}
\end{aligned}
\end{equation}
for the operator $A$, where $T = (\gamma B_1,\ldots,\gamma B_K)^t$.

\begin{definition}\label{cShapiroLopatinsky}
The boundary value problem \eqref{BVPCone} is called $c$-elliptic with parameter $\lambda \in \Lambda$ if
$$
\begin{pmatrix}
({}^c\gamma_0\otimes\id_{\cpi^*F_1})\csymb(B_1)(z',\zeta') \\
\vdots \\
({}^c\gamma_0\otimes\id_{\cpi^*F_K})\csymb(B_K)(z',\zeta')
\end{pmatrix} : \K_{(z',\zeta',\lambda)} \to \bigoplus\limits_{j=1}^K \cpi^*F_j
$$
is bijective for all $(z',\zeta',\lambda) \in \bigl(\cT^*\overline{N}\times\Lambda\bigr)\minus 0$,
where ${}^c\gamma_0 : {}^c\S_+ \to \C$ is evaluation at zero.

Similar to the case of totally characteristic operators, $c$-ellipticity is the appropriate
version of the Shapiro-Lopatinsky condition for cone operators.
\end{definition}

\begin{lemma}\label{TotConEllipt}
The boundary value problem \eqref{BVPCone} with cone operators $A$ and $B_j$, $j=1,\ldots,K$,
is $c$-elliptic with parameter $\lambda \in \Lambda$ if and only if the boundary value problem
\begin{alignat*}{2}
((x^mA) - \lambda)u &= f & &\text{in } \open{\overline{M}}, \\
\gamma (x^{m_j}B_j)u &= g_j \quad& &\text{on } \open{\overline{N}}, \; j=1,\ldots,K
\end{alignat*}
is $b$-elliptic with parameter $\lambda \in \Lambda$ in the sense of Definition \ref{bShapiroLopatinsky}.
\end{lemma}

Our primary concern is to investigate the spectral properties of $c$-elliptic boundary value problems
under the assumption of parameter-dependent ellipticity, i.e. we investigate the operator family
\begin{equation}\label{BVPConeSobsp}
\begin{pmatrix}
A - \lambda \\
T
\end{pmatrix} : \Dom^s\binom{A}{T} \subset x^{\mu}H_b^{s}(\overline{M};E) \to \begin{array}{c}
x^{\mu}H_b^{s}(\overline{M};E) \\ \oplus \\ \bigoplus\limits_{j=1}^K x^{\mu+m-m_j}H_b^{s+m-m_j-1/2}(\overline{N};F_j)
\end{array}
\end{equation}
for $s > -1/2$ and some weight $\mu \in \R$, as well as the behavior of the associated
family of unbounded operators
\begin{equation}\label{BVPConeReal}
A - \lambda : \Dom^s(A_T) \subset x^{\mu}H_b^{s}(\overline{M};E) \to
x^{\mu}H_b^{s}(\overline{M};E)
\end{equation}
with domain $\Dom^s(A_T) = \Dom^s\binom{A}{T} \cap \ker T$. Of particular interest is of course
the case $s = 0$, i.e. the $x^{\mu}L^2_b$-realization of the operator $A$ under the boundary
condition $Tu = 0$.

The domain in \eqref{BVPConeSobsp} can be any intermediate space
$$
\Dom_{\min}^s\binom{A}{T} \subset \Dom^s\binom{A}{T} \subset \Dom_{\max}^s\binom{A}{T}
$$
of the minimal and maximal $x^{\mu}H_b^{s+m}$-domains
\begin{align*}
\Dom_{\max}^s\binom{A}{T} &= \Biggl\{u \in x^{\mu}H_b^{s+m}(\overline{M};E) \st
\begin{aligned}[c]
Au &\in x^{\mu}H_b^s(\overline{M};E), \\
\gamma B_ju &\in x^{\mu+m-m_j}H_b^{s+m-m_j-1/2}(\overline{N};F_j) \\
&\qquad \textup{for } j=1,\ldots,K
\end{aligned}\Biggr\}, \\
\Dom_{\min}^s\binom{A}{T} &= \Dom_{\max}^s\binom{A}{T} \cap \bigcap\limits_{\eps > 0}
x^{\mu+m-\eps}H_b^{s+m}(\overline{M};E).
\end{align*}
As conjugation of \eqref{BVPConeSobsp} with the weight function $x^{\delta}$ for any
$\delta \in \R$ gives rise to a unitary equivalent parameter-dependent $c$-elliptic boundary
value problem of the form \eqref{BVPCone} in the corresponding shifted function spaces, we can
without loss of generality base all our investigations on the weight $\mu = -m/2$.
Moreover, we usually write $\Dom^s = \Dom^s\binom{A}{T}$ as well as $\Dom = \Dom^0\binom{A}{T}$.

A discussion of domains and adjoints of $c$-elliptic boundary value problems and normal
boundary conditions is given in \cite{CoSeiSchr04}, generalizing previous results in
\cite{GiMe01} in the boundaryless case. In contrast to the mere $c$-elliptic case, our situation of
parameter-dependent $c$-ellipticity makes it possible to avoid a technical discussion
of the issue of normality.

The next proposition follows analogous to the boundaryless case from a corresponding
analysis in the Mellin image using the conormal symbols of \eqref{BVPCone},
see \cite{Le97}, \cite{GiMe01}. The proof of the Fredholmness in part iv) follows
by employing a standard parametrix (without parameters) of elliptic boundary value problems
on the cone, see, e.g., \cite{SchroSchu}, \cite{CoSeiSchr04}.

\begin{proposition}\label{Domainproperties}
Assume that \eqref{BVPCone} is $c$-elliptic with parameter in some closed
sector $\Lambda \subset \C$.
\begin{enumerate}[i)]
\item $\Dom^s_{\max}$ is complete in the norm
\begin{equation}\label{NormwithBoundCond}
\|u\|_{A_T} = \|u\|_{x^{-m/2}H_b^{s+m}} + \|Au\|_{x^{-m/2}H_b^s} + \sum\limits_{j=1}^K \|\gamma B_ju\|_{x^{m/2-m_j}H_b^{s+m-m_j-1/2}},
\end{equation}
and $\Dom^s_{\min} \subset \Dom^s_{\max}$ is a closed subspace of finite codimension.
\item We have
$$
x^{m/2}H_b^{s+m} \embed \Dom^s_{\min} \embed \Dom^s_{\max} \embed x^{-m/2+\eps}H_b^{s+m}
$$
for some $\eps > 0$ with continuous embeddings. In particular, the embedding
$\Dom^s_{\max} \embed x^{-m/2}H_b^s$ is compact.
\item $C_0^{\infty}(\open{\overline{M}};E) \subset \Dom^s_{\min}$ is a dense subspace.
\item For every $\lambda \in \C$ the boundary value problem \eqref{BVPConeSobsp} is Fredholm
with index independent of $\lambda$ and $s > -1/2$, and we have the following relative
index formula
\begin{equation}\label{RelIndex}
\Ind\binom{A-\lambda}{T}_{\Dom^s} = \Ind\binom{A}{T}_{\Dom^s_{\min}} + \dim \Dom^s/\Dom^s_{\min}.
\end{equation}
Here the subscripts refer to the corresponding domains.
\end{enumerate}
\end{proposition}

The quotient $\Dom^s_{\max}/\Dom^s_{\min}$ is actually independent of $s > -1/2$ and
can be identified with a space of singular functions. We will come back to this in
Section \ref{sec-AssociatedDomains} soon (see also \cite{GKM2}, Section 6).

From Lemma \ref{TotConEllipt} and Theorem \ref{BVPTotallyThm} we obtain that the boundary
condition
$$
T : x^{m/2}H_b^{s+m}(\overline{M};E) \to \bigoplus\limits_{j=1}^K x^{m/2-m_j}H_b^{s+m-m_j-1/2}(\overline{N};F_j)
$$
is surjective, and necessarily so is its extension to $\Dom^s$ by Proposition \ref{Domainproperties}.
Consequently, for every $\lambda \in \C$,
$$
A_T - \lambda : \Dom^s(A_T) \to x^{-m/2}H_b^s(\overline{M};E)
$$
is Fredholm with index
$$
\Ind(A-\lambda)_{\Dom^s(A_T)} = \Ind A_{\Dom^s(A_T)} = \Ind\binom{A}{T}_{\Dom^s},
$$
and a necessary condition for $A$ with domain $\Dom(A_T)$ to admit nonempty resolvent set
is that $\Ind(A_{T,\min}) \leq 0$ and $\Ind(A_{T,\max}) \geq 0$ (in \cite{GKM1} such issues
are discussed from a fairly abstract perspective, and many of the results therefore apply also
to the situation under study in this article). Moreover, \eqref{BVPConeSobsp} is invertible for some
$\lambda$ if and only if \eqref{BVPConeReal} is bijective, i.e. if and only if $\lambda \notin \spec(A_T)$.

Let us formulate an immediate consequence of these observations (note, in particular, that this
constitutes a substantial difference to the totally characteristic case):

\begin{corollary}
Let \eqref{BVPCone} be $c$-elliptic with parameter in the closed sector $\Lambda \subset \C$. Then
either the spectrum of the operator $A_T : \Dom^s(A_T) \subset x^{-m/2}H_b^s(\overline{M};E) \to x^{-m/2}H_b^s(\overline{M};E)$
is discrete or it is all of $\C$, and a necessary condition for the spectrum to be discrete is that
$\Ind A_{\Dom^s(A_T)} = 0$.
\end{corollary}

\begin{lemma}\label{latticeisomanifold}
The mapping $\Dom^s\binom{A}{T} \to \Dom^s(A_T) = \Dom^s\binom{A}{T}\cap \ker T$ is a bijection of the lattice of intermediate
spaces
\begin{alignat*}{2}
\Dom^s_{\min}\binom{A}{T} &\subset \Dom^s\binom{A}{T} & &\subset \Dom^s_{\max}\binom{A}{T} \\
\intertext{onto the lattice of intermediate spaces}
\Dom^s_{\min}(A_T) &\subset \Dom^s(A_T) & &\subset \Dom^s_{\max}(A_T),
\end{alignat*}
where
\begin{align*}
\Dom^s_{\min}(A_T) &= \bigl\{u \in \bigcap\limits_{\eps>0}x^{m/2-\eps}H_b^{s+m}(\overline{M};E) \st
Au \in x^{-m/2}H_b^s(\overline{M};E) \textup{ and } Tu = 0\bigr\}, \\
\Dom^s_{\max}(A_T) &= \bigl\{u \in x^{-m/2}H_b^{s+m}(\overline{M};E) \st
Au \in x^{-m/2}H_b^s(\overline{M};E) \textup{ and } Tu = 0\bigr\}.
\end{align*}
We have $\dim\Dom^s_{\max}\binom{A}{T}/\Dom^s_{\min}\binom{A}{T} =
\dim\Dom^s_{\max}(A_T)/\Dom^s_{\min}(A_T)$. More precisely, given a basis
$$
s_j+\Dom^s_{\min}\binom{A}{T}, \quad j = 1,\ldots,M,
$$
of $\Dom^s_{\max}\binom{A}{T}/\Dom^s_{\min}\binom{A}{T}$, we pick
$u_j \in x^{m/2}H_b^{s+m}(\overline{M};E)$ with $Ts_j = Tu_j$ and obtain in this way
a basis
$$
(s_j-u_j) + \Dom^s_{\min}(A_T), \quad j = 1,\ldots,M,
$$
of $\Dom^s_{\max}(A_T)/\Dom^s_{\min}(A_T)$.
\end{lemma}

As already mentioned, the $s_j$ in Lemma \ref{latticeisomanifold} can be chosen to be singular functions (see also
Section \ref{sec-AssociatedDomains}), and the domains $\Dom^s\binom{A}{T}$ as well as the corresponding domains
$\Dom^s(A_T)$ are thus characterized in terms of a specified asymptotic behavior near $\overline{Y}$.

\begin{proposition}\label{Functionalanalyticclosed}
If $A_T - \lambda : \Dom^s(A_T) \to x^{-m/2}H_b^s(\overline{M};E)$ is invertible for some $\lambda \in \C$
and some domain $\Dom^s(A_T)$, then $A$ is closed in the functional analytic sense for every
domain $\Dom^s_{\min}(A_T) \subset \Dom^s \subset \Dom^s_{\max}(A_T)$, i.e.
$\Dom^s_{\max}(A_T)$ is complete in the graph norm
$$
\|u\|_A = \|u\|_{x^{-m/2}H_b^s} + \|Au\|_{x^{-m/2}H_b^s}.
$$
\end{proposition}
\begin{proof}
Let $(u_k)_k \subset \Dom^s_{\max}(A_T)$ be such that $u_k \to u$ in $x^{-m/2}H_b^s$ and
$Au_k \to v$ in $x^{-m/2}H_b^s$. Consequently, $(A-\lambda)u_k \to v-\lambda u$ in
$x^{-m/2}H_b^s$, and by the closed graph theorem the inverse
$$
(A-\lambda)^{-1} : x^{-m/2}H_b^s \to \Dom^s(A_T)
$$
is continuous, where $\Dom^s(A_T)$ is endowed with \eqref{NormwithBoundCond}. Thus there exists
a convergent sequence $(v_k)_k \subset \Dom^s(A_T) \subset \Dom^s_{\max}(A_T)$ with
$(A-\lambda)v_k = (A-\lambda)u_k \to v-\lambda u$ as $k \to \infty$. Let $\ker(A-\lambda) \subset \Dom^s_{\max}(A_T)$
be the eigenspace of $A$ with domain $\Dom^s_{\max}(A_T)$ associated with the eigenvalue $\lambda$. As
$A-\lambda : \Dom^s_{\max}(A_T) \to x^{-m/2}H_b^s$ is Fredholm, this eigenspace is finite dimensional, and so
the norm \eqref{NormwithBoundCond} and the $x^{-m/2}H_b^s$-norm are equivalent on this space. Consequently,
the sequence $(u_k-v_k)_k \subset \ker(A-\lambda)$ is convergent with respect to \eqref{NormwithBoundCond},
and thus $u_k = v_k + (u_k-v_k)$ is also convergent in $\Dom^s_{\max}(A_T)$ with respect to \eqref{NormwithBoundCond},
and the limit necessarily coincides with the $x^{-m/2}H_b^s$-limit $u$.
\end{proof}


\section{The associated boundary value problem on the model cone}
\label{sec-Modeloperator}

For convenience, we fix from now on a collar neighborhood $U_{\overline{Y}} = \overline{Y}\times[0,1)$ of
$\overline{Y}$ in $\overline{M}$. Let $x$ be such that in this neighborhood it coincides
with the projection to $[0,1)$, and the $b$-density $\m$ be such that its pull-back
equals $dy\otimes\frac{dx}{x}$. In this neighborhood, the vector bundles $E$ and $F_j$,
$j=1,\ldots,K$, are isomorphic to the pull-backs of their restrictions to $\overline{Y}$,
and we also fix such isomorphisms.

Every cone operator $B \in x^{-m}\Diff_b^{m}(\overline{M};E,F)$ can be written in the form
\begin{equation}\label{TaylorExpansion}
B = x^{-m}\sum\limits_{k=0}^{N-1}P_kx^k + x^{N-m}\tilde{P}_N,
\end{equation}
where $N \in \N$ is arbitrary, $\tilde{P}_N \in \Diff_b^m(\overline{M};E,F)$, and the
$P_k \in \Diff_b^{m}(\overline{M};E,F)$ have coefficients independent of $x$ near $\overline{Y}$.

Recall that an operator $P \in \Diff_b^{m}(\overline{M};E,F)$ is said to have coefficients independent
of $x$ near $\overline{Y}$, or simply constant coefficients, if
$$
\nabla_{x\partial_x}P(u) = P(\nabla_{x\partial_x}u)
$$
for any smooth section $u$ of $E$ supported in $U_{\overline{Y}}$. Here $\nabla$ denotes a
Hermitian connection on $E$ or $F$, respectively.

Let $\overline{Y}^{\wedge} = \overline{\R}_+\times\overline{Y}$ be the model cone, and
correspondingly let $(\partial \overline{Y})^{\wedge} = \overline{\R}_+\times(\partial \overline{Y})$
be the model cone associated with the boundary.

With $B \in x^{-m}\Diff_b^m(\overline{M};E,F)$ we associate on $\overline{Y}^{\wedge}$ the model operator
$B_{\wedge} = x^{-m}P_0$, where $P_0$ is the constant term in the expansion \eqref{TaylorExpansion}.
Moreover, if $B$ is part of a boundary condition, we let $\gamma_{\wedge}B_{\wedge}$ be the model
boundary condition associated with $\gamma B$, where for every sufficiently smooth section $u$
on $\overline{Y}^{\wedge} \setminus \overline{Y}$ we denote by $\gamma_{\wedge}u$ its restriction to
the regular part of the boundary $(\partial \overline{Y})^{\wedge} \setminus \partial\overline{Y}$.

Consequently, for the family of boundary value problems \eqref{BVPCone} for the operator $A$ there is
the following associated family of model boundary value problems
\begin{equation}\label{BVPModelCone}
\begin{aligned}
(A_{\wedge} - \lambda)u &= f \quad \text{in } \open{\overline{Y}}^{\wedge} = \overline{Y}^{\wedge}\minus\overline{Y}, \\
T_{\wedge}u &= g \quad \text{on } \partial\open{\overline{Y}}^{\wedge} = (\partial \overline{Y})^{\wedge} \minus \partial\overline{Y}
\end{aligned}
\end{equation}
for $A_{\wedge}$ on the model cone $\overline{Y}^{\wedge}$, where
$T_{\wedge} = (\gamma_{\wedge} B_{1,\wedge},\ldots,\gamma_{\wedge} B_{K,\wedge})^t$.

The problem \eqref{BVPModelCone} is naturally realized in the scale of $\K^{s,\alpha}$-spaces on the
model cone. We briefly recall the definition of these spaces (see, e.g., \cite{KaSchu03}):

\begin{definition}\label{ksspaces}
Let ${\mathbb D} \subset S^{n-1}$ be an embedded $(n-1)$-dimensional ball (with boundary). Let
$H^s_{\textup{cone}}(\overline{Y}^{\wedge};E)$ be the space of $H^s_{\textup{loc}}$-distributions $u$ such that given
any coordinate patch $\Omega$ on $\overline{Y}$ diffeomorphic to an open subset of ${\mathbb D} \subset
S^{n-1}$, and given any function $\varphi \in C_0^{\infty}(\Omega)$, we have
$(1-\omega)\varphi u \in H^s({\mathbb D}^{\wedge};E)$, where ${\mathbb D}^{\wedge} = 
\R_+\times{\mathbb D} \subset \R^n$ is regarded as the cone in $\R^{n}\setminus\{0\}$ over
${\mathbb D}$ in polar coordinates, and the Sobolev space on ${\mathbb D}^{\wedge}$ is the space
of $H^s$-distributions in $\R^n$ restricted to that cone.

Correspondingly, we have the space $H^s_{\textup{cone}}((\partial \overline{Y})^{\wedge};F)$ that
is defined in the same way via regarding $(\partial \overline{Y})^{\wedge}$ (locally) as a cone
in $\R^{n-1}$.

Here and in the sequel, $\omega \in C_0^{\infty}(\overline{\R}_+)$ denotes a cut-off function near zero, i.e. $\omega$
is supported near the origin with $\omega \equiv 1$ near zero, and we consider $\omega$ a function either
on $U_{\overline{Y}}$ or on $\overline{Y}^{\wedge}$ (or on $(\partial\overline{Y})^{\wedge}$) which depends
only on the variable $x$.

For $s,\alpha \in \R$ we define $\K^{s,\alpha}(\overline{Y}^{\wedge};E)$ as the space of distributions
$u$ such that
\begin{align*}
\omega u &\in x^{\alpha}H^s_b(\overline{Y}^{\wedge};E) \textup{ and } (1-\omega)u \in x^{(n-m)/2}H^s_{\textup{cone}}(\overline{Y}^{\wedge};E), \\
\intertext{and $\K^{s,\alpha}((\partial\overline{Y})^{\wedge};F)$ as the space of all $u$ with}
\omega u &\in x^{\alpha}H^s_b((\partial\overline{Y})^{\wedge};F) \textup{ and } (1-\omega)u \in x^{(n-1-m)/2}H^s_{\textup{cone}}((\partial\overline{Y})^{\wedge};F).
\end{align*}
Obviously, the $\K^{s,\alpha}$-spaces have natural Hilbert space structures. Note, in particular, that
$\K^{0,-m/2} = x^{-m/2}L^2_b$, and the $x^{-m/2}L^2_b$-inner product serves as the reference
inner product on the model cone.
\end{definition}

For $s > -\frac{1}{2}$ the model boundary value problem \eqref{BVPModelCone} is considered as
{\small
\begin{equation}\label{BVPModelConeinSpaces}
\binom{A_{\wedge}-\lambda}{T_{\wedge}} : \Dom_{\wedge}^s\binom{A_{\wedge}}{T_{\wedge}} \subset
\K^{s,-m/2}(\overline{Y}^{\wedge};E) \to
\begin{array}{c}
\K^{s,-m/2}(\overline{Y}^{\wedge};E) \\
\oplus \\
\bigoplus\limits_{j=1}^K\K^{s+m-m_j-1/2,m/2-m_j}((\partial\overline{Y})^{\wedge};F_j)
\end{array}
\end{equation}}
with $\Dom^s_{\wedge,\min}\binom{A_{\wedge}}{T_{\wedge}} \subset \Dom^s_{\wedge}\binom{A_{\wedge}}{T_{\wedge}}
\subset \Dom^s_{\wedge,\max}\binom{A_{\wedge}}{T_{\wedge}}$, where
\begin{align*}
\Dom_{\wedge,\max}^s\binom{A_{\wedge}}{T_{\wedge}} &= \Bigl\{u \in \K^{s+m,-m/2}(\overline{Y}^{\wedge};E) \st
A_{\wedge}u \in \K^{s,-m/2}(\overline{Y}^{\wedge};E), \\
&\qquad \gamma_{\wedge} B_{j,\wedge}u \in \K^{s+m-m_j-1/2,m/2-m_j}((\partial\overline{Y})^{\wedge};F_j)
\textup{ for } j=1,\ldots,K\Bigr\}, \\
\Dom_{\wedge,\min}^s\binom{A_{\wedge}}{T_{\wedge}} &= \Dom_{\wedge,\max}^s\binom{A_{\wedge}}{T_{\wedge}} \cap \bigcap\limits_{\eps > 0}
\K^{s+m,m/2-\eps}(\overline{Y}^{\wedge};E).
\end{align*}

Analogously to Proposition \ref{Domainproperties} we have:

\begin{proposition}\label{Domainpropertieswedge}
Let \eqref{BVPCone} be $c$-elliptic with parameter in the closed sector $\Lambda$.
\begin{enumerate}[i)]
\item $\Dom^s_{\wedge,\max}\binom{A_{\wedge}}{T_{\wedge}}$ is complete in the norm
$$
\|u\|=\|u\|_{\K^{s+m,-m/2}} + \|A_{\wedge}u\|_{\K^{s,-m/2}} + \sum\limits_{j=1}^K\|\gamma_{\wedge}B_{j,\wedge}u\|_{\K^{s+m-m_j-1/2,m/2-m_j}},
$$
and $\Dom^s_{\wedge,\min}\binom{A_{\wedge}}{T_{\wedge}}$ is a closed subspace of finite codimension.
\item We have
$$
\K^{s+m,m/2} \embed \Dom_{\wedge,\min}^s \embed \Dom_{\wedge,\max}^s \embed \K^{s+m,-m/2+\eps}
$$
for some $\eps > 0$ with continuous embeddings.
\end{enumerate}
\end{proposition}

The quotient $\Dom^s_{\wedge,\max}/\Dom^s_{\wedge,\min}$ is actually independent of $s > -\frac{1}{2}$
and can be described in terms of singular functions as in the boundaryless case, see Section \ref{sec-AssociatedDomains}.

\begin{notation}
For functions $\varphi,\psi$ we write $\varphi \prec \psi$ if $\psi \equiv 1$ in a neighborhood of
the support of $\varphi$.
\end{notation}

\begin{lemma}\label{Twedgesurjective}
Let \eqref{BVPCone} be $c$-elliptic with parameter in $\Lambda$. Then the model boundary condition
$$
T_{\wedge} : \K^{s+m,m/2}(\overline{Y}^{\wedge};E) \to \bigoplus\limits_{j=1}^K\K^{s+m-m_j-1/2,m/2-m_j}((\partial\overline{Y})^{\wedge};F_j)
$$
is surjective for every $s > -\frac{1}{2}$, and necessarily so is its extension to
$\Dom_{\wedge}^s\binom{A_{\wedge}}{T_{\wedge}}$ by Proposition \ref{Domainpropertieswedge}.
\end{lemma}
\begin{proof}
We consider the $b$-elliptic boundary value problem
\begin{equation}\label{AuxTotBVP}
\begin{aligned}
\bigl((x^mA)-\lambda\bigr)u &= f \quad \textup{in } \open{\overline{M}}, \\
\gamma(x^{m_j}B_j)u &= g_j \quad \textup{on } \open{\overline{N}}, \quad j=1,\ldots,K.
\end{aligned}
\end{equation}
Let $h(\sigma,\lambda)$ be the conormal symbol of \eqref{AuxTotBVP}. Thus $h(\sigma,\lambda)$ is for
every $\sigma \in \C$ and every $\lambda \in \Lambda$ an elliptic boundary value problem on $\overline{Y}$,
and the $b$-ellipticity with parameter $\lambda \in \Lambda$ of \eqref{AuxTotBVP} implies that
$h(\sigma,\lambda)$ is elliptic with parameter $(\sigma,\lambda) \in \{\Im(\sigma)=-m/2\}\times\Lambda$.
Consequently, $h(\sigma,\lambda)$ has a parameter-dependent parametrix in Boutet de Monvel's calculus
on $\overline{Y}$ which is an inverse of $h(\sigma,\lambda)$ for all $\Im(\sigma) = -m/2$ and $|\lambda| > R$
sufficiently large.

Let $k(\sigma,\lambda)$ be the row matrix of potential operators in this parametrix, and define
\begin{align*}
K_1(\lambda) &: C_0^{\infty}\Bigl(\R_+\times(\partial\overline{Y}),\bigoplus\limits_{j=1}^KF_j\Bigr) \to
C^{\infty}(\R_+\times\overline{Y};E), \\
\bigl(K_1(\lambda)u\bigr)(x) &:= \frac{1}{2\pi}\int\limits_{\Im(\sigma)=-m/2}\int\limits_0^{\infty}\Bigl(\frac{x}{x'}\Bigr)^{i\sigma}
k(\sigma,\lambda) \begin{pmatrix}x'^{m_1} & \cdots & 0 \\ \vdots & \ddots & \vdots \\ 0 & \cdots & x'^{m_K}\end{pmatrix} u(x')\,\frac{dx'}{x'}\,d\sigma
\end{align*}
for $u \in C_0^{\infty}\Bigl(\R_+\times(\partial\overline{Y}),\bigoplus\limits_{j=1}^KF_j\Bigr) \cong
C_0^{\infty}\Bigl(\R_+,\bigoplus\limits_{j=1}^KC^{\infty}(\partial\overline{Y};F_j)\Bigr)$.

We now have $T_{\wedge}K_1(\lambda) = I$ for $|\lambda| > R$, and by continuity this identity holds on
$\bigoplus\limits_{j=1}^Kx^{m/2-m_j}H_b^{s+m-m_j-1/2}((\partial\overline{Y})^{\wedge};F_j)$. For
cut-off functions $\omega \prec \tilde{\omega}$ near zero we thus have $T_{\wedge}\tilde{\omega}K_1(\lambda)\omega = \omega + R_1(\lambda)$
on $\bigoplus\limits_{j=1}^K\K^{s+m-m_j-1/2,m/2-m_j}$ with a term $R_1(\lambda)$ which decreases
rapidly in the norm as $|\lambda| \to \infty$.

On the other hand, the $c$-ellipticity with parameter $\lambda \in \Lambda$ of \eqref{BVPCone} implies
that the boundary value problem \eqref{BVPModelCone} is away from $x=0$ elliptic with
parameter $\lambda \in \Lambda$ on the cone $\overline{Y}^{\wedge}$ as $x \to \infty$.
Consequently, there exists a parameter-dependent parametrix of $\binom{A_{\wedge}-\lambda}{T_{\wedge}}$
in the SG-calculus of boundary value problems (near infinity), and for the row matrix $K_2(\lambda)$ of
potential operators of this parametrix and a suitable cut-off function $\hat{\omega} \prec \omega$ we have
$T_{\wedge}(1-\hat{\omega})K_2(\lambda)(1-\omega) = (1-\omega) + R_2(\lambda)$ on
$\bigoplus\limits_{j=1}^K\K^{s+m-m_j-1/2,m/2-m_j}$ with a term $R_2(\lambda)$ which decreases
rapidly in the norm as $|\lambda| \to \infty$ (information on parameter-dependent ellipticity
and parametrices of (classical) boundary value problems on manifolds with conical exits to
infinity can be found, e.g., in \cite{KaSchu03}).

Thus for $K(\lambda):= \tilde{\omega}K_1(\lambda)\omega + (1-\hat{\omega})K_2(\lambda)(1-\omega)$ we
have $T_{\wedge}K(\lambda) = I + \tilde{R}(\lambda)$ on $\bigoplus\limits_{j=1}^K\K^{s+m-m_j-1/2,m/2-m_j}$,
and for $|\lambda| > 0$ sufficiently large $I + \tilde{R}(\lambda)$ is invertible.
\end{proof}

According to Lemma \ref{Twedgesurjective} it makes again sense to associate with the boundary value problem
\eqref{BVPModelConeinSpaces} a corresponding family of unbounded operators
\begin{equation}\label{BVPModelConeUnbounded}
A_{\wedge} - \lambda : \Dom^s_{\wedge}(A_{\wedge,T_{\wedge}}) \subset \K^{s,-m/2}(\overline{Y}^{\wedge};E) \to
\K^{s,-m/2}(\overline{Y}^{\wedge};E),
\end{equation}
where $\Dom^s_{\wedge}(A_{\wedge,T_{\wedge}}) = \Dom^s_{\wedge}\binom{A_{\wedge}}{T_{\wedge}}\cap \ker T_{\wedge}$.

Then \eqref{BVPModelConeinSpaces} is invertible for some $\lambda \in \C$ if and only if \eqref{BVPModelConeUnbounded}
is invertible, i.e. if and only if $\lambda \notin \spec(A_{\wedge,T_{\wedge}})$. The analogue of
Lemma \ref{latticeisomanifold} is true also for the associated problem on the model cone.

\begin{proposition}\label{ModelFredholm}
Let \eqref{BVPCone} be $c$-elliptic with parameter $\lambda \in \Lambda$. Then, for $\lambda \neq 0$,
the operator \eqref{BVPModelConeinSpaces} is Fredholm for every intermediate
domain $\Dom_{\wedge,\min}^s\binom{A_{\wedge}}{T_{\wedge}} \subset \Dom_{\wedge}^s\binom{A_{\wedge}}{T_{\wedge}}
\subset \Dom_{\wedge,\max}^s\binom{A_{\wedge}}{T_{\wedge}}$ with index independent of $s > -\frac{1}{2}$.
More precisely, we have
\begin{align*}
\Ind\binom{A_{\wedge}-\lambda}{T_{\wedge}}_{\Dom^s_{\wedge}} &=
\Ind\binom{A_{\wedge}-\lambda}{T_{\wedge}}_{\Dom^s_{\wedge,\min}} + \dim\Dom^s_{\wedge}/\Dom^s_{\wedge,\min} \\
&= \Ind\binom{A}{T}_{\Dom^s_{\min}} + \dim\Dom^s_{\wedge}/\Dom^s_{\wedge,\min} \\
&= \Ind(A_{T,\min}) + \dim\Dom^s_{\wedge}/\Dom^s_{\wedge,\min},
\end{align*}
and correspondingly the operator $A_{\wedge,T_{\wedge}}-\lambda : \Dom^s_{\wedge}(A_{\wedge,T_{\wedge}}) \to \K^{s,-m/2}(\overline{Y}^{\wedge};E)$
is Fredholm for $\lambda \neq 0$ with the same index
\begin{align*}
\Ind (A_{\wedge,T_{\wedge}}-\lambda)_{\Dom^s_{\wedge}} &=
\Ind (A_{\wedge,T_{\wedge}}-\lambda)_{\Dom^s_{\wedge,\min}} + \dim\Dom^s_{\wedge}/\Dom^s_{\wedge,\min} \\
&= \Ind(A_{T,\min}) + \dim\Dom^s_{\wedge}/\Dom^s_{\wedge,\min}.
\end{align*}
\end{proposition}
\begin{proof}
The Fredholmness follows from the parametrix construction in Section \ref{sec-Parametrix}, the
index formula is then elementary except for the assertion that
$$
\Ind\binom{A_{\wedge}-\lambda}{T_{\wedge}}_{\Dom^s_{\wedge,\min}} = \Ind\binom{A}{T}_{\Dom^s_{\min}}.
$$
Under the assumption that $\binom{A_{\wedge}-\lambda}{T_{\wedge}}$ is injective on $\Dom^s_{\wedge,\min}\binom{A_{\wedge}}{T_{\wedge}}$,
this equality is a by-product of Theorem \ref{MainTheoremParametrix}. However, the general case also
follows by the same methods that lead to Theorem \ref{MainTheoremParametrix} by possibly enlarging
the matrices of additional abstract conditions.
\end{proof}

\begin{proposition}\label{Functionalanalyticclosed2}
If $A_{\wedge} - \lambda : \Dom_{\wedge}^s(A_{\wedge,T_{\wedge}}) \to \K^{s,-m/2}(\overline{Y}^{\wedge};E)$ is invertible for
some $\lambda \in \Lambda$ and some domain $\Dom_{\wedge}^s(A_{\wedge,T_{\wedge}})$, then $A_{\wedge}$ is closed
in the functional analytic sense for every domain $\Dom_{\wedge,\min}^s(A_{\wedge,T_{\wedge}}) \subset \Dom_{\wedge}^s
\subset \Dom_{\wedge,\max}^s(A_{\wedge,T_{\wedge}})$, i.e. $\Dom_{\wedge,\max}^s(A_{\wedge,T_{\wedge}})$
is complete in the graph norm
$$
\|u\|_{A_{\wedge}} = \|u\|_{\K^{s,-m/2}} + \|A_{\wedge}u\|_{\K^{s,-m/2}}.
$$
\end{proposition}
\begin{proof}
The proof follows along the lines of the proof of Proposition \ref{Functionalanalyticclosed}, noting
that without loss of generality we may assume $\lambda \neq 0$, and hence
$$
A_{\wedge} - \lambda : \Dom_{\wedge,\max}^s(A_{\wedge,T_{\wedge}}) \to \K^{s,-m/2}(\overline{Y}^{\wedge};E)
$$
is Fredholm according to Proposition \ref{ModelFredholm}.
\end{proof}

\begin{definition}\label{kappagroup}
\begin{enumerate}[i)]
\item For $\varrho > 0$ we define the normalized dilation group action for sections on $\overline{Y}^{\wedge}$
and $(\partial\overline{Y})^{\wedge}$ via
$$
\bigl(\kappa_{\varrho}u\bigr)(x,y) = \varrho^{m/2}u(\varrho x,y).
$$
$\kappa_{\varrho}$ is a strongly continuous group action on the $\K^{s,\alpha}$-spaces, and
the normalization factor $\varrho^{m/2}$ makes it an isometry on $x^{-m/2}L^2_b$.
\item A family of operators $A(\lambda)$ defined on a $\kappa$-invariant space of distributions
on the model cone is called $\kappa$-homogeneous of degree $\nu$ if
$$
A(\varrho^m\lambda) = \varrho^{\nu}\kappa_{\varrho}A(\lambda)\kappa_{\varrho}^{-1}
$$
for every $\varrho > 0$.
\end{enumerate}
It is known that the dilation group action and the notion of $\kappa$-homogeneity play an important
role when dealing with parameter-dependent cone operators, and they are systematically employed
in Schulze's edge pseudodifferential calculus.
\end{definition}

Observe that $\Dom_{\wedge,\min}^s\binom{A_{\wedge}}{T_{\wedge}}$ and $\Dom_{\wedge,\max}^s\binom{A_{\wedge}}{T_{\wedge}}$
as well as the associated domains $\Dom_{\wedge,\min}^s(A_{\wedge,T_{\wedge}})$ and $\Dom_{\wedge,\max}^s(A_{\wedge,T_{\wedge}})$
of the unbounded operator $A_{\wedge}$ under the boundary condition $T_{\wedge}u = 0$
are $\kappa$-invariant. This follows immediately from the $\kappa$-homogeneity
\begin{equation}\label{AwedgeTwedgekappahomogeneous}
\binom{A_{\wedge}}{T_{\wedge}} = \left(\!\begin{array}{c|ccc}
\varrho^m & 0 & \cdots & 0 \\ \hline
\Vsp 0 & \varrho^{m_1} & \cdots & 0 \\
\vdots & \vdots & \ddots & \vdots \\
0 & 0 & \cdots & \varrho^{m_K}
\end{array}\!\right)
\kappa_{\varrho}\binom{A_{\wedge}}{T_{\wedge}}\kappa_{\varrho}^{-1}, \quad \varrho > 0,
\end{equation}
of $\binom{A_{\wedge}}{T_{\wedge}}$. Moreover, this $\kappa$-homogeneity makes it possible to get a fairly
complete picture of what it means to be a sector of minimal growth for realizations of the
operator $A_{\wedge}$ under the boundary condition $T_{\wedge}u = 0$ as the following
Proposition \ref{WedgeSectorMinimalGrowth} shows. Note that the case of $\kappa$-invariant
domains is particularly simple.
In view of the characterization of the domains in terms of singular functions given in
Section \ref{sec-AssociatedDomains}, the estimate \eqref{kappamatrixestimatequotient} below
can be regarded as a condition about the asymptotics of solutions of \eqref{BVPModelCone} as
$|\lambda| \to \infty$.

\begin{proposition}\label{WedgeSectorMinimalGrowth}
Let \eqref{BVPCone} be $c$-elliptic with parameter $\lambda \in \Lambda$. Then the following are
equivalent:
\begin{enumerate}[i)]
\item $\Lambda$ is a sector of minimal growth for the operator $A_{\wedge}$ with domain
$\Dom_{\wedge}(A_{\wedge,T_{\wedge}}) \subset x^{-m/2}L^2_b$.
\item $A_{\wedge} - \lambda : \Dom_{\wedge}(A_{\wedge,T_{\wedge}}) \to x^{-m/2}L_b^2$ is invertible
for large $\lambda \in \Lambda$, and the inverse satisfies the estimate
$$
\|\kappa_{|\lambda|^{1/m}}^{-1}(A_{\wedge}-\lambda)^{-1}\|_{\L(x^{-m/2}L^2_b,\Dom_{\wedge,\max})} = O(|\lambda|^{-1})
$$
as $|\lambda| \to \infty$.
\item
\begin{equation}\label{AwedgemitTwedge}
\binom{A_{\wedge}-\lambda}{T_{\wedge}} : \Dom_{\wedge}\binom{A_{\wedge}}{T_{\wedge}} \to
\begin{array}{c} x^{-m/2}L^2_b \\ \oplus \\ \bigoplus\limits_{j=1}^K\K^{m-m_j-1/2,m/2-m_j} \end{array}
\end{equation}
is invertible for large $\lambda \in \Lambda$, and
\begin{equation}\label{kappamatrixestimate}
\|\kappa^{-1}_{|\lambda|^{1/m}}\binom{A_{\wedge}-\lambda}{T_{\wedge}}^{-1}\kappa_{|\lambda|^{1/m}}\| =
O\begin{pmatrix} |\lambda|^{-1} & |\lambda|^{-m_1/m} & \cdots & |\lambda|^{-m_K/m}\end{pmatrix}
\end{equation}
as $|\lambda| \to \infty$, where the bounds are to be understood componentwise (with values
in $\Dom_{\wedge,\max}\binom{A_{\wedge}}{T_{\wedge}}$).
\item \eqref{AwedgemitTwedge} is bijective for large $\lambda \in \Lambda$, and
\begin{equation}\label{kappamatrixestimatequotient}
\|\kappa^{-1}_{|\lambda|^{1/m}}q_{\wedge}\binom{A_{\wedge}-\lambda}{T_{\wedge}}^{-1}\kappa_{|\lambda|^{1/m}}\| =
O\begin{pmatrix} |\lambda|^{-1} & |\lambda|^{-m_1/m} & \cdots & |\lambda|^{-m_K/m}\end{pmatrix}
\end{equation}
as $|\lambda| \to \infty$, where the bounds are to be understood componentwise with values
in the quotient $\Dom_{\wedge,\max}\binom{A_{\wedge}}{T_{\wedge}}/\Dom_{\wedge,\min}\binom{A_{\wedge}}{T_{\wedge}}$.
Here $q_{\wedge} : \Dom_{\wedge,\max} \to \Dom_{\wedge,\max}/\Dom_{\wedge,\min}$
denotes the canonical projection.

Note that the group action $\kappa_{\varrho}$ descends to the quotient because both $\Dom_{\wedge,\max}$
and $\Dom_{\wedge,\min}$ are $\kappa$-invariant.

\end{enumerate}
If the domain $\Dom_{\wedge}(A_{\wedge,T_{\wedge}})$ is $\kappa$-invariant, i)--iv) are equivalent
to
\begin{enumerate}[v)]
\item[v)] $A_{\wedge} - \lambda : \Dom_{\wedge}(A_{\wedge,T_{\wedge}}) \to x^{-m/2}L^2_b$ is
bijective for all $\lambda \in \Lambda$ with $|\lambda| = 1$.
\end{enumerate}
\end{proposition}
\begin{proof}
\emph{ii) $\Rightarrow$ i)}\ follows immediately because the group action $\kappa_{\varrho}$,
$\varrho > 0$, is unitary on $x^{-m/2}L^2_b$.

\emph{i) $\Rightarrow$ ii)}\/: Note first that $\Dom^s_{\wedge,\max}$ is by assumption complete
in the graph norm, see Proposition \ref{Functionalanalyticclosed2}. Consequently, as $\kappa_{\varrho}$
is an isometry on $x^{-m/2}L_b^2$, we only have to prove that
$$
\|A_{\wedge}\kappa_{|\lambda|^{1/m}}^{-1}(A_{\wedge}-\lambda)^{-1}\|_{\L(x^{-m/2}L^2_b)} = O(|\lambda|^{-1})
$$
as $|\lambda| \to \infty$. From the $\kappa$-homogeneity of $A_{\wedge}$ we obtain
$A_{\wedge}\kappa_{|\lambda|^{1/m}}^{-1} = |\lambda|^{-1}\kappa_{|\lambda|^{1/m}}^{-1}A_{\wedge}$,
and thus the desired estimate follows from the boundedness of the operator family
$A_{\wedge}(A_{\wedge}-\lambda)^{-1}$ in $\L(x^{-m/2}L_b^2)$ (as $|\lambda| \to \infty$), which is
part of our present assumption i).

\emph{ii) $\Rightarrow$ iii)}\/: From the invertibility of $A_{\wedge}-\lambda : \Dom_{\wedge} \to x^{-m/2}L_b^2$
for large $|\lambda| > 0$ and the surjectivity of the boundary condition $T_{\wedge}$ (see
Lemma \ref{Twedgesurjective}) we obtain that \eqref{AwedgemitTwedge} is invertible for large $\lambda$.
Let
$$
\begin{pmatrix} P_{\wedge}(\lambda) & K_{\wedge}(\lambda) \end{pmatrix} = \binom{A_{\wedge}-\lambda}{T_{\wedge}}^{-1} :
\begin{array}{c} x^{-m/2}L^2_b \\ \oplus \\ \bigoplus\limits_{j=1}^K\K^{m-m_j-1/2,m/2-m_j} \end{array}
\to \Dom_{\wedge}\binom{A_{\wedge}}{T_{\wedge}}
$$
be the inverse. The norm estimate in ii) implies the asserted norm estimate for
$P_{\wedge}(\lambda) = (A_{\wedge,T_{\wedge}}-\lambda)^{-1}$ in iii) as $|\lambda| \to \infty$,
noting that $\kappa_{\varrho}$ is an isometry on $x^{-m/2}L_b^2$.

Let $\tilde{K} : \bigoplus\limits_{j=1}^K\K^{m-m_j-1/2,m/2-m_j} \to \K^{m,m/2}(\overline{Y}^{\wedge};E)$
be any right inverse of $T_{\wedge}$ according to Lemma \ref{Twedgesurjective}, and define
$$
\tilde{K}(\varrho):= \kappa_{\varrho}\tilde{K}\kappa_{\varrho}^{-1}
\left(\!\begin{array}{ccc}
\varrho^{-m_1} & \cdots & 0 \\
\vdots & \ddots & \vdots \\
0 & \cdots & \varrho^{-m_K}
\end{array}\!\right)
$$
for $\varrho > 0$. In view of the $\kappa$-homogeneity of $T_{\wedge}$, see
\eqref{AwedgeTwedgekappahomogeneous}, we conclude that $\tilde{K}(\varrho)$ is a right inverse
of $T_{\wedge}$ for every $\varrho > 0$.

From $P_{\wedge}(\lambda)(A_{\wedge}-\lambda)+K_{\wedge}(\lambda)T_{\wedge} = 1$ we get
for large $|\lambda| > 0$
$$
K_{\wedge}(\lambda) = K_{\wedge}(\lambda)T_{\wedge}\tilde{K}(|\lambda|^{1/m}) = 
\tilde{K}(|\lambda|^{1/m})-P_{\wedge}(\lambda)(A_{\wedge}-\lambda)\tilde{K}(|\lambda|^{1/m}).
$$
Using the $\kappa$-homogeneity of $A_{\wedge}-\lambda$ we obtain by conjugation with the group action that
$\kappa^{-1}_{|\lambda|^{1/m}}K_{\wedge}(\lambda)\kappa_{|\lambda|^{1/m}}$ equals
$$
\Bigl(1-
\bigl(\kappa^{-1}_{|\lambda|^{1/m}}P_{\wedge}(\lambda)\kappa_{|\lambda|^{1/m}}\bigr)|\lambda|
\Bigl(A_{\wedge}-\frac{\lambda}{|\lambda|}\Bigr)
\Bigr)\tilde{K}
\left(\!\begin{array}{ccc}
|\lambda|^{-m_1/m} & \cdots & 0 \\
\vdots & \ddots & \vdots \\
0 & \cdots & |\lambda|^{-m_K/m}
\end{array}\!\right),
$$
and thus the asserted norm estimate in iii) holds for $K_{\wedge}(\lambda)$.

\emph{iii) $\Rightarrow$ ii)}\ and \emph{iii) $\Rightarrow$ iv)}\ are immediate.

\emph{iv) $\Rightarrow$ ii)}\/: We just have to worry about the norm estimate.
Let ${\mathcal B}_{T,\wedge}(\lambda) : x^{-m/2}L^2_b \to \Dom_{\wedge,\min}$
be the principal component of the interior part of the parametrix ${\mathcal B}(\lambda)$
from Theorem \ref{MainTheoremParametrix}. Then
$$
1 - {\mathcal B}_{T,\wedge}(\lambda)(A_{\wedge}-\lambda) \equiv 0 \quad \textup{on } \Dom_{\wedge,\min}(A_{\wedge,T_{\wedge}})
$$
for $\lambda \in \Lambda\setminus\{0\}$, and consequently the operator descends to
$$
1 - {\mathcal B}_{T,\wedge}(\lambda)(A_{\wedge}-\lambda) : \Dom_{\wedge,\max}(A_{\wedge,T_{\wedge}})/\Dom_{\wedge,\min}(A_{\wedge,T_{\wedge}}) \to \Dom_{\wedge,\max}(A_{\wedge,T_{\wedge}}).
$$
We may write
$$
\bigl(A_{\wedge,\Dom_{\wedge}}-\lambda\bigr)^{-1} = {\mathcal B}_{T,\wedge}(\lambda) +
\bigl(1-{\mathcal B}_{T,\wedge}(\lambda)(A_{\wedge}-\lambda)\bigr)q_{\wedge}\bigl(A_{\wedge,\Dom_{\wedge}}-\lambda\bigr)^{-1}
$$
as operators $x^{-m/2}L^2_b \to \Dom_{\wedge,\max}(A_{\wedge,T_{\wedge}})$.
By $\kappa$-homogeneity,
\begin{align*}
\kappa^{-1}_{|\lambda|^{1/m}}\bigl(A_{\wedge,\Dom_{\wedge}}-\lambda\bigr)^{-1} &=
|\lambda|^{-1}{\mathcal B}_{T,\wedge}\Bigl(\frac{\lambda}{|\lambda|}\Bigr)\kappa^{-1}_{|\lambda|^{1/m}} \\
&+ \Bigl(1-{\mathcal B}_{T,\wedge}\Bigl(\frac{\lambda}{|\lambda|}\Bigr)\Bigl(A_{\wedge}-\frac{\lambda}{|\lambda|}\Bigr)\Bigr)
\Bigl(\kappa^{-1}_{|\lambda|^{1/m}}q_{\wedge}\bigl(A_{\wedge,\Dom_{\wedge}}-\lambda\bigr)^{-1}\Bigr),
\end{align*}
and so the norm estimate in ii) follows. Recall that the group action $\kappa_{\varrho}$ is
unitary on $x^{-m/2}L^2_b$.

\medskip

If the domain $\Dom_{\wedge}(A_{\wedge,T_{\wedge}})$ is $\kappa$-invariant, then the invertibility
of
\begin{equation}\label{Awedminlamb}
A_{\wedge} - \lambda : \Dom_{\wedge}(A_{\wedge,T_{\wedge}}) \to x^{-m/2}L_b^2
\end{equation}
for large $\lambda \in \Lambda$ is by means of the $\kappa$-homogeneity
$$
A_{\wedge} - \varrho^m\lambda = \varrho^{m}\kappa_{\varrho}\bigl(A_{\wedge} - \lambda\bigr)\kappa_{\varrho}^{-1}
$$
equivalent to the invertiblity of \eqref{Awedminlamb} for all $\lambda \in \Lambda \setminus\{0\}$ or, equivalently,
only for $\lambda \in \Lambda$ with $|\lambda| = 1$. Moreover, from the identity
$$
\kappa_{|\lambda|^{1/m}}^{-1}(A_{\wedge}-\lambda)^{-1}\kappa_{|\lambda|^{1/m}} = |\lambda|^{-1}
\Bigl(A_{\wedge}-\frac{\lambda}{|\lambda|}\Bigr)^{-1} : x^{-m/2}L_b^2 \to \Dom_{\wedge}(A_{\wedge,T_{\wedge}})
$$
for $\lambda \neq 0$ we automatically obtain the norm estimate in ii),
and consequently the equivalence of i)--iv) and v) is proved.
\end{proof}

While all conditions in Proposition \ref{WedgeSectorMinimalGrowth} have an analytic
flavour, it should be noted that one can give a characterization of these properties
in geometric terms that involve the Grassmannian of domains of $A_{\wedge,T_{\wedge}}$
and the flow induced by $\kappa_{\varrho}$ on that Grassmannian. We refer to
\cite{GKM1, GKM3} for details in the boundaryless situation.


\section{Domains, associated domains, and singular functions}
\label{sec-AssociatedDomains}

In this section we give a description of domains of the realizations
of $A$ and $A_{\wedge}$ under the boundary condition $Tu = 0$ and
$T_{\wedge}u = 0$, respectively, in terms of singular functions, i.e.
the domains are characterized by the asymptotic behavior of their
elements near the ``singular boundary'' $\overline{Y}$.

Moreover, we explicitly construct an isomorphism
$$
\theta : \Dom_{\max}/\Dom_{\min} \to \Dom_{\wedge,\max}/\Dom_{\wedge,\min}
$$
that will be used to associate with a domain $\Dom$ of $A$ under the
boundary condition $Tu = 0$ a corresponding domain $\Dom_{\wedge}$ of $A_{\wedge}$
under the boundary condition $T_{\wedge}u = 0$ via
\begin{equation}\label{AssocDomain}
\theta\bigl(\Dom/\Dom_{\min}\bigr) = \Dom_{\wedge}/\Dom_{\wedge,\min}.
\end{equation}
The ellipticity condition for the resolvent constructions for the operator
$A_T$ with domain $\Dom(A_T)$ in Section \ref{sec-Resolvents} then involves a
spectral condition on the model operator $A_{\wedge,T_{\wedge}}$ with the
associated domain $\Dom(A_{\wedge,T_{\wedge}})$ to $\Dom(A_T)$ according to
\eqref{AssocDomain}. As the boundaryless case shows, a condition
of such type is to be expected also in this more advanced situation.

Our approach to consider domains with inhomogeneous boundary conditions
$\Dom\binom{A}{T}$ as well as $\Dom\binom{A_{\wedge}}{T_{\wedge}}$ makes
it possible to transfer the methods from \cite{GKM1} and \cite{GKM2}.

According to \eqref{TaylorExpansion} we write (near $\overline{Y}$)
\begin{alignat}{2}
\label{exp1} A &\equiv x^{-m}\sum\limits_{k=0}^{m-1}A_kx^k & &\mod \Diff_b^m(\overline{M};E), \\
\label{exp2} B_j &\equiv x^{-m_j}\sum\limits_{k=0}^{m-1}B_{j,k}x^k & &\mod x^{m-m_j}\Diff_b^{m_j}(\overline{M};E,F_j),
\quad j=1,\ldots,K,
\end{alignat}
with totally characteristic operators $A_k, B_{j,k}$ with coefficients independent of $x$
near $\overline{Y}$, and therefore they can be regarded also as operators acting in
sections on the model cone $\overline{Y}^{\wedge}$. Thus
\begin{equation}\label{exp3}
\begin{pmatrix} A \\ \gamma B_1 \\ \vdots \\ \gamma B_K \end{pmatrix} =
\begin{pmatrix} x^{-m} & 0 & \cdots & 0 \\ 0 & x^{-m_1} & \cdots & 0 \\ \vdots & \vdots & \ddots & \vdots \\
0 & 0 & \cdots & x^{-m_K} \end{pmatrix}
\sum\limits_{k=0}^{m-1}
\begin{pmatrix} A_k \\ \gamma B_{1,k} \\ \vdots \\ \gamma B_{K,k} \end{pmatrix} x^k + \tilde{R},
\end{equation}
and set ${\mathcal A}_k = \begin{pmatrix} A_k & \gamma B_{1,k} & \cdots & \gamma B_{K,k} \end{pmatrix}^t$,
$k = 0,\ldots,m-1$. Let
\begin{equation}\label{Konormalensymbol}
\hat{{\mathcal A}}_k(\sigma) : H^{s+m}(\overline{Y};E) \to
\begin{array}{c} H^{s}(\overline{Y};E) \\ \oplus \\ \bigoplus\limits_{j=1}^KH^{s+m-m_j-1/2}(\partial\overline{Y};F_j) \end{array}, \quad s > -\frac{1}{2},
\end{equation}
$\sigma \in \C$, be the conormal symbol of ${\mathcal A}_k$, see Section \ref{sec-Totally}.
From our standing assumption that \eqref{BVPCone} is $c$-elliptic (with parameter $\lambda \in \Lambda$),
we obtain that the leading term $\hat{{\mathcal A}}_0(\sigma)$ is a holomorphic Fredholm family in
\eqref{Konormalensymbol} which has a finitely meromorphic inverse $\hat{{\mathcal A}}^{-1}_0(\sigma)$.

\medskip

In the sequel, we make use of the following notion of Mellin transform for sections $u$ on
$\overline{M}$ or $\overline{Y}^{\wedge}$, respectively, which employs apriori a cut-off near $\overline{Y}$:

Fix a cut-off function $\omega \in C_0^{\infty}([0,1))$ near zero, i.e. $\omega$ is real valued and supported
near the origin with $\omega \equiv 1$ near zero. As usual, we regard $\omega$ as a function on
$\overline{M}$ supported in the collar neighborhood $U_{\overline{Y}} \subset \overline{M}$,
or on $\overline{Y}^{\wedge}$. Then the Mellin transform of a section $u \in C_0^{\infty}(\open{\overline{M}};E)$ is defined to be
the entire function $\hat{u} : \C \to C^{\infty}(\overline{Y};E|_{\overline{Y}})$ such that
for any $v \in C^{\infty}(\overline{Y};E|_{\overline{Y}})$
\begin{equation}\label{Mellintrafo}
\bigl(x^{-i\sigma}\omega u,\pi^*_{\overline{Y}}v\bigr)_{L^2_b(\overline{M};E)} =
\bigl(\hat{u}(\sigma),v\bigr)_{L^2(\overline{Y};E|_{\overline{Y}})},
\end{equation}
where $\pi^*_{\overline{Y}}v$ is the section of $E$ over $U_{\overline{Y}}$ obtained by parallel
transport of $v$ along the fibers of the projection $\pi_{\overline{Y}}$.
The Mellin transform of sections $u \in C_0^{\infty}(\open{\overline{Y}}^{\wedge};E)$ is defined in
the same way, but the pairing in \eqref{Mellintrafo} is the inner product in
$L^2_b(\overline{Y}^{\wedge};E)$ (where, as before, we identify the bundle $E \to \overline{Y}^{\wedge}$
with the pull-back $\pi^*_{\overline{Y}}E|_{\overline{Y}}$).

The Mellin transform extends to the spaces $x^{\alpha}H^s_b$ and $\K^{s,\alpha}$ in such a way that
$\hat{u}(\sigma)$ is a holomorphic $H^s(\overline{Y};E)$-valued function in $\{\Im(\sigma) > -\alpha\}$
with well known integrability conditions along lines parallel to the real axis.

In the same way we also define the Mellin transform for sections on the boundary $\overline{N}$,
and on the model cone $(\partial\overline{Y})^{\wedge}$ associated with the boundary, respectively.

\medskip

Let
$$
\spec_b\binom{A}{T} = \{\sigma \in \C \st \hat{{\mathcal A}}_0(\sigma) \textup{ is not invertible}\} \subset \C
$$
be the boundary spectrum of $\binom{A}{T}$. Then
$$
\spec_b\binom{A}{T}\cap \{\sigma \in \C \st \alpha < \Im(\sigma) < \beta\}
$$
is finite for all $\alpha,\beta \in \R$, $\alpha < \beta$, and let
\begin{equation}\label{SigmaDefinition}
\Sigma:= \spec_b\binom{A}{T} \cap \{\sigma \in \C \st -m/2 < \Im(\sigma) < m/2\}
\end{equation}
be the part of the boundary spectrum in the critical strip that is associated with realizations of $A$ and $A_{\wedge}$
in $x^{-m/2}L_b^2$ under the boundary condition $Tu = 0$ and $T_{\wedge}u = 0$, respectively.

For $\sigma_0 \in \Sigma$ let $\tilde\Sing_{\wedge,\sigma_0}$ be the space of all singular
functions of the form
$$
q = \Bigl(\sum\limits_{k=0}^{m_{\sigma_0}}c_{\sigma_0,k}(y)\log^kx\Bigr)x^{i\sigma_0} \in
C^{\infty}(\open{\overline{Y}}^{\wedge};E),
$$
where $c_{\sigma_0,k} \in C^{\infty}(\overline{Y};E)$ and $m_{\sigma_0} \in \N_0$, such that
${\mathcal A}_0q = 0$. Using the Mellin transform, this is equivalent to the holomorphicity
of $\hat{{\mathcal A}}_0(\sigma)\hat{q}(\sigma)$ on the whole complex plane, and as the inverse
$\hat{{\mathcal A}}_0^{-1}(\sigma)$ is finitely meromorphic (with regularizing principal parts of
Laurent expansions) we see that the space $\tilde\Sing_{\wedge,\sigma_0}$ is finite dimensional.

We set
$$
\tilde\Sing_{\wedge,\max} = \bigoplus\limits_{\sigma_0 \in \Sigma}\tilde\Sing_{\wedge,\sigma_0} \subset
C^{\infty}(\open{\overline{Y}}^{\wedge};E).
$$

Let $u \in \Dom^s_{\wedge,\max}\binom{A_{\wedge}}{T_{\wedge}}$. By Mellin transform and the definition
of the maximal domain, we thus obtain that $\hat{{\mathcal A}}_0(\sigma)\hat{u}(\sigma)$ is the Mellin
transform of a vector of functions
$$
v \in \begin{array}{c}
\K^{s,m/2}(\overline{Y}^{\wedge};E) \\
\oplus \\
\bigoplus\limits_{j=1}^K\K^{s+m-m_j-1/2,m/2}((\partial\overline{Y})^{\wedge};F_j)
\end{array}.
$$
In particular, $\hat{{\mathcal A}}_0(\sigma)\hat{u}(\sigma)$ is holomorphic in $\{\Im(\sigma) > -m/2\}$,
and by the meromorphic structure of $\hat{{\mathcal A}}_0^{-1}(\sigma)$ we see that there is a
singular function $q \in \tilde\Sing_{\wedge,\max}$ such that $\hat{u}(\sigma) - \hat{q}(\sigma)$ is
holomorphic in the critical strip $\{\sigma \in \C \st -m/2 < \Im(\sigma) < m/2\}$. Consequently,
$u - \omega q \in \Dom^s_{\wedge,\max}\binom{A_{\wedge}}{T_{\wedge}}$ with holomorphic Mellin
transform, and thus $u - \omega q \in \Dom^s_{\wedge,\min}\binom{A_{\wedge}}{T_{\wedge}}$. Note
that the minimal domain as a subspace of the maximal domain is characterized by the property that
the Mellin transforms of its elements are holomorphic in the critical strip.

Let us summarize this in the following proposition:

\begin{proposition}\label{SingularFunctionModel}
Every class
$$
u + \Dom^s_{\wedge,\min}\binom{A_{\wedge}}{T_{\wedge}} \in
\Dom^s_{\wedge,\max}\binom{A_{\wedge}}{T_{\wedge}}/\Dom^s_{\wedge,\min}\binom{A_{\wedge}}{T_{\wedge}}
$$
contains a representative of the form $\omega q$ with $q \in \tilde\Sing_{\wedge,\max}$, and
the singular function $q$ is uniquely determined by its class modulo $\Dom^s_{\wedge,\min}$.

In this way we obtain an isomorphism
$$
\Dom^s_{\wedge,\max}\binom{A_{\wedge}}{T_{\wedge}}/\Dom^s_{\wedge,\min}\binom{A_{\wedge}}{T_{\wedge}} \cong
\tilde\Sing_{\wedge,\max},
$$
and the quotient $\Dom^s_{\wedge,\max}/\Dom^s_{\wedge,\min}$ is independent of $s > -\frac{1}{2}$.

Consequently, specifying a domain $\Dom_{\wedge,\min}^s \subset \Dom_{\wedge}^s \subset
\Dom_{\wedge,\max}^s$ is equivalent to specifying a subspace of $\tilde\Sing_{\wedge,\max}$ of
admissible conormal asymptotics for the elements $u \in \Dom_{\wedge}^s$ near $\overline{Y}$.
\end{proposition}

In view of
$$
\Dom^s_{\wedge,\max}\binom{A_{\wedge}}{T_{\wedge}}/\Dom^s_{\wedge,\min}\binom{A_{\wedge}}{T_{\wedge}} \cong
\Dom^s_{\wedge,\max}(A_{\wedge,T_{\wedge}})/\Dom^s_{\wedge,\min}(A_{\wedge,T_{\wedge}}),
$$
see also Lemma \ref{latticeisomanifold}, we also obtain
$$
\Dom^s_{\wedge,\max}(A_{\wedge,T_{\wedge}})/\Dom^s_{\wedge,\min}(A_{\wedge,T_{\wedge}}) \cong
\tilde\Sing_{\wedge,\max},
$$
and the domains of the unbounded operator $A_{\wedge}$ under the boundary condition $T_{\wedge}u = 0$
are characterized in terms of the asymptotics near $\overline{Y}$.

\medskip

Now let $u \in \Dom^s_{\max}\binom{A}{T}$. Then we obtain analogously to the case of the
model cone that $\sum\limits_{k=0}^{m-1}\hat{\mathcal A}_k(\sigma)\hat{u}(\sigma+ik)$ is the
Mellin transform of a vector of functions
$$
v \in \begin{array}{c}
x^{m/2}H_b^{s}(\overline{M};E) \\
\oplus \\
\bigoplus\limits_{j=1}^Kx^{m/2}H_b^{s+m-m_j-1/2}(\overline{N};F_j)
\end{array},
$$
and consequently is holomorphic in $\{\Im(\sigma) > -m/2\}$. By inductively arguing for the strips
$\{m/2-k < \Im(\sigma) < m/2\}$, $k = 1,\ldots,m$, using thereby the meromorphic structure of the
inverse $\hat{{\mathcal A}}_0^{-1}(\sigma)$ and the apriori holomorphicity of $\hat{u}(\sigma)$ in
$\{\Im(\sigma) > m/2\}$, we conclude that $\hat{u}(\sigma)$ has a meromorphic extension to the critical strip
$\{-m/2 < \Im(\sigma) < m/2\}$, and there exists a singular function of the form
\begin{equation}\label{Singulaerfunktion}
q = \!\!\sum_{-\frac{m}{2} < \Im(\sigma) < \frac{m}{2}}
\Bigl(\sum_{k=0}^{m_{\sigma}}c_{\sigma,k}(y)\log^kx\Bigr)x^{i\sigma}
\end{equation}
with $c_{\sigma,k} \in C^{\infty}(\overline{Y};E)$, $m_\sigma \in \N_0$,
such that $\hat{u}(\sigma) - \hat{q}(\sigma)$ is holomorphic in this strip. Note that the
sum in \eqref{Singulaerfunktion} is actually only a finite sum. Consequently, as also
$\omega q \in \Dom^s_{\max}\binom{A}{T}$, we conclude that $u - \omega q \in \Dom^s_{\min}\binom{A}{T}$.
We hence obtain an isomorphism
$$
\Dom^s_{\max}\binom{A}{T}/\Dom^s_{\min}\binom{A}{T} \cong \tilde\Sing_{\max}
$$
to a finite dimensional space of singular functions $\tilde\Sing_{\max} \subset C^{\infty}(\open{\overline{Y}}^{\wedge};E)$
similar to the case of the model operator in Proposition \ref{SingularFunctionModel}.

\medskip

Let us be more precise about the structure of the space $\tilde\Sing_{\max}$ of singular functions:
We may write
$$
\tilde\Sing_{\max} = \bigoplus\limits_{\sigma_0 \in \Sigma}\tilde\Sing_{\sigma_0},
$$
and the elements $q \in \tilde\Sing_{\sigma_0}$ are of the form
$$
q = \sum\limits_{\vartheta=0}^{N(\sigma_0)}\Bigl(\sum\limits_{k=0}^{m_{\sigma_0-i\vartheta}}
c_{\sigma_0-i\vartheta,k}(y)\log^kx\Bigr)x^{i(\sigma_0-i\vartheta)}
$$
with $c_{\sigma_0-i\vartheta,k} \in C^{\infty}(\overline{Y};E)$, $m_{\sigma_0-i\vartheta} \in \N_0$,
and $N(\sigma_0) \in \N_0$ the largest integer such that $\Im(\sigma_0)-N(\sigma_0) > -m/2$.

More precisely, there is an isomorphism
$$
\theta : \tilde\Sing_{\max} \to \tilde\Sing_{\wedge,\max}
$$
that was already mentioned in the introduction of this section, which restricts to isomorphisms
$\theta|_{\tilde\Sing_{\sigma_0}} : \tilde\Sing_{\sigma_0} \to \tilde\Sing_{\wedge,\sigma_0}$.
The inverse $\theta^{-1}|_{\tilde\Sing_{\sigma_0}}$ is of the form
$$
\theta^{-1}|_{\tilde\Sing_{\sigma_0}} = \sum\limits_{k=0}^{N(\sigma_0)}\e_{\sigma_0,k} :
\tilde\Sing_{\wedge,\sigma_0} \to \tilde\Sing_{\sigma_0},
$$
where the $\e_{\sigma_0,k} : \tilde\Sing_{\wedge,\sigma_0} \to C^{\infty}(\open{\overline{Y}}^{\wedge};E)$
are inductively defined as follows:
\begin{itemize}
\item $\e_{\sigma_0,0} = \id$, the identity map.
\item Given $\e_{\sigma_0,0},\ldots,\e_{\sigma_0,\vartheta-1}$ for some
$\vartheta \in \{1,\ldots,N(\sigma_0)-1\}$, we define 
$\e_{\sigma_0,\vartheta}(\psi)$ for $\psi \in \tilde\Sing_{\wedge,\sigma_0}$ 
to be the unique singular function of the form
$$
\Bigl(\sum_{k=0}^{m_{\sigma_0-i\vartheta}}c_{\sigma_0-i\vartheta,k}(y)
\log^kx\Bigr)x^{i(\sigma_0-i\vartheta)}
$$
such that
\begin{equation*}
(\e_{\sigma_0,\vartheta}(\psi))^{\wedge}(\sigma) + 
\hat{{\mathcal A}}_0(\sigma)^{-1}\Bigl(\sum_{k=1}^{\vartheta}\hat{{\mathcal A}}_k(\sigma)
\s_{\sigma_0-i\vartheta}(\e_{\sigma_0,\vartheta-k}(\psi))^{\wedge}
(\sigma+ik)\Bigr)
\end{equation*}
is holomorphic at $\sigma = \sigma_0 - i\vartheta$, where 
$(\e_{\sigma_0,\vartheta-k}(\psi))^{\wedge}(\sigma)$ is
the Mellin transform of the function $\e_{\sigma_0,\vartheta-k}(\psi)$, 
and $\s_{\sigma_0-i\vartheta}(\e_{\sigma_0,\vartheta-k}(\psi))^{\wedge}
(\sigma+ik)$ is the singular part of its Laurent expansion at $\sigma_0 - 
i\vartheta$.
Recall that our notion of Mellin transform involves apriori a cut-off near $\overline{Y}$,
and so $(\e_{\sigma_0,\vartheta-k}(\psi))^{\wedge}(\sigma)$
is meromorphic in $\C$ with only one pole at $\sigma_0 - i(\vartheta-k)$.
\end{itemize}

It is of interest to note that this construction yields
$$
\sum_{k=0}^{\vartheta}\bigl({\mathcal A}_kx^k\bigr)(\e_{\sigma_0,\vartheta-k}(\psi))=0
$$
for every $\psi \in \tilde\Sing_{\wedge,\sigma_0}$ and every
$\vartheta = 0,\ldots,N(\sigma_0)$.

In conclusion, every space $\tilde\Sing_{\sigma_0}$ consists indeed of singular
functions of the form
\begin{equation*}
q = \sum_{\vartheta=0}^{N(\sigma_0)}
\Bigl(\sum_{k=0}^{m_{\sigma_0-i\vartheta}}c_{\sigma_0-i\vartheta,k}(y)
\log^kx\Bigr)x^{i(\sigma_0-i\vartheta)},
\end{equation*}
and we have
\begin{equation} \label{ThetaOperator}
\theta q = \Bigl(\sum_{k=0}^{m_{\sigma_0}}c_{\sigma_0,k}(y)
\log^kx\Bigr)x^{i\sigma_0}.
\end{equation}
It is more tedious than hard to verify that this furnishes an isomorphism $\theta : \tilde\Sing_{\max} \to \tilde\Sing_{\wedge,\max}$
as desired (see also \cite{GKM1} for further information in the boundaryless context).

\medskip

Let us summarize the above in the following proposition:

\begin{proposition}\label{SingularFunctionManifold}
\begin{enumerate}[i)]
\item There is a natural isomorphism
$$
\Dom^s_{\max}\binom{A}{T}/\Dom^s_{\min}\binom{A}{T} \cong \tilde\Sing_{\max}, \quad
u + \Dom^s_{\min}\binom{A}{T} \mapsto q,
$$
that is characterized by the property that $u - \omega q \in \Dom^s_{\min}\binom{A}{T}$,
where $\omega \in C_0^{\infty}([0,1))$ is any cut-off function near zero.

Consequently, the quotient $\Dom^s_{\max}/\Dom^s_{\min}$
is independent of $s > -\frac{1}{2}$, and its elements are characterized by their asymptotic behavior
near $\overline{Y}$.
\item By Lemma \ref{latticeisomanifold},
$$
\Dom^s_{\max}\binom{A}{T}/\Dom^s_{\min}\binom{A}{T} \cong \Dom^s_{\max}(A_T)/\Dom^s_{\min}(A_T),
$$
and consequently also the quotient $\Dom^s_{\max}(A_T)/\Dom^s_{\min}(A_T)$ of the maximal and minimal
domains of the unbounded operator $A$ under the boundary condition $Tu = 0$ is characterized by
the conormal asymptotics in $\tilde\Sing_{\max}$.
\item There is a natural isomorphism $\theta : \tilde\Sing_{\max} \to \tilde\Sing_{\wedge,\max}$ that
by i), ii), and Proposition \ref{SingularFunctionModel} gives rise to isomorphisms
$$
\theta : \begin{cases}
\Dom^s_{\max}\binom{A}{T}/\Dom^s_{\min}\binom{A}{T} \to \Dom^s_{\wedge,\max}\binom{A_{\wedge}}{T_{\wedge}}/\Dom^s_{\wedge,\min}\binom{A_{\wedge}}{T_{\wedge}}, \\
\Dom^s_{\max}(A_T)/\Dom^s_{\min}(A_T) \to \Dom^s_{\wedge,\max}(A_{\wedge,T_{\wedge}})/\Dom^s_{\wedge,\min}(A_{\wedge,T_{\wedge}}).
\end{cases}
$$
For a domain $\Dom^s_{\min} \subset \Dom^s \subset \Dom^s_{\max}$ we therefore have an associated domain
$$
\Dom^s_{\wedge,\min} \subset \theta(\Dom^s)= \Dom^s_{\wedge} \subset \Dom^s_{\wedge,\max}
$$
via $\theta\bigl(\Dom^s/\Dom^s_{\min}\bigr) = \Dom^s_{\wedge}/\Dom^s_{\wedge,\min}$.
\end{enumerate}
\end{proposition}


\section{Parametrix construction}
\label{sec-Parametrix}

Let ${\mathcal A}(\lambda)$ denote the boundary value problem \eqref{BVPCone}.
Our goal in this section is the construction of a parametrix under the assumption
that ${\mathcal A}(\lambda)$ is $c$-elliptic with parameter $\lambda \in \Lambda$, and
that the model operator
$$
{\mathcal A}_{\wedge}(\lambda) = \begin{pmatrix} A_{\wedge}-\lambda \\ T_{\wedge} \end{pmatrix} :
\Dom^s_{\wedge,\min}\binom{A_{\wedge}}{T_{\wedge}} \to
\begin{array}{c}
\K^{s,-m/2}(\overline{Y}^{\wedge};E) \\ \oplus \\
\bigoplus\limits_{j=1}^K\K^{s+m-m_j-1/2,m/2-m_j}((\partial\overline{Y})^{\wedge};F_j)
\end{array}
$$
is injective for some $s > -\frac{1}{2}$ and all $\lambda \in \Lambda\setminus\{0\}$.

More precisely, we will construct a parametrix
$$
{\mathcal B}(\lambda) = \begin{pmatrix} {\mathcal B}_T(\lambda) & K(\lambda) \end{pmatrix} :
\begin{array}{c}
x^{-m/2}H_b^{s}(\overline{M};E) \\ \oplus \\
\bigoplus\limits_{j=1}^K x^{m/2-m_j}H_b^{s+m-m_j-1/2}(\overline{N};F_j)
\end{array} \to \Dom^s_{\min}\binom{A}{T}
$$
such that
$$
{\mathcal B}(\lambda){\mathcal A}(\lambda) - 1 : \Dom^s_{\min}\binom{A}{T} \to \Dom^s_{\min}\binom{A}{T}
$$
is regularizing and compactly supported in $\lambda \in \Lambda$. In particular, for $\lambda$
sufficiently large, the boundary value problem
\begin{equation}\label{BVPConemin}
{\mathcal A}(\lambda) : \Dom^s_{\min}\binom{A}{T} \to
\begin{array}{c}
x^{-m/2}H_b^{s}(\overline{M};E) \\ \oplus \\
\bigoplus\limits_{j=1}^K x^{m/2-m_j}H_b^{s+m-m_j-1/2}(\overline{N};F_j)
\end{array}
\end{equation}
is injective and ${\mathcal B}(\lambda)$ is a left inverse. Moreover, the regularizing remainder
$$
\Pi(\lambda) = 1 - {\mathcal A}(\lambda){\mathcal B}(\lambda)
$$
is a finite dimensional projection to a complement of the range of \eqref{BVPConemin}.

For the actual construction of this parametrix we employ some ideas from pseudodifferential operator
theory of Shapiro-Lopatinsky elliptic edge-degenerate boundary value problems, the central topic
of the monograph \cite{KaSchu03}.

\bigskip

Choose local coordinates on $\overline{Y}$ centered at zero, and let $(0,1)\times\Omega$ be corresponding
coordinates in the collar neighborhood $U_{\overline{Y}} \subset \overline{M}$ of $\overline{Y}$.
In these coordinates, the operator $A-\lambda$ takes the form
$$
A-\lambda = x^{-m}\Bigl(\sum\limits_{k+|\alpha|\leq m}a_{k,\alpha}(x,y)D_y^{\alpha}(xD_x)^k-x^m\lambda\Bigr),
$$
and thus its complete symbol $\tilde{a}(x,y,\xi,\eta,\lambda)$ is given by
$$
\tilde{a}(x,y,\xi,\eta,\lambda) = x^{-m}a(x,y,x\xi,\eta,x^m\lambda)
$$
with a symbol $a(x,y,\xi,\eta,\lambda)$ that is smooth up to $x = 0$.
The $c$-ellipticity with parameter $\lambda \in \Lambda$ of $A-\lambda$ is equivalent to the
invertibility of the principal component $a_{(m)}(x,y,\xi,\eta,\lambda)$ for all covectors
$(\xi,\eta,\lambda)$ different from zero, and all $(x,y) \in [0,1)\times\Omega$ (up to $x=0$).
Note that the principal component $a_{(m)}$ is (anisotropic) homogeneous, i.e.
$$
a_{(m)}(x,y,\varrho\xi,\varrho\eta,\varrho^m\lambda) = \varrho^m a_{(m)}(x,y,\xi,\eta,\lambda)
$$
for $\varrho > 0$.

Assume for a moment that $\Omega \subset \R^{n-1}$ corresponds to an interior chart on $\overline{Y}$.
Then the parametrix construction from Section 5 in \cite{GKM2} implies that there exists a symbol
$p(x,y,\xi,\eta,\lambda)$ with the following properties:
\begin{enumerate}[i)]
\item $p$ is smooth in all variables up to $x = 0$.
\item We have
$$
\bigl|\partial^{\alpha}_{(x,y)}\partial^{\beta}_{(\xi,\eta)}\partial^{\gamma}_{\lambda}p(x,y,\xi,\eta,\lambda)\bigr| =
O\bigl((1+|\xi|+|\eta|+|\lambda|^{1/m})^{-m-|\beta|-m|\gamma|}\bigr)
$$
as $|(\xi,\eta,\lambda)| \to \infty$, locally uniformly for $(x,y) \in [0,1)\times\Omega$.
\item $p$ is a classical symbol, i.e. it admits an asymptotic expansion
$$
p \sim \sum\limits_{j=0}^{\infty}\chi(\xi,\eta,\lambda)p_{(-m-j)}(x,y,\xi,\eta,\lambda),
$$
where $\chi \in C^{\infty}(\R\times\R^{n-1}\times\Lambda)$ is a function such that $\chi \equiv 0$
near the origin and $\chi \equiv 1$ for large $|(\xi,\eta,\lambda)|$, and the components
$p_{(-m-j)}$ are anisotropic homogeneous, i.e. we have
$$
p_{(-m-j)}(x,y,\varrho\xi,\varrho\eta,\varrho^m\lambda) = \varrho^{-m-j}p_{(-m-j)}(x,y,\xi,\eta,\lambda)
$$
for $\varrho > 0$.
\item $(A-\lambda)\textup{Op}(x^mp(x,y,x\xi,\eta,x^m\lambda)) - 1$ and
$\textup{Op}(x^mp(x,y,x\xi,\eta,x^m\lambda))(A-\lambda) - 1$ are parameter-dependent smoothing
pseudodifferential operators on $(0,1)\times\Omega$, where
$\textup{Op}(x^mp(x,y,x\xi,\eta,x^m\lambda))$ denotes the standard Kohn-Nirenberg quantized
pseudodifferential operator in $(0,1)\times\Omega$ with symbol $x^mp(x,y,x\xi,\eta,x^m\lambda)$.
\end{enumerate}
Now let $\Omega \subset \overline{\R}_+^{n-1}$ correspond to a boundary chart on $\overline{Y}$.
We slightly extend $\Omega$ as well as $\tilde{a}$ over the boundary $\R^{n-2} \subset \overline{\R}_+^{n-1}$
such that the structure of the complete symbol $\tilde{a}$ of $A-\lambda$ and the $c$-ellipticity
with parameter remains preserved (by possibly shrinking $\Omega$ to a relatively compact chart,
this is always possible).

Now we can use the beforementioned results about the existence of a parametrix in the extended
domain $(0,1)\times\Omega'$, where the symbol $p$ in addition has the transmission property
with respect to the boundary $\R^{n-2} \subset \overline{\R}_+^{n-1}$.

Passing, as is usual in pseudodifferential boundary value problems, to
$$
\textup{Op}^+(x^mp(x,y,x\xi,\eta,x^m\lambda)) =
r^+\textup{Op}(x^mp(x,y,x\xi,\eta,x^m\lambda))e^+,
$$
where $e^+$ denotes the operator of extension by zero from the original domain
$(0,1)\times\Omega$ to the extended domain $(0,1)\times\Omega'$ and $r^+$ denotes restriction,
we obtain that
$$
\textup{Op}^+(x^mp(x,y,x\xi,\eta,x^m\lambda))(A-\lambda) = 1 + G_0(\lambda) + G_{-\infty}(\lambda),
$$
where $G_{-\infty}(\lambda)$ is a parameter-dependent regularizing singular Green operator
in Boutet de Monvel's calculus, and $G_0(\lambda)$ is a parameter-dependent singular Green operator
in Boutet de Monvel's calculus of order zero whose boundary symbol has the form
$$
\tilde{g}(x,y',\xi,\eta',\lambda) = g(x,y',x\xi,\eta',x^m\lambda)
$$
with an (anisotropic) parameter-dependent singular Green symbol $g(x,y',\xi,\eta',\lambda)$ of
order zero. The structure of the composition
$$
(A-\lambda)\textup{Op}^+(x^mp(x,y,x\xi,\eta,x^m\lambda))
$$
is the same.

By combining the standard parametrix constructions on a manifold with boundary away from $\overline{Y}$
with the above considerations, we arrive at the following

\begin{lemma}
$A-\lambda$ has a parametrix $P^+(\lambda)$ of order $-m$ and type zero in the parameter-dependent Boutet
de Monvel's calculus ${\mathcal B}^{-m,0}(\open{\overline{M}};\Lambda)$ on
$\open{\overline{M}} = \overline{M}\setminus\overline{Y}$. When restricted to the collar neighborhood
$U_{\overline{Y}} \cong (0,1)\times\overline{Y}$, this parametrix takes the form
$$
P^+(\lambda)u(x) = \frac{1}{2\pi}\iint e^{i(x-x')\xi}\tilde{p}(x,\xi,\lambda)u(x')\,dx'\,d\xi + C(\lambda)u(x)
$$
for $u \in C_0^{\infty}((0,1),C^{\infty}(\overline{Y};E))$,
where $C(\lambda) \in {\mathcal B}^{-\infty,0}((0,1)\times\overline{Y};\Lambda)$ is a parameter-dependent
regularizing singular Green operator of type zero on $(0,1)\times\overline{Y}$, and
$\tilde{p}(x,\xi,\lambda) = x^m p(x,x\xi,x^m\lambda)$ with a symbol
$$
p(x,\xi,\lambda) \in C^{\infty}([0,1),{\mathcal B}^{-m,0}(\overline{Y};\R\times\Lambda)),
$$
i.e. $p$ is a smooth function in $x \in [0,1)$ taking values in the space of
operators of order $-m$ and type zero in the parameter-dependent Boutet de Monvel's calculus on
$\overline{Y}$ (depending on the isotropic parameter $\xi$ and the anisotropic parameter $\lambda$).

The remainders $(A-\lambda)P^+(\lambda) - 1$ and $P^+(\lambda)(A-\lambda)-1$ are parameter-dependent
singular Green operators in Boutet de Monvel's calculus on $\open{\overline{M}}$ of order zero
and appropriate types (given by the standard type formula for the composition of operators).
When restricted to $U_{\overline{Y}}$, they take the form
$$
G_0(\lambda)u(x) = \frac{1}{2\pi}\iint e^{i(x-x')\xi}\tilde{g}(x,\xi,\lambda)u(x')\,dx'\,d\xi
$$
modulo parameter-dependent regularizing singular Green operators on $(0,1)\times\overline{Y}$,
where $\tilde{g}(x,\xi,\lambda) = g(x,x\xi,x^m\lambda)$ with a symbol
$$
g(x,\xi,\lambda) \in C^{\infty}([0,1),{\mathcal B}_{G}^{0,d}(\overline{Y};\R\times\Lambda)),
$$
i.e. $g$ is smooth in $x \in [0,1)$ taking values in the space of parameter-dependent singular Green
operators in Boutet de Monvel's calculus on $\overline{Y}$ of order zero and type $d$, where
$d = 0$ or $d = m$, respectively, and $\xi \in \R$ is again the isotropic parameter, while $\lambda \in \Lambda$
is the anisotropic parameter.

Observe, in particular, that $P^+(\lambda)$ has a well-defined homogeneous principal $c$-symbol $\csym(P^+)(z,\zeta,\lambda)$
on $\bigl(\cT^*\overline{M}\times\Lambda\bigr)\setminus 0$, as well as a principal
$c$-boundary symbol, which is a (twisted) homogeneous section
$$
\csymb(P^+)(z',\zeta',\lambda) : {}^c\S_+ \otimes \cpi^*E \to {}^c\S_+ \otimes \cpi^*E
$$
on $\bigl(\cT^*\overline{N}\times\Lambda\bigr)\setminus 0$.
\end{lemma}

\begin{proposition}\label{Parametrix1}
There exists a matrix of parameter-dependent generalized singular Green operators
$$
{\mathcal G}_1(\lambda) = \begin{pmatrix} G_1(\lambda) & K_1(\lambda) & \cdots & K_K(\lambda) \end{pmatrix} :
\begin{array}{c} C_0^{\infty}(\open{\overline{M}};E) \\ \oplus \\ \bigoplus\limits_{j=1}^KC_0^{\infty}(\open{\overline{N}};F_j) \end{array}
\to C^{\infty}(\open{\overline{M}};E)
$$
in Boutet de Monvel's calculus
of orders $-m$, $-m_1-\frac{1}{2}$, \ldots, $-m_K-\frac{1}{2}$ and type zero in
${\mathcal B}^{\ast}_G(\open{\overline{M}};\Lambda)$, such that its restriction to the collar neighborhood
$U_{\overline{Y}}$ of $\overline{Y}$ is (modulo a regularizing parameter-dependent
generalized singular Green operator) of the form
$$
{\mathcal G}_1(\lambda)\begin{pmatrix} u \\ v_1 \\ \vdots \\ v_K \end{pmatrix}(x) =
\frac{1}{2\pi}\iint e^{i(x-x')\xi}\tilde{g}(x,\xi,\lambda)\begin{pmatrix} u(x') \\ v_1(x') \\ \vdots \\ v_K(x') \end{pmatrix}\,dx'\,d\xi
$$
for $u \in C_0^{\infty}((0,1),C^{\infty}(\overline{Y};E))$ and
$v_j \in C_0^{\infty}((0,1),C^{\infty}(\partial\overline{Y};F_j))$, $j = 1,\ldots,K$, where
$$
\tilde{g}(x,\xi,\lambda) = \begin{pmatrix} x^{m}g(x,x\xi,x^m\lambda) & x^{m_1}k_1(x,x\xi,x^m\lambda) &
\cdots & x^{m_K}k_K(x,x\xi,x^m\lambda) \end{pmatrix},
$$
and $g(x,\xi,\lambda)$ as well as the $k_j(x,\xi,\lambda)$, $j = 1,\ldots,K$, are smooth with
respect to $x \in [0,1)$ taking values in the parameter-dependent generalized singular Green operators of orders
$-m$ and $-m_j-\frac{1}{2}$, $j=1,\ldots,K$, and type zero in Boutet de Monvel's
calculus on $\overline{Y}$ (depending on the isotropic parameter $\xi \in \R$ and the anisotropic
parameter $\lambda \in \Lambda$).

The operator family
$$
{\mathcal B}_1(\lambda) = \begin{pmatrix} P^+(\lambda) + G_1(\lambda) & K_1(\lambda) & \cdots & K_K(\lambda) \end{pmatrix} :
\begin{array}{c} C_0^{\infty}(\open{\overline{M}};E) \\ \oplus \\ \bigoplus\limits_{j=1}^KC_0^{\infty}(\open{\overline{N}};F_j) \end{array}
\to C^{\infty}(\open{\overline{M}};E)
$$
is a parameter-dependent parametrix in Boutet de Monvel's calculus of the boundary value problem
\eqref{BVPCone} in the sense that the remainders
$$
\begin{pmatrix} A-\lambda \\ T \end{pmatrix}{\mathcal B}_1(\lambda) - 1, \quad
{\mathcal B}_1(\lambda)\begin{pmatrix} A-\lambda \\ T \end{pmatrix} - 1 \in
{\mathcal B}^{-\infty}(\open{\overline{M}};\Lambda)
$$
are parameter-dependent regularizing generalized singular Green operators in Boutet de Monvel's calculus on
$\open{\overline{M}}$.
\end{proposition}
\begin{proof}
For any $d \in \N_0$ we consider the space $x^{-\varrho}{\mathcal B}^{\mu,d}_G$ of parameter-dependent
generalized singular Green operators in Boutet de Monvel's calculus on $\open{\overline{M}}$
which consist of matrix entries of the following form:
\begin{itemize}
\item \emph{Operators in the interior}\/:
$$
G(\lambda) : C_0^{\infty}(\open{\overline{M}};E) \to C^{\infty}(\open{\overline{M}};E)
$$
is a parameter-dependent singular Green operator of order $\mu$ and type $d$ in
Boutet de Monvel's calculus on $\open{\overline{M}}$, and when restricted to $U_{\overline{Y}}$
it takes the form
$$
G(\lambda)u(x) = \frac{1}{2\pi}\iint e^{i(x-x')\xi}\tilde{g}(x,\xi,\lambda)u(x')\,dx'\,d\xi
$$
(modulo a regularizing parameter-dependent singular Green operator of type $d$ in $U_{\overline{Y}}$)
with a symbol $\tilde{g}(x,\xi,\lambda) = x^{-\varrho}g(x,x\xi,x^m\lambda)$, where
$g(x,\xi,\lambda)$ is smooth with respect to $x \in [0,1)$ taking values in the
space ${\mathcal B}_G^{\mu,d}(\overline{Y};\R\times\Lambda)$ of parameter-dependent
singular Green operators of order $\mu$ and type $d$ in Boutet de Monvel's calculus
on $\overline{Y}$, depending on the isotropic parameter $\xi \in \R$ and the
anisotropic parameter $\lambda \in \Lambda$.
\item \emph{Trace operators}\/:
$$
T(\lambda) : C_0^{\infty}(\open{\overline{M}};E) \to C^{\infty}(\open{\overline{N}};F)
$$
is a parameter-dependent trace operator of order $\mu$ and type $d$ in Boutet de
Monvel's calculus on $\open{\overline{M}}$, and when restricted to $U_{\overline{Y}}$ it
takes the form
$$
T(\lambda)u(x) = \frac{1}{2\pi}\iint e^{i(x-x')\xi}\tilde{t}(x,\xi,\lambda)u(x')\,dx'\,d\xi
$$
(modulo a regularizing parameter-dependent trace operator of type $d$ in $U_{\overline{Y}}$)
with a symbol $\tilde{t}(x,\xi,\lambda) = x^{-\varrho+\frac{1}{2}}t(x,x\xi,x^m\lambda)$,
where $t(x,\xi,\lambda)$ is smooth with respect to $x \in [0,1)$ taking values in
the space of parameter-dependent trace operators of order $\mu$ and type $d$ in
Boutet de Monvel's calculus on $\overline{Y}$.
\item \emph{Potential operators}\/:
$$
K(\lambda) : C_0^{\infty}(\open{\overline{N}};F) \to C^{\infty}(\open{\overline{M}};E)
$$
is a parameter-dependent potential operator of order $\mu$ in Boutet de Monvel's
calculus on $\open{\overline{M}}$, and when restricted to the collar neighborhood it
takes the form
$$
K(\lambda)v(x) = \frac{1}{2\pi}\iint e^{i(x-x')\xi}\tilde{k}(x,\xi,\lambda)v(x')\,dx'\,d\xi
$$
(modulo a regularizing parameter-dependent potential operator in $U_{\overline{Y}}$)
with a symbol $\tilde{k}(x,\xi,\lambda) = x^{-\varrho-\frac{1}{2}}k(x,x\xi,x^m\lambda)$,
where $k(x,\xi,\lambda)$ is smooth with respect to $x \in [0,1)$ taking values
in the space of parameter-dependent potential operators of order $\mu$ in
Boutet de Monvel's calculus on $\overline{Y}$.
\item \emph{Operators on the boundary}\/:
$$
Q(\lambda) : C_0^{\infty}(\open{\overline{N}};F_1) \to C^{\infty}(\open{\overline{N}};F_2)
$$
is a parameter-dependent pseudodifferential operator of order $\mu$ on $\open{\overline{N}}$,
and when restricted to the collar neighborhood it takes the form
$$
Q(\lambda)v(x) = \frac{1}{2\pi}\iint e^{i(x-x')\xi}\tilde{q}(x,\xi,\lambda)v(x')\,dx'\,d\xi
$$
(modulo a regularizing parameter-dependent pseudodifferential operator in $(0,1)\times\partial\overline{Y}$)
with a symbol $\tilde{q}(x,\xi,\lambda) = x^{-\varrho}q(x,x\xi,x^m\lambda)$, where
$q(x,\xi,\lambda)$ is smooth with respect to $x \in [0,1)$ taking values in
the space of parameter-dependent pseudodifferential operators on $\partial\overline{Y}$
of order $\mu$ (where, as in all the other cases, $\xi \in \R$ is the isotropic
parameter, and $\lambda \in \Lambda$ is the anisotropic parameter).
\end{itemize}
Observe that every ${\mathcal G}(\lambda) \in x^{-\mu}{\mathcal B}^{\mu,d}_G$ has a well defined
principal $c$-boundary symbol $\csymb(G)(z',\zeta',\lambda)$ on $\bigl(\cT^*\overline{N}\times\Lambda\bigr) \setminus 0$,
which is (twisted) homogeneous of degree $\mu$, and we have in a canonical way a split exact sequence
\begin{equation}\label{bdrysymbsplitexact}
0 \longrightarrow x^{-\mu}{\mathcal B}^{\mu-1,d}_G \longrightarrow x^{-\mu}{\mathcal B}^{\mu,d}_G \underset{\csymb}{\longrightarrow} {}^c\Sigma \longrightarrow 0.
\end{equation}
Moreover, by standard arguments in the calculus of pseudodifferential operators with (twisted)
operator-valued symbols, we see that $x^{-\varrho}{\mathcal B}^{\ast,d}_G$ is asymptotically complete, i.e.
asymptotic summation is possible within the class.
Recall that the boundary symbolic calculus in Boutet de Monvel's algebra can be formulated in terms
of twisted operator-valued symbols, where the function spaces in the normal direction are equipped
with suitable dilation group actions. The principal boundary symbols are twisted homogeneous, i.e.
homogeneous up to conjugation with the groups (cf. \cite{SzWiley98}).

Notice that the assertion of the proposition about the structure of ${\mathcal G}_1(\lambda)$ just means that
$G_1(\lambda) \in x^{m}{\mathcal B}^{-m,0}_G$, while $K_j(\lambda) \in x^{m_j+\frac{1}{2}}{\mathcal B}^{-m_j-\frac{1}{2},0}_G$
for $j=1,\ldots,K$. Moreover, for the boundary conditions in \eqref{BVPCone} we find
$\gamma B_j \in x^{-m_j-\frac{1}{2}}{\mathcal B}^{m_j+\frac{1}{2},m}_G$, $j = 1,\ldots,K$.

Let
{\small
$$
\csymb({\mathcal A})(z',\zeta',\lambda) =
\begin{pmatrix}
\csymb(A)(z',\zeta') - \lambda  \\
({}^c\gamma_0\otimes\id_{\cpi^*F_1})\csymb(B_1)(z',\zeta') \\
\vdots \\
({}^c\gamma_0\otimes\id_{\cpi^*F_K})\csymb(B_K)(z',\zeta')
\end{pmatrix} : {}^c\S_+ \otimes \cpi^*E \to
\begin{array}{c} {}^c\S_+ \otimes \cpi^*E \\ \oplus \\ \bigoplus\limits_{j=1}^K \cpi^*F_j
\end{array}
$$}
be the principal $c$-boundary symbol of \eqref{BVPCone} on $\bigl(\cT^*\overline{N}\times\Lambda\bigr)\setminus 0$.
By choosing a smooth positive definite metric on $\cT^*\overline{N}$ we can consider
$\csymb({\mathcal A})(z',\zeta',\lambda)$ for $\bigl(|\zeta'|^{2m}+ |\lambda|^2\bigr)^{1/2m} = 1$
only, and by (twisted) homogeneity we still recover the full information.
As \eqref{BVPCone} is assumed to be $c$-elliptic with parameter $\lambda \in \Lambda$
we obtain that $\csymb({\mathcal A})(z',\zeta',\lambda)$ is invertible, and as
we consider this function now only on a compact sphere bundle the standard arguments
in Boutet de Monvel's calculus can be applied to show that its inverse is of the
form
$$
\csymb({\mathcal A})(z',\zeta',\lambda)^{-1} = 
\begin{pmatrix} \csymb(P^+)(z',\zeta',\lambda) & 0 \end{pmatrix} +
\csymb({\mathcal G}_1)(z',\zeta',\lambda),
$$
where $\csymb({\mathcal G}_1)(z',\zeta',\lambda)$ already is the principal
$c$-boundary symbol of the parame\-ter-dependent singular Green operator
${\mathcal G}_1(\lambda)$ of the assertion of the proposition.
 
In order to find ${\mathcal G}_1(\lambda)$, we first use (for each entry) the split exactness of
\eqref{bdrysymbsplitexact} to obtain a matrix ${\mathcal G}_2(\lambda)$ of
generalized parameter-dependent singular Green operators with
$\csymb({\mathcal G}_1)(z',\zeta',\lambda) = \csymb({\mathcal G}_2)(z',\zeta',\lambda)$,
and set ${\mathcal B}'(\lambda) = \begin{pmatrix} P^+(\lambda) & 0 \end{pmatrix} + {\mathcal G}_2(\lambda)$.
The composition rules imply that ${\mathcal B}'(\lambda){\mathcal A}(\lambda) = 1 + {\mathcal G}'(\lambda)$
with a parameter-dependent singular Green operator
${\mathcal G}'(\lambda) \in x^0{\mathcal B}^{-1,m}_G$. The standard formal Neumann series
argument now shows that there exists ${\mathcal G}(\lambda) \in x^0{\mathcal B}^{-1,m}_G$
(properly supported, uniformly for $\lambda \in \Lambda$) such that with
$$
{\mathcal B}_1(\lambda) = (1+{\mathcal G}(\lambda)){\mathcal B}'(\lambda) =
\begin{pmatrix} P^+(\lambda) & 0 \end{pmatrix} + {\mathcal G}_1(\lambda)
$$
we have that ${\mathcal B}_1(\lambda){\mathcal A}(\lambda)-1$ is a parameter-dependent
regularizing singular Green operator in Boutet de Monvel's calculus on $\open{\overline{M}}$
of type (at most) $m$.

A right parametrix is obtained in the same way. Note that the composition
$$
{\mathcal A}(\lambda){\mathcal B}'(\lambda) - 1 =
\left(\begin{array}{c|ccc}
 G_{0,0}(\lambda) & G_{0,1}(\lambda) & \cdots & G_{0,K}(\lambda) \\ \hline
\Vsp G_{1,0}(\lambda) & G_{1,1}(\lambda) & \cdots & G_{1,K}(\lambda) \\
\vdots & \vdots & \ddots & \vdots \\
G_{K,0}(\lambda) & G_{K,1}(\lambda) & \cdots & G_{K,K}(\lambda)
\end{array}\right)
$$
with $G_{i,j}(\lambda) \in x^{-(n_i-n_j)}{\mathcal B}^{n_i-n_j-1,0}_G$, where $n_0 = m$ and
$n_j = m_j+\frac{1}{2}$, $j=1,\ldots,K$, i.e. the matrix $\bigl(G_{i,j}(\lambda)\bigr)_{i,j}$
is of order $-1$ with respect to an order convention of Douglis-Nirenberg type (see also the proof
of Lemma \ref{EinsplusGreen}). Consequently, the formal
Neumann series argument also applies in this situation (with Douglis-Nirenberg order
convention), and the proposition is proved.
\end{proof}

Modulo a regularizing parameter-dependent generalized singular Green operator of type
zero, the restriction of the parametrix ${\mathcal B}_1(\lambda)$ from Proposition
\ref{Parametrix1} to $U_{\overline{Y}}$ can also be written in the form
\begin{equation}\label{B1Darstellung}
{\mathcal B}_1(\lambda)\begin{pmatrix} u \\ v_1 \\ \vdots \\ v_K \end{pmatrix}(x) =
\frac{1}{2\pi}\iint e^{i(x-x')\xi}\tilde{b}(x,\xi,\lambda)\begin{pmatrix} {x'}^m u(x') \\ {x'}^{m_1}v_1(x') \\ \vdots \\ {x'}^{m_K}v_K(x') \end{pmatrix}\,dx'\,d\xi
\end{equation}
for $u \in C_0^{\infty}((0,1),C^{\infty}(\overline{Y};E))$ and
$v_j \in C_0^{\infty}((0,1),C^{\infty}(\partial\overline{Y};F_j))$, $j = 1,\ldots,K$,
where $\tilde{b}(x,\xi,\lambda) = b(x,x\xi,x^m\lambda)$,
$$
b(x,\xi,\lambda) = \begin{pmatrix} (p'+g')(x,\xi,\lambda) & k'_1(x,\xi,\lambda) &
\cdots & k'_K(x,\xi,\lambda) \end{pmatrix},
$$
and the entries of $b(x,\xi,\lambda)$ are smooth with respect to $x \in [0,1)$ taking values in
the parameter-dependent Boutet de Monvel's calculus on $\overline{Y}$ of type zero
and corresponding orders.
By Mellin quantization (see, e.g., \cite{KaSchu03}), the operator \eqref{B1Darstellung} has a representation
$$
{\mathcal Q}(\lambda)\begin{pmatrix} u \\ v_1 \\ \vdots \\ v_K \end{pmatrix}(x) =
\frac{1}{2\pi}\int\limits_{\Im(\sigma)=-m/2}\int\limits_{(0,1)}
\Bigl(\frac{x}{x'}\Bigr)^{i\sigma}h(x,\sigma,x^m\lambda)\begin{pmatrix} {x'}^m u(x') \\ {x'}^{m_1}v_1(x') \\ \vdots \\ {x'}^{m_K}v_K(x') \end{pmatrix}
\,\frac{dx'}{x'}\,d\sigma
$$
modulo a regularizing parameter-dependent generalized singular Green operator in Boutet de Monvel's calculus,
where the Mellin symbol $h(x,\sigma,\lambda)$ is given by the formula
\begin{equation}\label{MellinQuantization}
h(x,\sigma,\lambda) = \frac{1}{2\pi}\iint e^{-i(r-1)\xi}r^{i\sigma}\varphi(r)b(x,\xi,\lambda)\,dr\,d\xi
\end{equation}
for $r,\xi \in \R$, $\sigma \in \C$, and $\varphi \in C_0^{\infty}(\R_+)$ is a
function such that $\varphi \equiv 1$ near $r = 1$.

Pick cut-off functions $\omega,\tilde{\omega},\hat{\omega} \in C_0^{\infty}([0,1))$ near
zero with $\hat{\omega} \prec \omega \prec \tilde{\omega}$, and consider these functions
as functions on $\overline{M}$ (or $\overline{N}$) supported in $U_{\overline{Y}}$.
With the parametrix ${\mathcal B}_1(\lambda)$ from Proposition \ref{Parametrix1} we define
\begin{equation}\label{Parametrix2Def}
{\mathcal B}_2(\lambda) = \omega{\mathcal Q}(\lambda)\tilde{\omega} +
(1-\omega){\mathcal B}_1(\lambda)(1-\hat{\omega}) :
\begin{array}{c} C_0^{\infty}(\open{\overline{M}};E) \\ \oplus \\ \bigoplus\limits_{j=1}^KC_0^{\infty}(\open{\overline{N}};F_j) \end{array}
\to C^{\infty}(\open{\overline{M}};E).
\end{equation}
By construction, we obtain the following

\begin{proposition}\label{Parametrix2}
${\mathcal B}_2(\lambda)$ is a parameter-dependent parametrix in Boutet de Monvel's
calculus of \eqref{BVPCone} which is properly supported, uniformly with respect to
$\lambda \in \Lambda$.

We have ${\mathcal B}_2(\lambda) - {\mathcal B}_1(\lambda) \in {\mathcal B}^{-\infty,0}(\open{\overline{M}};\Lambda)$,
and
$$
{\mathcal B}_2(\lambda) : \begin{array}{c} x^{-m/2}H_b^s(\overline{M};E) \\ \oplus \\
\bigoplus\limits_{j=1}^Kx^{m/2-m_j}H_b^{s+m-m_j-1/2}(\overline{N};F_j) \end{array} \to
x^{m/2}H_b^{s+m}(\overline{M};E) \hookrightarrow \Dom^s_{\min}\binom{A}{T}
$$
is continuous for all $s > -\frac{1}{2}$.
\end{proposition}

In order to further refine the parametrix ${\mathcal B}_2(\lambda)$, we first recall the
notion of operator-valued symbols on the sector $\Lambda$ (general information
about such symbol classes can be found in \cite{SzNorthHolland, SzWiley98}):

\begin{definition}\label{OpvaluedSymbol}
Let ${\mathbf H}$ and $\tilde{{\mathbf H}}$ be Hilbert spaces endowed with strongly continuous
groups of isomorphisms $\{\kappa_{\varrho}\}$ and $\{\tilde{\kappa}_{\varrho}\}$,
$\varrho > 0$, respectively.

A function $g \in C^{\infty}(\Lambda,\L({\mathbf H},\tilde{{\mathbf H}}))$ is called an operator-valued
symbol of order $\mu \in \R$, if for all multi-indices $\alpha \in \N_0^2$
\begin{equation}\label{opsymbabsch}
\bigl\|\tilde{\kappa}^{-1}_{[\lambda]^{1/m}}\partial^{\alpha}_{\lambda}g(\lambda)\kappa_{[\lambda]^{1/m}}\bigr\|_{\L({\mathbf H},\tilde{{\mathbf H}})}
= O\bigl(|\lambda|^{\mu/m-|\alpha|}\bigr)
\end{equation}
as $|\lambda| \to \infty$, where $[\cdot]$ is a smooth function on $\C$ with $[\lambda] > 1$ for all
$\lambda \in \C$, and $[\lambda] = |\lambda|$ for $|\lambda| > 2$.

Moreover, for every $j \in \N_0$ there should exist (twisted) homogeneous (or $\kappa$-homogeneous) components $g_{(\mu-j)} \in
C^{\infty}(\Lambda\setminus\{0\},\L({\mathbf H},\tilde{{\mathbf H}}))$, i.e.
$$
g_{(\mu-j)}(\varrho^m\lambda) = \varrho^{\mu-j}\tilde{\kappa}_{\varrho}g_{(\mu-j)}(\lambda)\kappa_{\varrho}^{-1}
$$
for $\varrho > 0$, such that for some excision function $\chi \in C^{\infty}(\C)$
($\chi \equiv 0$ near zero and $\chi \equiv 1$ near infinity) and all $N \in \N_0$ the symbol estimates
\eqref{opsymbabsch} hold for $g(\lambda) - \sum\limits_{j=0}^{N-1}\chi(\lambda)g_{(\mu-j)}(\lambda)$
in the place of $g(\lambda)$, and $\mu$ replaced by $\mu-N$.
We sometimes write
$$
g(\lambda) \sim \sum\limits_{j=0}^{\infty}g_{(\mu-j)}(\lambda),
$$
and $g_{\wedge}(\lambda):= g_{(\mu)}(\lambda)$ is called the principal component of $g(\lambda)$.

We call the operator-valued symbol $g(\lambda)$ compact, if in all conditions above we may replace the
space of all bounded operators $\L({\mathbf H},\tilde{{\mathbf H}})$ by the ideal of compact
operators $\K({\mathbf H},\tilde{{\mathbf H}})$.

In the considerations below the Hilbert spaces ${\mathbf H}$ and $\tilde{{\mathbf H}}$ will either be function spaces on the model cone
$\overline{Y}^{\wedge}$ and $(\partial\overline{Y})^{\wedge}$, or just $\C^{N}$, and the group action
is either the normalized dilation $\kappa_{\varrho}$ from Definition \ref{kappagroup} on the function spaces, or
the trivial group action ($\tilde{\kappa}_{\varrho} \equiv \id$) on $\C^{N}$.
\end{definition}

The following Definition \ref{Greenremainders} of generalized Green remainders is essential for
understanding the structure of remainders of the parametrix construction, and for further necessary
refinement of the parametrix itself.

Let $\partial_+$ denote normal differentiation on $2\overline{M}$ (near the boundary $\partial(2\overline{M})$),
where ``normal'' refers to some Riemannian metric, smooth up to $\partial(2\overline{M})$,
that coincides near $\overline{Y}$ with $dx^2+dy^2$.

\begin{definition}\label{Greenremainders}
Let $d \in \N_0$ and $\mu \in \R$. An operator family
\begin{equation}\label{Greenremaindereq}
{\mathcal G}(\lambda) = \sum\limits_{j=0}^d{\mathcal G}_{j}(\lambda)\begin{pmatrix} \partial_+ & 0 & 0 \\ 0 & 0 & 0 \\ 0 & 0 & 0 \end{pmatrix}^j :
\begin{array}{c} C_0^{\infty}(\open{\overline{M}};E) \\ \oplus \\ C_0^{\infty}(\open{\overline{N}};F) \\ \oplus \\ \C^{N_-} \end{array} \to
\begin{array}{c} C^{\infty}(\open{\overline{M}};E) \\ \oplus \\ C^{\infty}(\open{\overline{N}};F') \\ \oplus \\ \C^{N_+} \end{array}
\end{equation}
is called a generalized Green remainder of order $\mu$ and type $d$ in the scales
$\begin{pmatrix} x^{\alpha}H^s_b \\ x^{\beta}H^s_b \\ \C^{N_-}\end{pmatrix}$ to
$\begin{pmatrix} x^{\alpha'}H^s_b \\ x^{\beta'}H^s_b \\ \C^{N_+}\end{pmatrix}$
if for all cut-off functions $\omega,\tilde{\omega} \in C_0^{\infty}([0,1))$ near zero the following
holds:
\begin{enumerate}[i)]
\item For all $j=0,\ldots,d$
$$
(1-\omega){\mathcal G}_j(\lambda),\; {\mathcal G}_j(\lambda)(1-\tilde{\omega}) \in
\bigcap\limits_{s,t \in \R}\S\left(\Lambda,\K\left(
\begin{array}{c}x^{\alpha}H^s_{b,0}(\overline{M};E) \\ \oplus \\ x^{\beta}H^s_b(\overline{N};F) \\ \oplus \\ \C^{N_-}\end{array},
\begin{array}{c}x^{\alpha'}H^t_{b}(\overline{M};E) \\ \oplus \\ x^{\beta'}H^t_b(\overline{N};F') \\ \oplus \\ \C^{N_+}\end{array}
\right)\right).
$$
\item For all $j=0,\ldots,d$,
$$
g_j(\lambda) = \omega{\mathcal G}_j(\lambda)\tilde{\omega} : 
\begin{array}{c} C_0^{\infty}(\open{\overline{Y}}^{\wedge};E) \\ \oplus \\ C_0^{\infty}(\partial\open{\overline{Y}}^{\wedge};F) \\ \oplus \\ \C^{N_-} \end{array} \to
\begin{array}{c} C^{\infty}(\open{\overline{Y}}^{\wedge};E) \\ \oplus \\ C^{\infty}(\partial\open{\overline{Y}}^{\wedge};F') \\ \oplus \\ \C^{N_+} \end{array}
$$
is a Green symbol, i.e. a compact operator-valued symbol
$$
\begin{array}{c} {}_{\delta}\K_0^{s,\alpha}(\overline{Y}^{\wedge};E) \\ \oplus \\ {}_{\delta}\K^{s,\beta}((\partial\overline{Y})^{\wedge};F) \\ \oplus \\ \C^{N_-} \end{array} \to
\begin{array}{c} {}_{\delta'}\K^{t,\alpha'}(\overline{Y}^{\wedge};E) \\ \oplus \\ {}_{\delta'}\K^{t,\beta'}((\partial\overline{Y})^{\wedge};F') \\ \oplus \\ \C^{N_+} \end{array}
$$
of order $\mu \in \R$ for all $s,t,\delta,\delta' \in \R$.

\medskip

Here, for $s,\delta,\alpha \in \R$, we write
${}_{\delta}\K^{s,\alpha} = \omega\K^{s,\alpha} + (1-\omega)x^{-\delta}\K^{s,\alpha}$, as well as
${}_{\delta}\K_0^{s,\alpha} = \bigl({}_{-\delta}\K^{-s,-\alpha-m}\bigr)'$,
where the dual space is to be understood with respect to the pairing induced by the
$x^{-m/2}L_b^2$-inner product.
\end{enumerate}

\bigskip

\noindent
Correspondingly, the operator family \eqref{Greenremaindereq} is called a generalized Green remainder
of order $\mu$ and type $d$ in the scales
$\begin{pmatrix} x^{\alpha}H^s_b \\ x^{\beta}H^s_b \\ \C^{N_-}\end{pmatrix}$ to
$\begin{pmatrix} \Dom^s_{\min} \\ x^{\beta'}H^s_b \\ \C^{N_+}\end{pmatrix}$
if for all cut-off functions $\omega,\tilde{\omega} \in C_0^{\infty}([0,1))$ near zero the following
holds:
\begin{enumerate}[iii)]
\item[iii)] For all $j=0,\ldots,d$
$$
(1-\omega){\mathcal G}_j(\lambda),\; {\mathcal G}_j(\lambda)(1-\tilde{\omega}) \in
\bigcap\limits_{\substack{s \in \R \\ t > -\frac{1}{2}}}\S\left(\Lambda,\K\left(
\begin{array}{c}x^{\alpha}H^s_{b,0}(\overline{M};E) \\ \oplus \\ x^{\beta}H^s_b(\overline{N};F) \\ \oplus \\ \C^{N_-}\end{array},
\begin{array}{c}\Dom^t_{\min}\binom{A}{T} \\ \oplus \\ x^{\beta'}H^t_b(\overline{N};F') \\ \oplus \\ \C^{N_+}\end{array}
\right)\right).
$$
\item[iv)] For all $j=0,\ldots,d$, $g_j(\lambda) = \omega{\mathcal G}_j(\lambda)\tilde{\omega}$
is a compact operator-valued symbol
$$
\begin{array}{c} {}_{\delta}\K_0^{s,\alpha}(\overline{Y}^{\wedge};E) \\ \oplus \\ {}_{\delta}\K^{s,\beta}((\partial\overline{Y})^{\wedge};F) \\ \oplus \\ \C^{N_-} \end{array} \to
\begin{array}{c} {}_{\delta'}\Dom^t_{\wedge,\min}\binom{A_{\wedge}}{T_{\wedge}} \\ \oplus \\ {}_{\delta'}\K^{t,\beta'}((\partial\overline{Y})^{\wedge};F') \\ \oplus \\ \C^{N_+} \end{array}
$$
of order $\mu$ for all $s,t,\delta,\delta' \in \R$, $t > -\frac{1}{2}$, where analogously to the above
$$
{}_{\delta'}\Dom^t_{\wedge,\min}\binom{A_{\wedge}}{T_{\wedge}} =
\omega\Dom^t_{\wedge,\min}\binom{A_{\wedge}}{T_{\wedge}} + (1-\omega)x^{-\delta'}\Dom^t_{\wedge,\min}\binom{A_{\wedge}}{T_{\wedge}}.
$$
\end{enumerate}

In ii) and iv), the property of being an operator-valued symbol always refers to the
group action
$\begin{pmatrix} \kappa_{\varrho} & 0 & 0 \\ 0 & \kappa_{\varrho} & 0 \\ 0 & 0 & \id \end{pmatrix}$,
i.e. we consider the normalized dilation group $\kappa_{\varrho}$ from Definition \ref{kappagroup}
on the function spaces, and the trivial group action on $\C^{N_{\pm}}$.
Moreover, multiplication of the ${\mathcal G}_j(\lambda)$ by the cut-off function $\omega$ (or $\tilde{\omega}$) above is to be understood as
multiplication by the diagonal matrix
$\begin{pmatrix} \omega & 0 & 0 \\ 0 & \omega & 0 \\ 0 & 0 & \id \end{pmatrix}$,
while $1$ always is the identity matrix. In particular, $1 - \omega$ is to be understood as multiplication by
the matrix $\begin{pmatrix} 1-\omega & 0 & 0 \\ 0 & 1-\omega & 0 \\ 0 & 0 & 0 \end{pmatrix}$.

It is needless to say that these definitions also apply to each entry of the matrix ${\mathcal G}(\lambda)$
separately (which corresponds, e.g., to $N_{\pm} = 0$ or $F=0$, $F'=0$).
For $N_- = N_+ = 0$, every generalized Green remainder of order $\mu$ and type $d$ is an element
of ${\mathcal B}^{-\infty,d}(\open{\overline{M}};\Lambda)$, the class of regularizing parameter-dependent
generalized singular Green operators in Boutet de Monvel's calculus on $\open{\overline{M}}$ of type $d$,
i.e. we pass to a specific class of admissible remainders here. Both meanings of
``Green'' should not be mixed up, and it will always be clear from the context which notion applies.

It is not hard to prove that every generalized Green remainder ${\mathcal G}(\lambda)$ of order
$\mu$ has an associated sequence ${\mathcal G}_{(\mu-j)}(\lambda)$ of (twisted) homogeneous
components of order $\mu-j$, $j \in \N_0$ --- namely the components of the operator-valued
symbol $\omega{\mathcal G}(\lambda)\tilde{\omega}$ --- and these components are unique, i.e. they
do not depend on the choice of cut-off functions $\omega,\tilde{\omega} \in C_0^{\infty}([0,1))$
(the argument is similar to Lemma 5.19 in \cite{GKM2}).
Consequently, we write
$$
{\mathcal G}(\lambda) \sim \sum\limits_{j=0}^{\infty}{\mathcal G}_{(\mu-j)}(\lambda),
$$
and call ${\mathcal G}_{\wedge}(\lambda):= {\mathcal G}_{(\mu)}(\lambda)$ the principal component
of ${\mathcal G}(\lambda)$.
Moreover, the intersection of all generalized Green remainders of order
$\mu$ and fixed type $d$ consists of the so called regularizing generalized Green remainders of type $d$.
\end{definition}

Observe, in particular, that the operator
$$
{\mathcal G}(\lambda) :
\begin{array}{c}x^{\alpha}H^s_{b}(\overline{M};E) \\ \oplus \\ x^{\beta}H^s_b(\overline{N};F) \\ \oplus \\ \C^{N_-}\end{array} \to
\begin{array}{c}x^{\alpha'}H^t_{b}(\overline{M};E) \\ \oplus \\ x^{\beta'}H^t_b(\overline{N};F') \\ \oplus \\ \C^{N_+}\end{array}
$$
is compact for all $\lambda \in \Lambda$ and all $s,t \in \R$, $s > d-\frac{1}{2}$,
where ${\mathcal G}(\lambda)$ is any generalized Green remainder of
type $d \in \N_0$ in the scales
$\begin{pmatrix} x^{\alpha}H^s_b \\ x^{\beta}H^s_b \\ \C^{N_-}\end{pmatrix}$ to
$\begin{pmatrix} x^{\alpha'}H^s_b \\ x^{\beta'}H^s_b \\ \C^{N_+}\end{pmatrix}$.
Analogously,
$$
{\mathcal G}(\lambda) :
\begin{array}{c}x^{\alpha}H^s_{b}(\overline{M};E) \\ \oplus \\ x^{\beta}H^s_b(\overline{N};F) \\ \oplus \\ \C^{N_-}\end{array} \to
\begin{array}{c}\Dom^t_{\min}\binom{A}{T} \\ \oplus \\ x^{\beta'}H^t_b(\overline{N};F') \\ \oplus \\ \C^{N_+}\end{array}
$$
is compact for all $\lambda \in \Lambda$ and all $s > d-\frac{1}{2}$, $t > -\frac{1}{2}$, for any
generalized Green remainder ${\mathcal G}(\lambda)$ of type $d \in \N_0$ in the scales
$\begin{pmatrix} x^{\alpha}H^s_b \\ x^{\beta}H^s_b \\ \C^{N_-}\end{pmatrix}$ to
$\begin{pmatrix} \Dom^s_{\min} \\ x^{\beta'}H^s_b \\ \C^{N_+}\end{pmatrix}$.

It is easy to see from the definition that the generalized Green remainders form an algebra,
i.e. the composition ${\mathcal G}_1(\lambda){\mathcal G}_2(\lambda)$ of generalized Green remainders
${\mathcal G}_j(\lambda)$ of orders $\mu_j$ and types $d_j$, $j=1,2$, is a generalized
Green remainder of order $\mu_1+\mu_2$ and type $d_2$, and the principal component of
the composition equals the product of the principal components of the factors (here it is of
course assumed that the scales fit together such that the composition makes sense).

\begin{lemma}\label{EinsplusGreen}
Let
$$
{\mathcal G}(\lambda) = \left(\begin{array}{c|ccc|c} {\mathcal G}_{0,0}(\lambda) & {\mathcal G}_{0,1}(\lambda) & \cdots & {\mathcal G}_{0,K}(\lambda) & {\mathcal G}_{0,K+1}(\lambda) \\ \hline
\Vsp
{\mathcal G}_{1,0}(\lambda) & {\mathcal G}_{1,1}(\lambda) & \cdots & {\mathcal G}_{1,K}(\lambda) & {\mathcal G}_{1,K+1}(\lambda) \\
\vdots & \vdots & \ddots & \vdots & \vdots \\
{\mathcal G}_{K,0}(\lambda) & {\mathcal G}_{K,1}(\lambda) & \cdots & {\mathcal G}_{K,K}(\lambda) & {\mathcal G}_{K,K+1}(\lambda) \\ \hline
\Vsp
{\mathcal G}_{K+1,0}(\lambda) & {\mathcal G}_{K+1,1}(\lambda) & \cdots & {\mathcal G}_{K+1,K}(\lambda) & {\mathcal G}_{K+1,K+1}(\lambda)
\end{array}\right)
$$
be a matrix of generalized Green remainders of fixed type $d \in \N_0$, and let $n_i-n_j$ be the
order of ${\mathcal G}_{i,j}(\lambda)$, $i,j=0,\ldots,K+1$. Here ${\mathcal G}(\lambda)$
is an operator
$$
\begin{array}{c}
x^{-m/2}H_b^s(\overline{M};E) \\ \oplus \\ \bigoplus\limits_{j=1}^K x^{m/2-m_j}H_b^{s+m-m_j-1/2}(\overline{N};F_j) \\ \oplus \\ \C^{N} \end{array}
\to
\begin{array}{c}
x^{-m/2}H_b^s(\overline{M};E) \\ \oplus \\ \bigoplus\limits_{j=1}^K x^{m/2-m_j}H_b^{s+m-m_j-1/2}(\overline{N};F_j) \\ \oplus \\ \C^{N} \end{array}
$$
for $s > d-\frac{1}{2}$, and the ${\mathcal G}_{i,j}(\lambda)$ are assumed to be Green in the
corresponding scales of spaces.

Let ${\mathcal G}_{\wedge}(\lambda)$ be the matrix of principal parts of ${\mathcal G}(\lambda)$, an
operator family in the spaces
$$
\begin{array}{c}
\K^{s,-m/2}(\overline{Y}^{\wedge};E) \\ \oplus \\ \bigoplus\limits_{j=1}^K \K^{s+m-m_j-1/2,m/2-m_j}((\partial\overline{Y})^{\wedge};F_j) \\ \oplus \\ \C^{N} \end{array}
\to
\begin{array}{c}
\K^{s,-m/2}(\overline{Y}^{\wedge};E) \\ \oplus \\ \bigoplus\limits_{j=1}^K \K^{s+m-m_j-1/2,m/2-m_j}((\partial\overline{Y})^{\wedge};F_j) \\ \oplus \\ \C^{N} \end{array}
$$
for $s > d-\frac{1}{2}$ and $\lambda \in \Lambda\setminus\{0\}$, and assume that
$1+{\mathcal G}_{\wedge}(\lambda)$ is invertible for $\lambda \in \Lambda\setminus\{0\}$.

Then $1+{\mathcal G}(\lambda)$ is invertible for large $\lambda \in \Lambda$, and there exists a
matrix $\tilde{{\mathcal G}}(\lambda) = \bigl(\tilde{{\mathcal G}}_{i,j}(\lambda)\bigr)$
of generalized Green remainders $\tilde{{\mathcal G}}_{i,j}(\lambda)$ of the same type $d \in \N_0$ and
order $n_i-n_j$, $i,j=0,\ldots,K+1$, such that $\bigl(1+{\mathcal G}(\lambda)\bigr)^{-1} = 1 + \tilde{{\mathcal G}}(\lambda)$
for $|\lambda|>0$ sufficiently large.
\end{lemma}
\begin{proof}
Note that the matrix ${\mathcal G}(\lambda)$ is of order zero with respect to an order convention of
Douglis-Nirenberg type: A matrix $\bigl({\mathcal G}_{i,j}(\lambda)\bigr)$ is to be considered of
order $\mu \in \R$ if ${\mathcal G}_{i,j}(\lambda)$ has order $n_i-n_j+\mu$, and
correspondingly a matrix ${\mathcal G}_{\wedge}(\lambda)$ is $\kappa$-homogeneous of Douglis-Nirenberg order $\mu$ if
{\small
$$
{\mathcal G}_{\wedge}(\varrho^m\lambda) = \varrho^{\mu}\begin{pmatrix} \varrho^{n_0}\kappa_{\varrho} & \cdots & 0 & 0 \\
\vdots & \ddots & \vdots & \vdots \\
0 & \cdots & \varrho^{n_K}\kappa_{\varrho} & 0 \\
 0 & \cdots & 0 & \varrho^{n_{K+1}} \end{pmatrix}
\!{\mathcal G}_{\wedge}(\lambda)\!
\begin{pmatrix} \varrho^{n_0}\kappa_{\varrho} & \cdots & 0 & 0 \\
\vdots & \ddots & \vdots & \vdots \\
0 & \cdots & \varrho^{n_K}\kappa_{\varrho} & 0 \\
 0 & \cdots & 0 & \varrho^{n_{K+1}} \end{pmatrix}^{-1}
$$}
for $\varrho > 0$, where in our situation $\mu = 0$.

Therefore we see that the inverse of $1 + {\mathcal G}_{\wedge}(\lambda)$ is of the form
$1 + \tilde{{\mathcal G}}_{\wedge}(\lambda)$, where $\tilde{{\mathcal G}}_{\wedge}(\lambda)$ is
$\kappa$-homogeneous of (Douglis-Nirenberg) order zero, and from the identity
$$
\bigl(1+{\mathcal G}_{\wedge}(\lambda)\bigr)^{-1} = 1 - {\mathcal G}_{\wedge}(\lambda) + {\mathcal G}_{\wedge}(\lambda)\bigl(1+{\mathcal G}_{\wedge}(\lambda)\bigr)^{-1}{\mathcal G}_{\wedge}(\lambda)
$$
we see that $\tilde{{\mathcal G}}_{\wedge}(\lambda)$ is a principal Green symbol of type $d$ and
Douglis-Nirenberg order zero.

With a cut-off function $\omega \in C_0^{\infty}([0,1))$ and a function $\chi \in C^{\infty}(\C)$
with $\chi \equiv 0$ near zero and $\chi \equiv 1$ near infinity define
$$
{\mathcal G}'(\lambda):= \begin{pmatrix} \omega & 0 & 0 \\ 0 & \omega & 0 \\ 0 & 0 & 1 \end{pmatrix}\chi(\lambda)\tilde{{\mathcal G}}_{\wedge}(\lambda)
\begin{pmatrix} \omega & 0 & 0 \\ 0 & \omega & 0 \\ 0 & 0 & 1 \end{pmatrix}.
$$
Then ${\mathcal G}'(\lambda)$ is a generalized Green remainder of Douglis-Nirenberg order zero and
type $d$, and
$$
\bigl(1+{\mathcal G}(\lambda)\bigr)\bigl(1+{\mathcal G}'(\lambda)\bigr) - 1, \;
\bigl(1+{\mathcal G}'(\lambda)\bigr)\bigl(1+{\mathcal G}(\lambda)\bigr) - 1
$$
are generalized Green remainders of Douglis-Nirenberg order $-1$ and type $d$.

As the classes of generalized Green remainders are asymptotically complete, a standard
formal Neumann series argument now shows that $1 + {\mathcal G}(\lambda)$ has an inverse of the
asserted form modulo regularizing generalized Green remainders of type $d$, and as these are
rapidly decreasing in the norm as $|\lambda| \to \infty$ the assertion of the lemma regarding
the invertibility of $1+{\mathcal G}(\lambda)$ for large $\lambda$ follows.

Finally, it remains to note that if $\hat{{\mathcal G}}(\lambda)$ is a regularizing generalized Green remainder
of type $d$, the inverse of $1 + \hat{{\mathcal G}}(\lambda)$ for large $\lambda$ is of the form
$$
\bigl(1 + \hat{{\mathcal G}}(\lambda)\bigr)^{-1} = 1 - \hat{{\mathcal G}}(\lambda) +
\hat{{\mathcal G}}(\lambda)\chi(\lambda)\bigl(1 + \hat{{\mathcal G}}(\lambda)\bigr)^{-1}\hat{{\mathcal G}}(\lambda),
$$
where $\chi \in C^{\infty}(\C)$ is an excision function as above, and
$$
-\hat{{\mathcal G}}(\lambda) + \hat{{\mathcal G}}(\lambda)\chi(\lambda)\bigl(1 + \hat{{\mathcal G}}(\lambda)\bigr)^{-1}\hat{{\mathcal G}}(\lambda)
$$
is obviously a regularizing generalized Green remainder of type $d$.
\end{proof}

Let $\hat{{\mathcal A}}_0(\sigma)$ be the conormal symbol of \eqref{BVPCone}.
The $c$-ellipticity implies that the inverse $\hat{{\mathcal A}}_0^{-1}(\sigma)$
is a finitely meromorphic Fredholm function on $\C$,
and there exists a sufficiently small $\eps_0>0$
such that $\hat{{\mathcal A}}_0(\sigma)$ is invertible in
$$
 \set{\sigma \in \C \st -m/2-\eps_0 < \Im\sigma <
 -m/2+\eps_0,\; \Im\sigma \neq -m/2}.
$$
Define
\begin{equation}\label{SmoothMellRest}
h_0(\sigma)= \hat{{\mathcal A}}_0^{-1}(\sigma) - h(0,\sigma,0),
\end{equation}
where $h$ is the holomorphic Mellin symbol from \eqref{MellinQuantization}.
Then $h_0(\sigma)$ is finitely meromorphic in $\C$ taking values in
${\mathcal B}^{-\infty,0}(\overline{Y})$, and it is rapidly decreasing as $|\Re\sigma| 
\to \infty$, uniformly for $\Im\sigma$ in compact intervals (this is subject to the
composition behavior of cone pseudodifferential operators without parameters; a proof
can be found in \cite{SchroSchu}).
Moreover, the set
$$
 \set{\sigma \in \C \st -m/2-\eps_0 < \Im\sigma <
 -m/2+\eps_0,\; \Im\sigma \neq -m/2}
$$
is free of poles of $h_0(\sigma)$.

Let $\omega \in C_0^{\infty}([0,1))$ be a cut-off function and $0 < \eps < \eps_0$. Define
\begin{align*}
m(x,x',\sigma,\lambda)&:= \omega(x[\lambda]^{1/m})h_0(\sigma)\omega(x'[\lambda]^{1/m})
\begin{pmatrix} {x'}^{m} & 0 & \cdots & 0 \\ 0 & {x'}^{m_1} & \cdots & 0 \\ \vdots & \vdots & \ddots & \vdots \\
0 & 0 & \cdots & {x'}^{m_K} \end{pmatrix}, \\
m_{\wedge}(x,x',\sigma,\lambda)&:= \omega(x|\lambda|^{1/m})h_0(\sigma)\omega(x'|\lambda|^{1/m})
\begin{pmatrix} {x'}^{m} & 0 & \cdots & 0 \\ 0 & {x'}^{m_1} & \cdots & 0 \\ \vdots & \vdots & \ddots & \vdots \\
0 & 0 & \cdots & {x'}^{m_K} \end{pmatrix},
\end{align*}
and associated operators
\begin{align}
M(\lambda) &: \begin{array}{c} C_0^{\infty}(\open{\overline{M}};E) \\ \oplus \\ \bigoplus\limits_{j=1}^K C_0^{\infty}(\open{\overline{N}};F_j) \end{array}
\to C^{\infty}(\open{\overline{M}};E), \\
M_{\wedge}(\lambda) &: \begin{array}{c} C_0^{\infty}(\open{\overline{Y}}^{\wedge};E) \\ \oplus \\ \bigoplus\limits_{j=1}^K C_0^{\infty}(\partial\open{\overline{Y}}^{\wedge};F_j) \end{array}
\to C^{\infty}(\open{\overline{Y}}^{\wedge};E)
\end{align} 
via 
\begin{equation*}
\begin{pmatrix} u \\ v_1 \\ \vdots \\ v_K \end{pmatrix} \mapsto \biggl(\frac{1}{2{\pi}}\!
\int\limits_{\Im\sigma = -m/2+\eps\,} 
\int\limits_{\R_+}\Bigl(\frac{x}{x'}\Bigr)^{i\sigma} m_{(\wedge)}(x,x',\sigma,\lambda)
\begin{pmatrix} u(x') \\ v_1(x') \\ \vdots \\ v_K(x') \end{pmatrix} \, \frac{dx'}{x'}\,d\sigma\biggr).
\end{equation*}
$M(\lambda)$ is a regularizing parameter-dependent generalized singular Green operator in Boutet de Monvel's
calculus of type zero, and since the function $\omega(x[\lambda]^{1/m})$ is supported
in the collar $[0,1){\times}\overline{Y}$, $M(\lambda)$ can be
regarded as an operator both on the manifold and the model cone.
Observe, moreover, that the components of the matrix $M_\wedge(\lambda)$ are $\kappa$-homogeneous
of degrees $-m$, $-m_1$, \ldots, $-m_K$.

We define a refinement of the parameter-dependent parametrix ${\mathcal B}_2(\lambda)$ of \eqref{BVPCone}
from Proposition \ref{Parametrix2} via
\begin{equation}\label{Parametrix3Definition}
{\mathcal B}_3(\lambda):= {\mathcal B}_2(\lambda) + M(\lambda) :
\begin{array}{c} C_0^{\infty}(\open{\overline{M}};E) \\ \oplus \\ \bigoplus\limits_{j=1}^K C_0^{\infty}(\open{\overline{N}};F_j) \end{array}
\to C^{\infty}(\open{\overline{M}};E),
\end{equation}
and correspondingly let
\begin{equation}\label{Parametrix3principalpart}
{\mathcal B}_{3,\wedge}(\lambda):= {\mathcal B}_{2,\wedge}(\lambda) + M_{\wedge}(\lambda),
\end{equation}
$\lambda \in \Lambda\setminus\{0\}$, be the principal part of ${\mathcal B}_3(\lambda)$, where
$$
{\mathcal B}_{2,\wedge}(\lambda)\begin{pmatrix} u \\ v_1 \\ \vdots \\ v_K \end{pmatrix}(x) =
\frac{1}{2\pi}\int\limits_{\Im(\sigma)=-m/2}\int\limits_{\R_+}
\Bigl(\frac{x}{x'}\Bigr)^{i\sigma}h(0,\sigma,x^m\lambda)\begin{pmatrix} {x'}^m u(x') \\ {x'}^{m_1}v_1(x') \\ \vdots \\ {x'}^{m_K}v_K(x') \end{pmatrix}
\,\frac{dx'}{x'}\,d\sigma.
$$

\begin{proposition}\label{Parametrix3}
Let
$$
{\mathcal G}(\lambda) = \bigl({\mathcal G}_{i,j}(\lambda)\bigr)_{\begin{subarray}{l}i=0,\ldots,K+1 \\ j=0,1\end{subarray}} :
\begin{array}{c} \Dom^s_{\min}\binom{A}{T} \\ \oplus \\ \C^{N_-} \end{array} \to
\begin{array}{c} x^{-m/2}H_b^s(\overline{M};E) \\ \oplus \\ \bigoplus\limits_{j=1}^K x^{m/2-m_j}H_b^{s+m-m_j-1/2}(\overline{N};F_j) \\ \oplus \\ \C^{N_+} \end{array}
$$
be a matrix of generalized Green remainders, where ${\mathcal G}_{i,0}(\lambda) \equiv 0$
and ${\mathcal G}_{i,1}(\lambda)$ has order $m_i$, $i=0,\ldots,K+1$, with $m_0 = m$ and arbitrary
$m_{K+1} \in \R$.

Moreover, let
$$
\tilde{{\mathcal G}}(\lambda) = \bigl(\tilde{{\mathcal G}}_{i,j}(\lambda)\bigr)_{\begin{subarray}{l} i=0,1 \\ j=0,\ldots,K+1\end{subarray}} :
\begin{array}{c} x^{-m/2}H_b^s(\overline{M};E) \\ \oplus \\ \bigoplus\limits_{j=1}^K x^{m/2-m_j}H_b^{s+m-m_j-1/2}(\overline{N};F_j) \\ \oplus \\ \C^{N_+} \end{array}
\to
\begin{array}{c} \Dom^s_{\min}\binom{A}{T} \\ \oplus \\ \C^{N_-} \end{array}
$$
be a matrix of generalized Green remainders of type zero, where $\tilde{\mathcal G}_{i,j}(\lambda)$
has order $-m_j$, $i=0,1$, $j=0,\ldots,K+1$.

Then
\begin{align*}
\Biggl[\left(\begin{array}{c|c} A-\lambda & 0 \\ \hline \Vsp
T & 0 \\ \hline \Vsp
0 & 0 \end{array}\right) + {\mathcal G}(\lambda)\Biggr] \cdot
\Biggl[\left(\begin{array}{cc} {\mathcal B}_3(\lambda) & 0 \\ \hline \Vsp
0 & 0 \end{array}\right) + \tilde{{\mathcal G}}(\lambda)\Biggr] &= 1 + \hat{{\mathcal G}}(\lambda), \\
\Biggl[\left(\begin{array}{c|c} A_{\wedge}-\lambda & 0 \\ \hline \Vsp
T_{\wedge} & 0 \\ \hline \Vsp
0 & 0 \end{array}\right) + {\mathcal G}_{\wedge}(\lambda)\Biggr] \cdot
\Biggl[\left(\begin{array}{cc} {\mathcal B}_{3,\wedge}(\lambda) & 0 \\ \hline \Vsp
0 & 0 \end{array}\right) + \tilde{{\mathcal G}}_{\wedge}(\lambda)\Biggr] &= 1 + \hat{{\mathcal G}}_{\wedge}(\lambda),
\end{align*}
where $\hat{{\mathcal G}}(\lambda) = \bigl(\hat{{\mathcal G}}_{i,j}(\lambda)\bigr)_{i,j=0,\ldots,K+1}$
is a matrix of generalized Green remainders of type zero, and $\hat{{\mathcal G}}_{i,j}(\lambda)$ has
order $m_i-m_j$.
\end{proposition}
\begin{proof}
The proof of this proposition amounts in understanding the structure of the following compositions:
\begin{enumerate}[i)]
\item $(A-\lambda)\tilde{{\mathcal G}}(\lambda)$ and $\gamma B_j\tilde{{\mathcal G}}(\lambda)$,
$j=1,\ldots,K$.
\item ${\mathcal G}(\lambda){\mathcal B}_3(\lambda)$.
\item $\binom{A-\lambda}{T}{\mathcal B}_3(\lambda)$.
\end{enumerate}
In i) and ii), ${\mathcal G}(\lambda)$ and $\tilde{{\mathcal G}}(\lambda)$ are appropriate
(matrices of) generalized Green remainders. Using the identity
$$
\Dom^s_{\min}\binom{A}{T} = \Dom^s_{\max}\binom{A}{T}\cap\bigcap\limits_{\eps > 0}x^{m/2-\eps}H_b^{s+m}(\overline{M};E),
$$
as well as (anisotropic modifications of) the results about the structure and composition behavior of parameter-dependent pseudodifferential
cone operators in the edge symbolic calculus in scales of Sobolev spaces from \cite{KaSchu03}, we
can employ here the same strategy as in the boundaryless case, see \cite{GKM2}:
\begin{itemize}
\item Using the expansions \eqref{exp1}, \eqref{exp2}, and \eqref{exp3} and a similar expansion for
${\mathcal B}_2(\lambda)$ (Taylor expansion of the symbol $h(x,\sigma,\lambda)$ from
\eqref{MellinQuantization} in $x=0$), an inspection
of the proof of Lemma 5.20 in \cite{GKM2} reveals that the analogue of this lemma also holds in our
present situation, i.e. the compositions i) and ii) above result in Green remainder terms as asserted,
and the principal components satisfy the desired multiplicative identity. Note, moreover, that each
component of $M(\lambda)$ gives rise to an operator-valued symbol taking values in 
$\Dom^s_{\wedge,\min}\binom{A_{\wedge}}{T_{\wedge}}$.
\item Composition iii) is of the form ``$1 + \textup{Green}$'', and the proof of this follows similar to
the corresponding result for the boundaryless case, see Theorem 5.24 in \cite{GKM2}:

The composition behavior of parameter-dependent cone operators in Sobolev spaces implies
\begin{gather*}
\binom{A-\lambda}{T}{\mathcal B}_2(\lambda) = 1 + \tilde{M}(\lambda) + {\mathcal G}'(\lambda), \\
\tilde{M}(\lambda)\begin{pmatrix} u \\ v_1 \\ \vdots \\ v_K \end{pmatrix}(x) =
\frac{1}{2{\pi}}\!
\int\limits_{\Im\sigma = -m/2}
\int\limits_{\R_+}\Bigl(\frac{x}{x'}\Bigr)^{i\sigma} \tilde{m}(x,x',\sigma,\lambda)
\begin{pmatrix} u(x') \\ v_1(x') \\ \vdots \\ v_K(x') \end{pmatrix} \, \frac{dx'}{x'}\,d\sigma,
\end{gather*}
where ${\mathcal G}'(\lambda)$ is a matrix of generalized Green remainders, and
$\tilde{m}(x,x',\sigma,\lambda)$ equals
{\small
$$
-\left(\begin{smallmatrix} {x}^{-m} & 0 & \cdots & 0 \\ 0 & {x}^{-m_1} & \cdots & 0 \\ \vdots & \vdots & \ddots & \vdots \\
0 & 0 & \cdots & {x}^{-m_K} \end{smallmatrix}\right)
\omega\bigl(x[\lambda]^{1/m}\bigr)\hat{{\mathcal A}}_0(\sigma)h_0(\sigma)\omega\bigl(x'[\lambda]^{1/m}\bigr)
\left(\begin{smallmatrix} {x'}^{m} & 0 & \cdots & 0 \\ 0 & {x'}^{m_1} & \cdots & 0 \\ \vdots & \vdots & \ddots & \vdots \\
0 & 0 & \cdots & {x'}^{m_K} \end{smallmatrix}\right)
$$}
with $h_0(\sigma)$ from \eqref{SmoothMellRest}. The composition $\binom{A-\lambda}{T}M(\lambda)$
compensates the Mellin term $\tilde{M}(\lambda)$, and the remainder is therefore Green.
\end{itemize}
\end{proof}

\begin{remark}\label{Parametrix3links}
In the situation of Proposition \ref{Parametrix3}, we also have that
$$
\Biggl[\left(\begin{array}{cc} {\mathcal B}_3(\lambda) & 0 \\ \hline \Vsp
0 & 0 \end{array}\right) + \tilde{{\mathcal G}}(\lambda)\Biggr] \cdot
\Biggl[\left(\begin{array}{c|c} A-\lambda & 0 \\ \hline \Vsp
T & 0 \\ \hline \Vsp
0 & 0 \end{array}\right) + {\mathcal G}(\lambda)\Biggr] = 1 + \check{{\mathcal G}}(\lambda),
$$
where $\check{{\mathcal G}}(\lambda) = \bigl(\check{{\mathcal G}}_{i,j}(\lambda)\bigr)_{i,j=0,1}$, and
$\check{{\mathcal G}}_{0,0}(\lambda) \in {\mathcal B}^{-\infty,m}(\open{\overline{M}};\Lambda)$ is a regularizing para\-me\-ter-dependent
singular Green operator in Boutet de Monvel's calculus on $\open{\overline{M}}$. In addition, we have
the following properties (see also Lemma 5.20, Proposition 5.22, and Theorem 5.24 in \cite{GKM2}):

For all cut-off functions $\omega,\tilde{\omega} \in C_0^{\infty}([0,1))$ near zero the following
holds:
\begin{enumerate}[i)]
\item {\small
$$
\begin{pmatrix} 1-\omega & 0 \\ 0 & 0 \end{pmatrix}\check{{\mathcal G}}(\lambda),\; \check{{\mathcal G}}(\lambda)\begin{pmatrix} 1-\tilde{\omega} & 0 \\ 0 & 0 \end{pmatrix} \in
\bigcap\limits_{s,t > -\frac{1}{2}}\S\left(\Lambda,\K\left(
\begin{array}{c}\Dom^s_{\min}\binom{A}{T} \\ \oplus \\ \C^{N_-}\end{array},
\begin{array}{c}\Dom^t_{\min}\binom{A}{T} \\ \oplus \\ \C^{N_-}\end{array}
\right)\right).
$$}
\item $g(\lambda) = \begin{pmatrix} \omega & 0 \\ 0 & 1 \end{pmatrix}\check{{\mathcal G}}(\lambda)\begin{pmatrix} \tilde{\omega} & 0 \\ 0 & 1 \end{pmatrix}$
is a compact operator-valued symbol
$$
\begin{array}{c} {}_{\delta}\Dom^s_{\wedge,\min}\binom{A_{\wedge}}{T_{\wedge}} \\ \oplus \\ \C^{N_-} \end{array} \to
\begin{array}{c} {}_{\delta'}\Dom^t_{\wedge,\min}\binom{A_{\wedge}}{T_{\wedge}} \\ \oplus \\ \C^{N_-} \end{array}
$$
of order zero for all $\delta,\delta' \in \R$, $s,t > -\frac{1}{2}$.
\end{enumerate}
The reader surely noticed that in Definition \ref{Greenremainders} of generalized Green remainders
the case of operators whose domain is the $\Dom_{\min}$-scale was excluded. For the purposes of this
paper, we can consider the abovementioned properties as a definition for generalized Green remainders
of type $m$ (and order zero). However, if one is interested in parametrices for parabolic equations
with time-dependent coefficients on conic manifolds, this definition will not describe
in an appropriate way the structure of the symbol kernels.

For $\check{{\mathcal G}}_{\wedge}(\lambda):= g_{\wedge}(\lambda)$, where $g(\lambda)$ is the
operator-valued symbol in ii) above, we also have
$$
\Biggl[\left(\begin{array}{cc} {\mathcal B}_{3,\wedge}(\lambda) & 0 \\ \hline \Vsp
0 & 0 \end{array}\right) + \tilde{{\mathcal G}}_{\wedge}(\lambda)\Biggr] \cdot
\Biggl[\left(\begin{array}{c|c} A_{\wedge}-\lambda & 0 \\ \hline \Vsp
T_{\wedge} & 0 \\ \hline \Vsp
0 & 0 \end{array}\right) + {\mathcal G}_{\wedge}(\lambda)\Biggr] = 1 + \check{{\mathcal G}}_{\wedge}(\lambda)
$$
for $\lambda \in \Lambda\setminus\{0\}$.
\end{remark}

Observe, in particular, that the parametrix ${\mathcal B}_3(\lambda)$ is a Fredholm inverse of
the boundary value problem \eqref{BVPCone} by Proposition \ref{Parametrix3} and
Remark \ref{Parametrix3links}, and the principal part ${\mathcal B}_{3,\wedge}(\lambda)$
is a Fredholm inverse of the associated problem on the model cone $\begin{pmatrix} A_{\wedge} - \lambda \\ T_{\wedge} \end{pmatrix}$
for $\lambda \in \Lambda\setminus\{0\}$.

The following Theorem \ref{MainTheoremParametrix} deals with the final refinement of the parametrix,
and constitutes the main result as regards the parametrix construction of \eqref{BVPCone}.

\begin{theorem}\label{MainTheoremParametrix}
Assume that \eqref{BVPCone} is $c$-elliptic with parameter $\lambda \in \Lambda$, and assume that the
model boundary value problem
\begin{equation}\label{Voraussetzung}
{\mathcal A}_{\wedge}(\lambda) = \begin{pmatrix} A_{\wedge}-\lambda \\ T_{\wedge} \end{pmatrix} :
\Dom^s_{\wedge,\min}\binom{A_{\wedge}}{T_{\wedge}} \to
\begin{array}{c}
\K^{s,-m/2}(\overline{Y}^{\wedge};E) \\ \oplus \\
\bigoplus\limits_{j=1}^K\K^{s+m-m_j-1/2,m/2-m_j}((\partial\overline{Y})^{\wedge};F_j)
\end{array}
\end{equation}
is injective for some $s > -\frac{1}{2}$ and all $\lambda \in \Lambda\setminus\{0\}$. Recall that, by
$\kappa$-homogeneity, the injectivity needs to be required for $|\lambda| = 1$ only, and
the injectivity of \eqref{Voraussetzung} is equivalent to the injectivity of
\begin{equation}\label{Voraussetzung1}
A_{\wedge,T_{\wedge}} - \lambda : \Dom^s_{\wedge,\min}(A_{\wedge,T_{\wedge}}) \to
\K^{s,-m/2}(\overline{Y}^{\wedge};E)
\end{equation}
for all $\lambda \in \Lambda\setminus\{0\}$ (or $|\lambda| = 1$).
\begin{enumerate}[i)]
\item There exists a generalized Green remainder
${\mathcal K}_0(\lambda) : \C^{d''} \to x^{-m/2}H_b^s(\overline{M};E)$
of order $m_0=m$, where $d'' = -\Ind\binom{A}{T}_{\Dom^s_{\min}}$, such that
\begin{equation}\label{ExtraCondInvertierbar}
\left(\begin{array}{c|c}
A-\lambda & {\mathcal K}_0(\lambda) \\ \hline
\Vsp
\gamma B_1 & 0 \\
\vdots & \vdots \\
\gamma B_K & 0
\end{array}\right) :
\begin{array}{c} \Dom_{\min}^s\binom{A}{T} \\ \oplus \\ \C^{d''} \end{array} \to
\begin{array}{c} x^{-m/2}H_b^s(\overline{M};E) \\ \oplus \\ \bigoplus\limits_{j=1}^K x^{m/2-m_j}H_b^{s+m-m_j-1/2}(\overline{N};F_j) \end{array}
\end{equation}
is invertible for all $s > -\frac{1}{2}$ and $\lambda \in \Lambda$ with $|\lambda| > 0$ sufficiently large.

Moreover, there exists a matrix
{\small
$$
\left(\begin{array}{c|ccc}
{\mathcal G}_0(\lambda) & {\mathcal G}_1(\lambda) & \cdots & {\mathcal G}_K(\lambda) \\ \hline
\Vsp
{\mathcal T}_0(\lambda) & {\mathcal T}_1(\lambda) & \cdots & {\mathcal T}_K(\lambda)
\end{array}\right) :
\begin{array}{c} x^{-m/2}H_b^s(\overline{M};E) \\ \oplus \\ \bigoplus\limits_{j=1}^K x^{m/2-m_j}H_b^{s+m-m_j-1/2}(\overline{N};F_j) \end{array}
\to
\begin{array}{c} \Dom_{\min}^s\binom{A}{T} \\ \oplus \\ \C^{d''} \end{array}
$$}
of generalized Green remainders ${\mathcal G}_j(\lambda)$ and ${\mathcal T}_j(\lambda)$ of
orders $-m_j$ and type zero, such that the inverse of \eqref{ExtraCondInvertierbar} is of the form
\begin{equation}\label{ExtraCondInverse}
\left(\begin{array}{c}
{\mathcal B}_3(\lambda) + {\mathcal G}(\lambda) \\ \hline
\Vsp
{\mathcal T}(\lambda)
\end{array}\right) :
\left[\begin{array}{c} x^{-m/2}H_b^s(\overline{M};E) \\ \oplus \\ \bigoplus\limits_{j=1}^K x^{m/2-m_j}H_b^{s+m-m_j-1/2}(\overline{N};F_j) \end{array}\right]
\to
\begin{array}{c} \Dom_{\min}^s\binom{A}{T} \\ \oplus \\ \C^{d''} \end{array},
\end{equation}
where ${\mathcal B}_3(\lambda)$ is the parametrix from \eqref{Parametrix3Definition}, and
\begin{align*}
{\mathcal G}(\lambda) &= \begin{pmatrix} {\mathcal G}_0(\lambda) & {\mathcal G}_1(\lambda) & \cdots & {\mathcal G}_K(\lambda) \end{pmatrix}, \\
{\mathcal T}(\lambda) &= \begin{pmatrix} {\mathcal T}_0(\lambda) & {\mathcal T}_1(\lambda) & \cdots & {\mathcal T}_K(\lambda) \end{pmatrix}.
\end{align*}
In particular, the boundary value problem
\begin{equation}\label{BVPaufMin}
{\mathcal A}(\lambda) = \left(\begin{array}{c}
A-\lambda \\ \hline
\Vsp
\gamma B_1 \\
\vdots \\
\gamma B_K
\end{array}\right) :
\Dom_{\min}^s\binom{A}{T} \to
\begin{array}{c} x^{-m/2}H_b^s(\overline{M};E) \\ \oplus \\ \bigoplus\limits_{j=1}^K x^{m/2-m_j}H_b^{s+m-m_j-1/2}(\overline{N};F_j) \end{array}
\end{equation}
is injective for large $\lambda \in \Lambda$, and the parametrix ${\mathcal B}(\lambda):= {\mathcal B}_3(\lambda) + {\mathcal G}(\lambda)$
is a left inverse.
\item Let
$$
\Pi(\lambda):= 1 - {\mathcal A}(\lambda){\mathcal B}(\lambda) =
\left(\begin{array}{c|ccc} \Pi_{0,0}(\lambda) & \Pi_{0,1}(\lambda) & \cdots & \Pi_{0,K}(\lambda) \\ \hline
\Vsp
\Pi_{1,0}(\lambda) & \Pi_{1,1}(\lambda) & \cdots & \Pi_{1,K}(\lambda) \\
\vdots & \vdots & \ddots & \vdots \\
\Pi_{K,0}(\lambda) & \Pi_{K,1}(\lambda) & \cdots & \Pi_{K,K}(\lambda)
\end{array}\right).
$$
Then $\Pi_{i,j}(\lambda)$ is a generalized Green remainder of order $m_i-m_j$ and
type zero, and $\Pi(\lambda)$ is for large $\lambda$ a finite-dimensional projection onto a complement of
the range of \eqref{BVPaufMin}. Whenever
$$
{\mathcal A}_{\Dom}(\lambda_0) = \left(\begin{array}{c}
A-\lambda_0 \\ \hline
\Vsp
\gamma B_1 \\
\vdots \\
\gamma B_K
\end{array}\right) :
\Dom^s\binom{A}{T} \to
\begin{array}{c} x^{-m/2}H_b^s(\overline{M};E) \\ \oplus \\ \bigoplus\limits_{j=1}^K x^{m/2-m_j}H_b^{s+m-m_j-1/2}(\overline{N};F_j) \end{array}
$$
is invertible for some $\lambda_0 \in \Lambda$ and some domain
$\Dom^s_{\min}\binom{A}{T} \subset \Dom^s\binom{A}{T} \subset \Dom^s_{\max}\binom{A}{T}$, the inverse
${\mathcal A}_{\Dom}(\lambda_0)^{-1}$ can be written in the form
${\mathcal B}(\lambda_0) + {\mathcal A}_{\Dom}(\lambda_0)^{-1}\Pi(\lambda_0)$.
\item Let ${\mathcal B}_T(\lambda) : x^{-m/2}H_b^s(\overline{M};E) \to \Dom^s_{\min}\binom{A}{T}$
be the interior part of the parametrix ${\mathcal B}(\lambda)$. Then, for large $|\lambda| > 0$,
$$
{\mathcal B}_T(\lambda) : x^{-m/2}H_b^s(\overline{M};E) \to \Dom^s_{\min}(A_T) = \Dom^s_{\min}\binom{A}{T}\cap\ker T,
$$
and ${\mathcal B}_T(\lambda)$ is a left inverse of
\begin{equation}\label{Unboundedmin}
A_T-\lambda : \Dom^s_{\min}(A_T) \to x^{-m/2}H_b^s(\overline{M};E).
\end{equation}

The operator $\Pi_T(\lambda):= \Pi_{0,0}(\lambda) = 1-(A_T-\lambda){\mathcal B}_T(\lambda)$ is a generalized Green remainder
of order and type zero, and for large $\lambda$ a (finite-dimensional) projection onto a complement of
the range of \eqref{Unboundedmin}. Whenever
$$
A_T-\lambda_0 : \Dom^s(A_T) \to x^{-m/2}H_b^s(\overline{M};E)
$$
is invertible for some $\lambda_0 \in \Lambda$ with $|\lambda_0| > 0$ sufficiently large
and some domain $\Dom^s_{\min}(A_T) \subset
\Dom^s(A_T) \subset \Dom^s_{\max}(A_T)$, the resolvent can be written as
$$
(A_T-\lambda_0)^{-1} = B_T(\lambda_0) + (A_T-\lambda_0)^{-1}\Pi_T(\lambda_0).
$$
\item The principal component
\begin{equation}\label{PrKomp}
\begin{pmatrix} A_{\wedge}-\lambda & {\mathcal K}_{0,\wedge}(\lambda) \\
T_{\wedge} & 0 \end{pmatrix} : \begin{array}{c}
\Dom^s_{\wedge,\min}\binom{A_{\wedge}}{T_{\wedge}} \\ \oplus \\ \C^{d''} \end{array}
\to
\begin{array}{c} \K^{s,-m/2}(\overline{Y}^{\wedge};E) \\ \oplus \\ \bigoplus\limits_{j=1}^K\K^{s+m-m_j-1/2,m/2-m_j}((\partial\overline{Y})^{\wedge};F_j) \end{array}
\end{equation}
of \eqref{ExtraCondInvertierbar} is invertible for all $\lambda \in \Lambda\setminus\{0\}$, and
the principal component
{\small
\begin{equation}\label{PrKompInv}
\left(\begin{array}{c} {\mathcal B}_{3,\wedge}(\lambda) + {\mathcal G}_{\wedge}(\lambda) \\ \hline
\Vsp
{\mathcal T}_{\wedge}(\lambda) \end{array}\right) :
\left[\begin{array}{c} \K^{s,-m/2}(\overline{Y}^{\wedge};E) \\ \oplus \\ \bigoplus\limits_{j=1}^K\K^{s+m-m_j-1/2,m/2-m_j}((\partial\overline{Y})^{\wedge};F_j) \end{array}\right]
\to \begin{array}{c} \Dom^s_{\wedge,\min}\binom{A_{\wedge}}{T_{\wedge}} \\ \oplus \\ \C^{d''} \end{array}
\end{equation}}
of \eqref{ExtraCondInverse} is the inverse.
\end{enumerate}
\end{theorem}
\begin{proof}
Let $X = \Lambda\cap S^1$, and consider the operator family $A_{\wedge,T_{\wedge}} - \lambda$
from \eqref{Voraussetzung1} as a smooth Fredholm function on $X$. By well known results about
Fredholm families on compact spaces and a density argument (see also the appendix of \cite{GKM2}),
there exists a function
$$
{\mathcal K}_{0,\wedge}(\lambda) \in C^{\infty}(X)\otimes (\C^{d''})^*\otimes C_0^{\infty}(\open{\overline{Y}}^{\wedge};E)
$$
such that
\begin{equation}\label{aufgefuellt}
\begin{pmatrix} A_{\wedge}-\lambda & {\mathcal K}_{0,\wedge}(\lambda) \end{pmatrix} :
\begin{array}{c} \Dom^s_{\wedge,\min}(A_{\wedge,T_{\wedge}}) \\ \oplus \\ \C^{d''} \end{array} \to
\K^{s,-m/2}(\overline{Y}^{\wedge};E)
\end{equation}
is invertible for $\lambda \in X$, and so is the extension of \eqref{aufgefuellt} to $\Lambda\setminus\{0\}$
by (twisted) homogeneity of degree $m$ (we will use the same notation ${\mathcal K}_{0,\wedge}(\lambda)$).
A simple calculation now shows that also \eqref{PrKomp} is invertible for $\lambda \in \Lambda\setminus\{0\}$
for this choice of ${\mathcal K}_{0,\wedge}(\lambda)$.

As $\Ind {\mathcal A}_{\wedge}(\lambda) = -\Ind {\mathcal B}_{3,\wedge}(\lambda)$, the same abstract
results about Fredholm families on compact spaces as applied before and extension by (twisted)
homogeneity now imply the existence of a matrix ${\mathcal C}_{\wedge}(\lambda)$
of principal Green symbols of type zero and suitable $N_-, N_+ \in \N_0$,
$$
N_+ - N_- = \Ind {\mathcal B}_{3,\wedge}(\lambda) = -d'',
$$
such that
{\small
\begin{equation}\label{auff}
\begin{pmatrix} {\mathcal B}_{3,\wedge}(\lambda) & 0 & 0 \\ 0 & 0 & 0 \end{pmatrix} + {\mathcal C}_{\wedge}(\lambda) :
\begin{array}{c} \K^{s,-m/2}(\overline{Y}^{\wedge};E) \\ \oplus \\ \bigoplus\limits_{j=1}^K\K^{s+m-m_j-1/2,m/2-m_j}((\partial\overline{Y})^{\wedge};F_j) \\ \oplus \\ \C^{N_-} \end{array}
\to
\begin{array}{c} \Dom^s_{\wedge,\min}\binom{A_{\wedge}}{T_{\wedge}} \\ \oplus \\ \C^{N_+} \end{array}
\end{equation}}
is invertible for $\lambda \in \Lambda\setminus\{0\}$.
By possibly enlarging $N_-$ and $N_+$ and the matrix ${\mathcal C}_{\wedge}(\lambda)$, we may assume
that $N_+ \geq d''$, and by possibly augmenting the matrix \eqref{PrKomp} by an invertible lower right
corner (if $N_+ > d''$), we can multiply \eqref{PrKomp} and \eqref{auff}. The product is an
invertible matrix of the form $1 + {\mathcal C}_{\wedge}'(\lambda)$, where ${\mathcal C}'_{\wedge}(\lambda)$ is
a matrix of principal Green symbols. As
$$
\bigl(1 + {\mathcal C}_{\wedge}'(\lambda)\bigr)^{-1} = 1 + {\mathcal C}_{\wedge}''(\lambda)
$$
with a matrix ${\mathcal C}_{\wedge}''(\lambda)$ of principal Green symbols, we conclude that
the inverse of \eqref{PrKomp} is indeed of the form \eqref{PrKompInv} for suitable matrices ${\mathcal G}_{\wedge}(\lambda)$
and ${\mathcal T}_{\wedge}(\lambda)$ of principal Green symbols of the asserted order and type zero.
Here we used the fact that the matrices of the form ``$1 + \textup{Green}$'' are spectrally invariant,
see the proof of Lemma \ref{EinsplusGreen}.

Choose a cut-off function $\omega \in C_0^{\infty}([0,1))$ near zero, and an excision function
$\chi \in C^{\infty}(\C)$, i.e. $\chi \equiv 0$ near zero, and $\chi \equiv 1$ near infinity.
Define ${\mathcal K}_0(\lambda) = \omega\chi(\lambda){\mathcal K}_{0,\wedge}(\lambda)$, and
${\mathcal G}'(\lambda) = \omega\chi(\lambda){\mathcal G}_{\wedge}(\lambda)\omega$,
${\mathcal T}'(\lambda) = \chi(\lambda){\mathcal T}_{\wedge}(\lambda)\omega$.
Then ${\mathcal K}_0(\lambda)$, ${\mathcal G}'(\lambda)$, and ${\mathcal T}'(\lambda)$ are
matrices of generalized Green remainders, and
$$
\left(\begin{array}{c|c}
A-\lambda & {\mathcal K}_0(\lambda) \\ \hline
\Vsp
T & 0
\end{array}\right)
\left(\begin{array}{c}
{\mathcal B}_3(\lambda) + {\mathcal G}'(\lambda) \\ \hline
\Vsp
{\mathcal T}'(\lambda)
\end{array}\right)
= 1 + {\mathcal G}''(\lambda),
$$
where ${\mathcal G}''(\lambda) = \bigl({\mathcal G}''_{i,j}(\lambda)\bigr)_{i,j=0,\ldots,K}$ is
a matrix of generalized Green remainders of type zero, and ${\mathcal G}''_{i,j}(\lambda)$ has
order $m_i-m_j$ (where $m_0 = m$). Moreover, by construction we have the situation of Lemma \ref{EinsplusGreen}
for $1 + {\mathcal G}''(\lambda)$, and thus $\begin{pmatrix} A-\lambda & {\mathcal K}_0(\lambda) \\ T & 0 \end{pmatrix}$ is invertible
from the right for $\lambda \in \Lambda$ with $|\lambda| > 0$ sufficiently large, and the right inverse is of the form
\eqref{ExtraCondInverse}. Hence both i) and iv) will be proved if we show that $\begin{pmatrix} A-\lambda & {\mathcal K}_0(\lambda) \\ T & 0 \end{pmatrix}$
is also invertible from the left.

To this end, note that
$$
\left(\begin{array}{c}
{\mathcal B}_3(\lambda) + {\mathcal G}'(\lambda) \\ \hline
\Vsp
{\mathcal T}'(\lambda)
\end{array}\right)
\left(\begin{array}{c|c}
A-\lambda & {\mathcal K}_0(\lambda) \\ \hline
\Vsp
T & 0
\end{array}\right)
= 1 + \check{{\mathcal G}}(\lambda)
$$
with an operator $\check{{\mathcal G}}(\lambda)$ which satisfies the conditions
of Remark \ref{Parametrix3links}. Moreover, by construction $\check{{\mathcal G}}_{\wedge}(\lambda) = 0$,
where $\check{{\mathcal G}}_{\wedge}(\lambda)$ is the principal part of $\check{{\mathcal G}}(\lambda)$.
Thus
$$
\Bigl(\sum\limits_{j=0}^N(-1)^j\check{{\mathcal G}}(\lambda)^j\Bigr)\bigl(1+\check{{\mathcal G}}(\lambda)\bigr) = 1 + \check{{\mathcal G}}_{N+1}(\lambda),
$$
and for $N > 0$ sufficiently large the operator norm of $\check{{\mathcal G}}_{N+1}(\lambda)$ in
$\L\Biggl(\begin{array}{c} \Dom^s_{\min}\binom{A}{T} \\ \oplus \\ \C^{d''} \end{array}\Biggr)$
is tending to zero as $|\lambda| \to \infty$. This shows that \eqref{ExtraCondInvertierbar} is invertible
from the left and completes the proof of i) and iv). Note that ii) and iii) follow immediately from i)
by simple algebraic calculations.
\end{proof}


\section{Resolvents}
\label{sec-Resolvents}

The final section is devoted to the main theorem of this article:

\begin{theorem}\label{Maintheorem}
Let \eqref{BVPCone} be $c$-elliptic with parameter in the closed sector $\Lambda \subset \C$, and consider
the unbounded operator $A$ in $x^{-m/2}L^2_b(\overline{M};E)$ under the boundary condition $Tu = 0$
on some intermediate domain $\Dom_{\min}(A_T) \subset \Dom(A_T) \subset \Dom_{\max}(A_T)$.

Let $\Dom_{\wedge}(A_{\wedge,T_{\wedge}}) = \theta(\Dom(A_T))$ be the associated domain for the
model operator $A_{\wedge}$ under the boundary condition $T_{\wedge}u = 0$ according to
Proposition \ref{SingularFunctionManifold}, and assume that $\Lambda$ is a sector
of minimal growth for $A_{\wedge}$ with this domain, i.e.
$$
A_{\wedge} - \lambda : \Dom_{\wedge}(A_{\wedge,T_{\wedge}}) \to x^{-m/2}L_b^2(\overline{Y}^{\wedge};E)
$$
is invertible for $\lambda \in \Lambda$ with $|\lambda| > 0$ sufficiently large, and the
resolvent satisfies the norm estimate
$$
\|\bigl(A_{\wedge,\Dom_{\wedge}} - \lambda\bigr)^{-1}\|_{\L(x^{-m/2}L_b^2)} = O(|\lambda|^{-1})
$$
as $|\lambda| \to \infty$.

Then $\Lambda$ is a sector of minimal growth for the operator $A$ in $x^{-m/2}L_b^2(\overline{M};E)$
with domain $\Dom(A_T)$, and for large $\lambda \in \Lambda$ the resolvent can be written in the form
\begin{equation}\label{InverseDarstellung}
\bigl(A_{\Dom}-\lambda\bigr)^{-1} = {\mathcal B}_T(\lambda) + \bigl(A_{\Dom}-\lambda\bigr)^{-1}\Pi_T(\lambda)
\end{equation}
with the parametrix ${\mathcal B}_T(\lambda)$ and projection $\Pi_T(\lambda)$ onto a complement
of the range of $A_{\min}-\lambda$ from Theorem \ref{MainTheoremParametrix}.
\end{theorem}

\noindent
The resolvent condition on $A_{\wedge}$ from Theorem \ref{Maintheorem} is an analogue of the
Shapiro-Lopatinsky condition and is associated with the ``singular boundary'' $\overline{Y}$ of
$\overline{M}$ (see Proposition \ref{WedgeSectorMinimalGrowth} for a discussion of this assumption).

With the preparations from the previous sections, we are able to follow the same idea as in the
boundaryless case in \cite{GKM2}.

\medskip

Let
$$
\theta : \tilde\Sing_{\max} = \bigoplus\limits_{\sigma_0 \in \Sigma}\tilde\Sing_{\sigma_0} \to
\bigoplus\limits_{\sigma_0 \in \Sigma}\tilde\Sing_{\wedge,\sigma_0} = \tilde\Sing_{\wedge,\max}
$$
be the isomorphism of the spaces of singular functions $\tilde\Sing_{\max} \cong \Dom_{\max}/\Dom_{\min}$
and $\tilde\Sing_{\wedge,\max} \cong \Dom_{\wedge,\max}/\Dom_{\wedge,\min}$ that was
constructed in Section \ref{sec-AssociatedDomains}. Recall that $\Sigma$ is
the part of the boundary spectrum of $\binom{A}{T}$ in $\{\sigma\in\C\st -m/2 < \Im(\sigma) < m/2\}$,
and for $\sigma_0 \in \Sigma$ let $N(\sigma_0) \in \N_0$ be the largest integer such that
$\Im(\sigma_0) -N(\sigma_0) > -m/2$.

The normalized dilation group $\kappa_{\varrho}$ respects the space $\tilde\Sing_{\wedge,\max}$, i.e.
$\kappa_{\varrho} : \tilde\Sing_{\wedge,\max} \to \tilde\Sing_{\wedge,\max}$ for $\varrho > 0$.
Consequently, we can define a group action $\tilde\kappa_{\varrho}$ on $\tilde\Sing_{\max}$ via
$$
\tilde\kappa_{\varrho} = \theta^{-1}\kappa_{\varrho}\theta : \tilde\Sing_{\max} \to \tilde\Sing_{\max}.
$$
We may write $\tilde\kappa_{\varrho} = \kappa_{\varrho}L_{\varrho}$, where
$$
L_{\varrho} = \kappa_{\varrho}^{-1}\theta^{-1}\kappa_{\varrho}\theta : \tilde\Sing_{\max} \to
C^{\infty}(\open{\overline{Y}}^{\wedge};E)
$$
is the direct sum of the operators $L_{\varrho}|_{\tilde\Sing_{\sigma_0}}$ which act as follows:

For $\tilde{u} \in \tilde{\Sing}_{\sigma_0}$ we have
\begin{equation}\label{AltLRho}
L_{\varrho}\tilde{u} = \sum\limits_{\vartheta=0}^{N(\sigma_0)}\varrho^{-\vartheta}\e_{\sigma_0,\vartheta}(\varrho)(\theta\tilde{u}),
\end{equation}
where $\e_{\sigma_0,\vartheta}(\varrho)$ is defined as
\begin{equation*}
\e_{\sigma_0,\vartheta}(\varrho) = \varrho^{\vartheta}\kappa_{\varrho}^{-1}
\e_{\sigma_0,\vartheta}\kappa_{\varrho}:
\tilde{\Sing}_{\wedge,\sigma_0}\to C^{\infty}(\open{\overline{Y}}^{\wedge};E)
\end{equation*}
with the operators $\e_{\sigma_0,\vartheta}$ from Section \ref{sec-AssociatedDomains}.
In particular, $\e_{\sigma_0,0}(\varrho)(\tilde u) = \tilde u$ for
all $\varrho\in\R_+$ and $\tilde u\in\tilde\Sing_{\wedge,\sigma_0}$.

\begin{lemma}\label{Ltheta} \
\begin{enumerate}[i)]
\item For every $\psi \in \tilde{\Sing}_{\wedge,\sigma_0}$ and every 
$\vartheta \in \{0,\ldots,N(\sigma_0)\}$ there exists a polynomial 
$q_\vartheta(y,\log x,\log \varrho)$ in $(\log x,\log \varrho)$ with
coefficients in $C^{\infty}(\overline{Y};E)$ such that
\begin{equation}\label{esigmavarrho}
\e_{\sigma_0,\vartheta}(\varrho)(\psi) = q_\vartheta(y,\log x,\log \varrho)
x^{i(\sigma_0 - i\vartheta)},
\end{equation}
and the degree of $q_\vartheta$ with respect to $(\log x,\log \varrho)$ is 
bounded by some $\mu \in \N_0$ which is independent of $\sigma_0\in\Sigma$, 
$\psi \in \tilde{\Sing}_{\wedge,\sigma_0}$, and
$\vartheta \in \{0,\ldots,N(\sigma_0)\}$.

\item Let $\omega \in C_0^{\infty}(\overline{\R}_+)$ be any cut-off
function near the origin, i.e.,  $\omega = 1$ near zero and $\omega =
0$ near infinity. Then the operator family
$$
\omega\bigl(L_{\varrho} - \theta\bigr) : \tilde{\Sing}_{\max} \to \K^{\infty,-m/2}(\overline{Y}^{\wedge};E)
$$
satisfies for every $s \in \R$ the norm estimate
$$
\|\omega\bigl(L_{\varrho} - \theta\bigr)\|_{\L(\tilde{\Sing}_{\max},
\K^{s,-m/2})} = O(\varrho^{-1}\log^{\mu}\varrho)
\;\text{ as } \varrho \to \infty, 
$$
where $\mu \in \N_0$ is the bound for the degrees of the polynomials
$q_\vartheta$ in i).
\end{enumerate}
\end{lemma}
\begin{proof}
The proof is literally the same as in Lemma 6.18 from \cite{GKM2}.
\end{proof}

\begin{lemma}\label{kappaSchlangelift}
Fix a cut-off function $\omega \in C_0^{\infty}([0,1))$ near $0$. 
For $\varrho > 1$ consider the operator family
\begin{equation*}
\tilde{K}(\varrho)= \omega_{\varrho} \tilde{\kappa}_{\varrho} :
\tilde\Sing_{\max} \to \Dom^{\infty}_{\max}\binom{A}{T}
= \bigcap_{t > -\frac{1}{2}} \Dom^t_{\max}\binom{A}{T},
\end{equation*}
where $\omega_{\varrho}(x)= \omega(\varrho x)$. If $q: \Dom_{\max}\binom{A}{T} \to 
\tilde\Sing_{\max}$ is the canonical projection, then 
\begin{equation*}
q{\circ}\tilde{K}(\varrho) = \tilde{\kappa}_{\varrho},
\end{equation*}
and we have the following norm estimates as $\varrho \to \infty$:
\begin{align} \label{Kestimate1}
&\|\tilde{K}(\varrho)\|_{\L(\tilde\Sing_{\max},x^{-m/2}L^2_b)} = O(1), \\ \label{Kestimate2}
&\|\kappa_{\varrho}^{-1}A\tilde{K}(\varrho)\|_{\L(\tilde\Sing_{\max},x^{-m/2}L^2_b)} = O(\varrho^m), \\ \label{Kestimate3}
&\|\kappa_{\varrho}^{-1}\gamma B_\ell\tilde{K}(\varrho)\|_{\L(\tilde\Sing_{\max},\K^{m-m_\ell-1/2,m/2-m_\ell})} = O(\varrho^{m_\ell}), \quad \ell=1,\ldots,K.
\end{align}
Note that $\tilde{K}(\varrho)\tilde{u}$ is supported in $(0,\varrho^{-1}]\times\overline{Y} \subset U_{\overline{Y}}$ for
all $\tilde{u} \in \tilde\Sing_{\max}$, and thus it makes sense to apply the group action $\kappa^{-1}_{\varrho}$
in the estimates \eqref{Kestimate2} and \eqref{Kestimate3}.
\end{lemma}
\begin{proof}
That $\tilde{K}(\varrho)$ is a lift of $\tilde{\kappa}_{\varrho}$ to 
$\Dom_{\max}^{\infty}\binom{A}{T}$ is evident from the definition.
In order to show the norm estimates, it is sufficient 
to consider for each $\sigma_0 \in \Sigma$ the restriction
$$
\tilde{K}_{\sigma_0}(\varrho)= \tilde{K}(\varrho)|_{\tilde\Sing_{\sigma_0}}: 
\tilde\Sing_{\sigma_0} \to \Dom_{\max}^{\infty}\binom{A}{T}
$$
and prove the estimates for this operator.
Recall that $\tilde\kappa_\varrho=\kappa_\varrho L_\varrho$ so that for 
$\tilde u\in \tilde\Sing_{\sigma_0}$ we have $\tilde K_{\sigma_0}(\varrho)
\tilde u = \kappa_\varrho(\omega L_\varrho\tilde u)$.

The norm estimates \eqref{Kestimate1} and \eqref{Kestimate2} follow in the
same way as the corresponding assertion in the boundaryless case, see
Lemma 6.20 in \cite{GKM2}. The same method of proof also gives \eqref{Kestimate3};
for sake of completeness, we give a proof of this estimate below,
i.e. we prove that there exists a constant $C > 0$, independent of 
$\tilde u\in \tilde\Sing_{\sigma_0}$ and $\varrho \geq 1$, such that
\begin{equation*}
\|\kappa_{\varrho}^{-1}\gamma B_\ell(\kappa_{\varrho}(\omega L_{\varrho}\tilde u))\|_{\K^{m-m_\ell-1/2,m/2-m_\ell}} \leq C\varrho^{m_\ell} \|\omega \tilde u\|_{\Dom_{\max}}.
\end{equation*}
To this end we split $B_\ell$ near $\overline{Y}$ as in \eqref{exp2}, i.e.
$$
B_\ell \equiv x^{-m_\ell}\sum\limits_{k=0}^{m-1}B_{\ell,k}x^k + \tilde{B}_{\ell,m}
$$
with totally characteristic operators $B_{\ell,k} \in \Diff_b^{m_\ell}(\overline{Y}^{\wedge};E,F_\ell)$ with
coefficients independent of $x$, and $\tilde{B}_{\ell,m} \in x^{m-m_\ell}\Diff_b^{m_\ell}(\overline{Y}^{\wedge};E,F_\ell)$.
As we are working exclusively near $\overline{Y}$, we may without loss of generality assume that
the coefficients of $\tilde{B}_{\ell,m}$ vanish near infinity.

By \eqref{AltLRho} we obtain
{\small\begin{equation}\label{PhiAPhiL}
\begin{aligned}
\kappa_{\varrho}^{-1}\gamma B_{\ell} & (\kappa_{\varrho}(\omega L_{\varrho}\tilde u)) \\
&= \kappa_{\varrho}^{-1}\Bigl(x^{-m_{\ell}}\sum_{k=0}^{m-1}\gamma B_{\ell,k}x^k\Bigr)\kappa_{\varrho}
\bigl({\omega}L_{\varrho}\tilde u\bigr) +
\kappa_{\varrho}^{-1}\gamma \tilde{B}_{\ell,m}\kappa_{\varrho}\bigl({\omega}L_{\varrho}\tilde u\bigr) \\
&= \varrho^{m_{\ell}}
\Bigl(x^{-m_{\ell}}\sum_{k=0}^{m-1}\varrho^{-k}\gamma B_{\ell,k}x^k\Bigr)
\Bigl({\omega}\sum_{j=0}^{N(\sigma_0)}\varrho^{-j}
\e_{\sigma_0,j}(\varrho)(\theta \tilde u)\Bigr) +
\kappa_{\varrho}^{-1}\gamma \tilde{B}_{\ell,m}\kappa_{\varrho}\bigl({\omega}L_{\varrho}\tilde u\bigr) \\
&= \sum_{\vartheta=0}^{2m-2} \varrho^{m_{\ell}-\vartheta}\Bigl(
x^{-m_{\ell}}\!\!\!\sum_{\substack{k+j=\vartheta \\ 0 \leq k,j \leq m-1}}\!\!\!
\bigl(\gamma B_{\ell,k}x^k\bigr)\bigl({\omega}\e_{\sigma_0,j}(\varrho)(\theta \tilde u)\bigr)
\Bigr) +\kappa_{\varrho}^{-1}\gamma \tilde{B}_{\ell,m}\kappa_{\varrho}\bigl({\omega}L_{\varrho}\tilde u\bigr)
\end{aligned}
\end{equation}}
with the convention that $\e_{\sigma_0,j}(\varrho)= 0$ for $j > N(\sigma_0)$.

For every $\vartheta\in\{0,\dots,2m-2\}$ we consider the family of linear maps 
\begin{equation}\label{KreisSumme}
\begin{aligned}
\tilde u\mapsto 
x^{-m_{\ell}}\!\!\sum_{\substack{k+j=\vartheta \\ 0 \leq k,j \leq m-1}}\!\!
\bigl(\gamma B_{\ell,k}x^k\bigr)&\bigl({\omega}\e_{\sigma_0,j}(\varrho)(\theta \tilde u)\bigr) \\
&: \tilde\Sing_{\sigma_0}\to \K^{m-m_{\ell}-1/2,m/2-m_{\ell}}((\partial\overline{Y})^{\wedge};F_{\ell}).
\end{aligned}
\end{equation}
We will prove that \eqref{KreisSumme} is well-defined, i.e., every $\tilde u \in \tilde{\Sing}_{\sigma_0}$ is 
indeed mapped into $\K^{m-m_{\ell}-1/2,m/2-m_{\ell}}$, and that the norms are 
bounded by a constant times $\log^\mu \varrho$ as $\varrho\to\infty$ with
$\mu$ as in Lemma~\ref{Ltheta}. 
Thus for every $\vartheta \in \{0,\ldots,2m-2\}$ we have
\begin{multline*}
\Bigl\|\varrho^{m_{\ell}-\vartheta}\Bigl(
x^{-m_{\ell}}\!\!\!\sum_{\substack{k+j=\vartheta \\ 0 \leq k,j \leq m-1}}
\!\!\bigl(\gamma B_{\ell,k}x^k\bigr)\bigl({\omega}\e_{\sigma_0,j}(\varrho)(\theta \tilde u)
\bigr)\Bigr)\Bigr\|_{\K^{m-m_{\ell}-1/2,m/2-m_{\ell}}} \\
\leq \textup{const}\cdot\bigl(\varrho^{m_{\ell}-\vartheta} \log^\mu
\varrho\bigr) \|\omega \tilde u\|_{\Dom_{\max}},
\end{multline*}
while for $\vartheta=0$, 
\begin{equation}\label{Theta0Estimate}
\varrho^{m_\ell}  x^{-m_{\ell}}\gamma B_{\ell,0}\omega 
\e_{\sigma_0,0}(\varrho)(\theta \tilde u) =
\varrho^{m_\ell} \gamma_{\wedge}B_{\ell,\wedge}\omega (\theta \tilde u),
\end{equation}
so for this term we have a norm estimate without $\log$.

Let $\tilde\omega\in C_0^\infty(\overline\R_+)$ be a cut-off function near 
$0$ with $\omega\prec\tilde\omega$. Then there exist suitable 
$\vp,\tilde \vp\in C_0^\infty(\R_+)$ such that for all 
$\tilde u\in\tilde\Sing_{\sigma_0}$,
{\small
\begin{multline}\label{SumSplitting}
x^{-m_{\ell}}\!\!\!\sum_{\substack{k+j=\vartheta \\ 0 \leq k,j \leq m-1}}\!\!
\big(\gamma B_{\ell,k}x^k\big)\big({\omega}\e_{\sigma_0,j}(\varrho)(\theta \tilde u)\big)\\
= \tilde\omega x^{-m_{\ell}}\!\!\!\sum_{\substack{k+j=\vartheta \\ 0 \leq k,j 
\leq m-1}}\!\! \bigl(\gamma B_{\ell,k}x^k\bigr)\e_{\sigma_0,j}(\varrho)(\theta \tilde u) +
\tilde\vp x^{-m_{\ell}}\!\!\!\sum_{\substack{k+j=\vartheta \\ 0 \leq k,j \leq m-1}}
\!\!  \bigl(\gamma B_{\ell,k}x^k\bigr)\vp\e_{\sigma_0,j}(\varrho)(\theta \tilde u).
\end{multline}}

According to Lemma \ref{Ltheta} the second sum in \eqref{SumSplitting} is a
polynomial in $\log \varrho$ of degree at most $\mu$ with coefficients in
$C_0^{\infty}(\partial\open{\overline{Y}}^{\wedge};F_{\ell})$. As both $\kappa_{\varrho}^{-1}\gamma B_{\ell}(\kappa_{\varrho}(\omega L_{\varrho}\tilde u))$
and $\kappa_{\varrho}^{-1}\gamma\tilde{B}_{\ell,m}(\kappa_{\varrho}({\omega}L_{\varrho}\tilde u))$ belong to
$\K^{m-m_{\ell}-1/2,m/2-m_{\ell}}$, we get from the equations \eqref{PhiAPhiL} and
\eqref{SumSplitting} that necessarily
$$
x^{-m_{\ell}}\!\!\!\sum_{\substack{k+j=\vartheta \\ 0 \leq k,j \leq m-1}}\!\!
\big(\gamma B_{\ell,k}x^k\big)\big({\omega}\e_{\sigma_0,j}(\varrho)(\theta \tilde u)\big)
\in \K^{m-m_{\ell}-1/2,m/2-m_{\ell}}
$$
for all $\varrho \in \R_+$ and all $\tilde{u} \in \tilde{\Sing}_{\sigma_0}$, and, moreover, that
$$
\tilde{\omega}x^{-m_{\ell}}\!\!\!\sum_{\substack{k+j=\vartheta \\ 0 \leq k,j \leq m-1}}\!\!
\big(\gamma B_{\ell,k}x^k\big)\e_{\sigma_0,j}(\varrho)(\theta \tilde u) = 0
$$
for $\Im(\sigma_0)-\vartheta \geq -m/2$ because these functions are of the form
$$
\tilde{\omega}\Bigl(\sum\limits_{\nu}c_{\sigma_0-i(\vartheta-m_{\ell}),\nu}(y')\log^{\nu} x\Bigr)x^{i(\sigma_0-i(\vartheta-m_{\ell}))}.
$$
For $\Im(\sigma_0)-\vartheta < -m/2$ every single summand
$\tilde{\omega}x^{-m_{\ell}}\big(\gamma B_{\ell,k}x^k\big)\e_{\sigma_0,j}(\varrho)(\theta \tilde u)$ belongs to the space
$\K^{m-m_{\ell}-1/2,m/2-m_{\ell}}$, and by Lemma \ref{Ltheta} is a polynomial in $\log \varrho$ of degree
at most $\mu$ with coefficients in $\K^{m-m_{\ell}-1/2,m/2-m_{\ell}}$.

Summing up, we have shown that for every $\tilde{u} \in \tilde{\Sing}_{\sigma_0}$ the function
$$
x^{-m_{\ell}}\!\!\!\sum_{\substack{k+j=\vartheta \\ 0 \leq k,j \leq m-1}}\!\!
\big(\gamma B_{\ell,k}x^k\big)\big({\omega}\e_{\sigma_0,j}(\varrho)(\theta \tilde u)\big)
$$
is a polynomial in $\log\varrho$ of degree at most $\mu$ with coefficients in 
$\K^{m-m_{\ell}-1/2,m/2-m_{\ell}}$, and from the Banach-Steinhaus theorem we now obtain the desired norm
estimates for the family of maps \eqref{KreisSumme}.

On the other hand, 
\begin{gather*}
\|\kappa_{\varrho}^{-1}\gamma\tilde{B}_{\ell,m}\kappa_{\varrho}\bigl({\omega}L_{\varrho}
\tilde u\bigr)\|_{\K^{m-m_{\ell}-1/2,m/2-m_{\ell}}} \leq \\
\|\check{\omega}\kappa_{\varrho}^{-1}\gamma\tilde{B}_{\ell,m}\kappa_{\varrho}\tilde{\omega}\|_{\L(\K^{m,-m/2},\K^{m-m_{\ell}-1/2,m/2-m_{\ell}})}
\|{\omega}L_{\varrho}\tilde u\|_{\K^{m,-m/2}}
\end{gather*}
for cut-off functions $\omega \prec \tilde{\omega} \prec \check{\omega}$.
Lemma~\ref{Ltheta} implies
$$
\|{\omega}L_{\varrho}\tilde u\|_{\K^{m,-m/2}}\le \textup{const} 
\|\omega \tilde u\|_{\Dom_{\max}},
$$
and so
\begin{equation*}
\|\kappa_{\varrho}^{-1}\gamma\tilde{B}_{\ell,m}\kappa_{\varrho}\bigl({\omega}L_{\varrho}\tilde u\bigr)\|_{
\K^{m-m_{\ell}-1/2,m/2-m_{\ell}}} \leq \textup{const} \|\omega \tilde u\|_{\Dom_{\max}}
\end{equation*}
since $\|\check{\omega}\kappa_{\varrho}^{-1}\gamma \tilde{B}_{\ell,m}\kappa_{\varrho}\tilde{\omega}\|=O(\varrho^{m_{\ell}-m})$ as $\varrho\to\infty$.
Thus \eqref{Kestimate3} is proved. 
\end{proof}

\bigskip 
\begin{proof}[Proof of Theorem~\ref{Maintheorem}]
Fix some complement $\Sing_{\max}$ of $\Dom_{\min}\binom{A}{T}$ in $\Dom_{\max}\binom{A}{T}$
and let $\Sing \subset \Sing_{\max}$ be a subspace such that
$\Dom\binom{A}{T} = \Dom_{\min}\binom{A}{T}\oplus \Sing$. With respect to this decomposition we
may write the boundary value problem \eqref{BVPCone} as
\begin{align*}
{\mathcal A}_{\Dom}(\lambda) = \binom{A-\lambda}{T} &= \begin{pmatrix} {\mathcal A}_{\Dom_{\min}}(\lambda) & {\mathcal A}(\lambda)|_{\Sing} \end{pmatrix} =
\begin{pmatrix} 
(A -\lambda)|_{\Dom_{\min}} & (A -\lambda)|_{\Sing} \\
T|_{\Dom_{\min}} & T|_{\Sing}
\end{pmatrix} \\ &:
\begin{array}{c} \Dom_{\min}\binom{A}{T} \\ \oplus \\ \Sing \end{array} \to
\begin{array}{c} x^{-m/2}L_b^2(\overline{M};E) \\ \oplus \\ \bigoplus\limits_{j=1}^K x^{m/2-m_j}H_b^{m-m_j-1/2}(\overline{N};F_j)\end{array}.
\end{align*}   
Let $d''=\dim \Sing$. Under the assumptions of Theorem~\ref{Maintheorem}
we may apply Theorem~\ref{MainTheoremParametrix} and obtain the existence of
a parametrix ${\mathcal B}(\lambda)$ and a 
generalized Green remainder ${\mathcal K}(\lambda) = \begin{pmatrix} {\mathcal K}_0(\lambda) \\ 0 \end{pmatrix}$
of order $m$ such that
\begin{equation*}
\begin{pmatrix} 
(A -\lambda)|_{\Dom_{\min}} & {\mathcal K}_0(\lambda) \\
T|_{\Dom_{\min}} & 0
\end{pmatrix} :
\begin{array}{c} \Dom_{\min}\binom{A}{T} \\ \oplus \\ \C^{d''} \end{array} 
\to
\begin{array}{c} x^{-m/2}L_b^2(\overline{M};E) \\ \oplus \\ \bigoplus\limits_{j=1}^K x^{m/2-m_j}H_b^{m-m_j-1/2}(\overline{N};F_j)\end{array}
\end{equation*}   
is invertible for $\lambda \in \Lambda$ sufficiently large with inverse
\begin{equation}\label{LeftDParametrix}
\begin{pmatrix} 
(A -\lambda)|_{\Dom_{\min}} & {\mathcal K}_0(\lambda) \\
T|_{\Dom_{\min}} & 0
\end{pmatrix}^{-1} = 
\left(\begin{array}{c}
{\mathcal B}(\lambda) \\ \hline \Vsp  {\mathcal T}(\lambda)
\end{array}\right), 
\end{equation}
where ${\mathcal T}(\lambda) = \begin{pmatrix} {\mathcal T}_0(\lambda) & \cdots & {\mathcal T}_K(\lambda)\end{pmatrix}$ 
is a matrix of generalized Green remainders of orders $-m,-m_1,\ldots,-m_K$ and type zero. Since
\begin{align*}
\id =
\begin{pmatrix} 
{\mathcal B}(\lambda) \\  {\mathcal T}(\lambda)
\end{pmatrix}
\begin{pmatrix} 
{\mathcal A}_{\Dom_{\min}}(\lambda) & {\mathcal K}(\lambda)
\end{pmatrix} 
= 
\begin{pmatrix} 
{\mathcal B}(\lambda){\mathcal A}_{\Dom_{\min}}(\lambda) & {\mathcal B}(\lambda){\mathcal K}(\lambda) \\
{\mathcal T}(\lambda){\mathcal A}_{\Dom_{\min}}(\lambda) & {\mathcal T}(\lambda){\mathcal K}(\lambda)
\end{pmatrix}, 
\end{align*}
we have ${\mathcal B}(\lambda){\mathcal A}_{\Dom_{\min}}(\lambda) = 1$, 
${\mathcal T}(\lambda){\mathcal A}_{\Dom_{\min}}(\lambda)=0$, and ${\mathcal T}(\lambda){\mathcal K}(\lambda)=1$. Then
\begin{equation}\label{CompWithDParametrix}
\begin{pmatrix} 
{\mathcal B}(\lambda) \\ {\mathcal T}(\lambda)
\end{pmatrix}
\begin{pmatrix} 
{\mathcal A}_{\Dom_{\min}}(\lambda) & {\mathcal A}(\lambda)|_{\Sing}
\end{pmatrix} = 
\begin{pmatrix} 
1 & {\mathcal B}(\lambda){\mathcal A}(\lambda)|_{\Sing} \\  
0 & {\mathcal T}(\lambda){\mathcal A}(\lambda)|_{\Sing}
\end{pmatrix} 
\end{equation}
which implies that $\begin{pmatrix} {\mathcal A}_{\Dom_{\min}}(\lambda) & {\mathcal A}(\lambda)|_{\Sing} \end{pmatrix}$
is invertible if and only if 
\begin{equation}\label{FLambda}
F(\lambda) = {\mathcal T}(\lambda){\mathcal A}(\lambda) = {\mathcal T}(\lambda)\binom{A-\lambda}{T}: \Sing\to \C^{d''}
\end{equation}    
is invertible. Moreover, we get the explicit representation
\begin{equation}\label{TheResolvent}
\binom{A-\lambda}{T}^{-1} = {\mathcal A}_{\Dom}(\lambda)^{-1} =
{\mathcal B}(\lambda) + \bigl(1-{\mathcal B}(\lambda){\mathcal A}(\lambda)\bigr)
F(\lambda)^{-1}{\mathcal T}(\lambda),
\end{equation}
and \eqref{InverseDarstellung} follows from Theorem \ref{MainTheoremParametrix}.

As $F(\lambda) = {\mathcal T}(\lambda){\mathcal A}(\lambda)$ and $1 - {\mathcal B}(\lambda){\mathcal A}(\lambda)$ vanish
on $\Dom_{\min}\binom{A}{T}$ for large $\lambda$, they descend to operators 
$F(\lambda) : \tilde\Sing_{\max} \to \C^{d''}$ and 
$1 - {\mathcal B}(\lambda){\mathcal A}(\lambda) : \tilde\Sing_{\max} \to \Dom_{\max}\binom{A}{T}$. 
If $\tilde\Sing=\Dom\binom{A}{T}/\Dom_{\min}\binom{A}{T}$, then the invertibility of \eqref{FLambda} 
is equivalent to the invertibility of
\begin{equation*}
F(\lambda) : \tilde\Sing \to \C^{d''},
\end{equation*}
and in this case, \eqref{TheResolvent} still makes sense in this context.

Let $q : \Dom_{\max}\binom{A}{T} \to \tilde \Sing_{\max}$ be the canonical projection.
The inverses ${\mathcal A}_{\Dom}(\lambda)^{-1}$ and $F(\lambda)^{-1} : \C^{d''} \to \tilde\Sing \subset
\tilde\Sing_{\max}$ are related by the formulas
\begin{align*}
F(\lambda)^{-1} = q{\mathcal A}_{\Dom}(\lambda)^{-1}{\mathcal K}(\lambda) &: 
\C^{d''} \to \tilde\Sing_{\max}, \\
q{\mathcal A}_{\Dom}(\lambda)^{-1} = F(\lambda)^{-1}{\mathcal T}(\lambda) &: 
\begin{array}{c} x^{-m/2}L_b^2(\overline{M};E)  \\ \oplus \\ \bigoplus\limits_{j=1}^Kx^{m/2-m_j}H_b^{m-m_j-1/2}(\overline{N};F_j)
\end{array} \to \tilde\Sing_{\max}
\end{align*}
in view of ${\mathcal T}(\lambda){\mathcal K}(\lambda) = 1$.

Under the assumptions of Theorem \ref{Maintheorem} we prove that
$F(\lambda) : \tilde\Sing \to \C^{d''}$ is invertible for large $\lambda \in \Lambda$, 
and that the inverse satisfies the norm estimate 
\begin{equation}\label{Festimate}
\|\tilde{\kappa}_{[\lambda]^{1/m}}^{-1}
F(\lambda)^{-1}\|_{\L(\C^{d''},\tilde\Sing_{\max})} = O(1)
\;\text{ as } |\lambda| \to \infty.
\end{equation}
Observe that the parametrix construction from Theorem~\ref{MainTheoremParametrix}
gives the relation
$$
\begin{pmatrix} 
(A_{\wedge} -\lambda)|_{\Dom_{\wedge,\min}} & {\mathcal K}_{0,\wedge}(\lambda) \\
T_{\wedge}|_{\Dom_{\wedge,\min}} & 0
\end{pmatrix}^{-1} = 
\begin{pmatrix}
{\mathcal B}_{\wedge}(\lambda) \\ {\mathcal T}_{\wedge}(\lambda)
\end{pmatrix}
$$
for the $\kappa$-homogeneous principal parts of \eqref{LeftDParametrix}. 
Thus with the same reasoning as above we conclude that
$$
{\mathcal A}_{\wedge,\Dom_{\wedge}}(\lambda) = \binom{A_{\wedge}-\lambda}{T_{\wedge}} : \Dom_{\wedge}\binom{A_{\wedge}}{T_{\wedge}}
\to \begin{array}{c} x^{-m/2}L_b^2(\overline{Y}^{\wedge};E) \\ \oplus \\ \bigoplus\limits_{j=1}^K
\K^{m-m_j-1/2,m/2-m_j}((\partial\overline{Y})^{\wedge};F_j) \end{array}
$$
is invertible if and only if the restriction of the induced operator
$$
F_{\wedge}(\lambda) = {\mathcal T}_{\wedge}(\lambda){\mathcal A}_{\wedge}(\lambda) = {\mathcal T}_{\wedge}(\lambda)\binom{A_{\wedge} - \lambda}{T_{\wedge}} : 
\tilde\Sing_{\wedge,\max}\to \C^{d''}
$$
to $\tilde\Sing_{\wedge}=\Dom_{\wedge}\binom{A_{\wedge}}{T_{\wedge}}/\Dom_{\wedge,\min}\binom{A_{\wedge}}{T_{\wedge}}$ is invertible. 
Let $q_{\wedge}: \Dom_{\wedge,\max}\binom{A_{\wedge}}{T_{\wedge}} \to \tilde\Sing_{\wedge,\max}$ be the 
canonical projection. From the relations
\begin{align*}
F_{\wedge}(\lambda)^{-1} = q_{\wedge}{\mathcal A}_{\wedge,\Dom_{\wedge}}(\lambda)^{-1}
{\mathcal K}_{\wedge}(\lambda) &: \C^{d''} \to \tilde\Sing_{\wedge,\max}, \\
q_{\wedge}{\mathcal A}_{\wedge,\Dom_{\wedge}}(\lambda)^{-1} = F_{\wedge}(\lambda)^{-1}
{\mathcal T}_{\wedge}(\lambda) &: \begin{array}{c} x^{-m/2}L_b^2(\overline{Y}^{\wedge};E) \\ \oplus \\
\bigoplus\limits_{j=1}^K \K^{m-m_j-1/2,m/2-m_j}((\partial\overline{Y})^{\wedge};F_j) \end{array}
\to \tilde\Sing_{\wedge,\max},
\end{align*}
and Proposition \ref{WedgeSectorMinimalGrowth}, we deduce that our assumption about the resolvent of
$A_\wedge$ under the boundary condition $T_{\wedge}u = 0$ is equivalent to the invertibility
of ${\mathcal A}_{\wedge,\Dom_{\wedge}}(\lambda)$ with domain $\Dom_{\wedge}\binom{A_{\wedge}}{T_{\wedge}}$,
where $\Dom_{\wedge} = \theta\bigl(\Dom\bigr)$ is the associated domain to $\Dom$, and
\begin{equation}\label{Fdachestimate}
\|\kappa^{-1}_{|\lambda|^{1/m}}F_{\wedge}(\lambda)^{-1}\|_{\L(\C^{d''},
\tilde\Sing_{\wedge,\max})} = O(1)
\;\text{ as } |\lambda| \to \infty.
\end{equation}
Note that $\|\kappa^{-1}_{|\lambda|^{1/m}}{\mathcal K}_{0,\wedge}(\lambda)\| = O(|\lambda|)$ and
$\|{\mathcal T}_{k,\wedge}(\lambda)\kappa_{|\lambda|^{1/m}}\| = O(|\lambda|^{-m_k/m})$ as $|\lambda| \to \infty$
for $k=0,\ldots,m_K$, where $m_0 = m$.
 
Write the operator $F(\lambda)\theta^{-1}
F_{\wedge}(\lambda)^{-1} : \C^{d''} \to \C^{d''}$ as
$$
F(\lambda)\theta^{-1}F_{\wedge}(\lambda)^{-1} = 
1 +\bigl(F(\lambda) - F_{\wedge}(\lambda)\theta\bigr)
\tilde{\kappa}_{|\lambda|^{1/m}} \theta^{-1}\kappa^{-1}_{|\lambda|^{1/m}}
F_{\wedge}(\lambda)^{-1},
$$
and let 
\begin{equation*}
R(\lambda)=\bigl(F(\lambda) - F_{\wedge}(\lambda)\theta\bigr)
\tilde{\kappa}_{|\lambda|^{1/m}} \theta^{-1}\kappa^{-1}_{|\lambda|^{1/m}}
F_{\wedge}(\lambda)^{-1}.
\end{equation*}
We will prove in Lemma \ref{FFdach} that
$$
\|(F(\lambda) - F_{\wedge}(\lambda)\theta)
\tilde{\kappa}_{|\lambda|^{1/m}}\|_{\L(\tilde\Sing_{\max},\C^{d''})} \to 0
\;\text{ as } |\lambda| \to \infty.
$$
Thus together with \eqref{Fdachestimate} we obtain that
$\|R(\lambda)\| \to 0$ as $|\lambda| \to \infty$.
Hence $1 + R(\lambda)$ is invertible for large $|\lambda| > 0$, and the
inverse is of the form $1 + \tilde{R}(\lambda)$ with $\|\tilde{R}(\lambda)\|
\to 0$ as $|\lambda| \to \infty$.
This shows that $F(\lambda) : \tilde\Sing \to \C^{d''}$ is invertible from 
the right for large $\lambda$, and by \eqref{Fdachestimate} the right inverse 
$\theta^{-1}F_{\wedge}(\lambda)^{-1}(1 + \tilde{R}(\lambda))$ satisfies the 
norm estimate \eqref{Festimate}. Since
$$
\dim \tilde\Sing = \dim \tilde\Sing_{\wedge} = d'',
$$
we conclude that $F(\lambda)$ is also injective, and
so the invertibility of $F(\lambda)$ is proved.
In particular, the operator ${\mathcal A}_{\Dom}(\lambda)$ is invertible for large $\lambda$, and
consequently also
$$
A_{\Dom} - \lambda : \Dom(A_T) \to x^{-m/2}L_b^2(\overline{M};E)
$$
is invertible for large $\lambda \in \Lambda$. It remains to show the resolvent estimates for
$\bigl(A_{\Dom} - \lambda\bigr)^{-1}$ as $|\lambda| \to \infty$, see Definition \ref{Sectorminimalgrowthdef}.

In order to prove these estimates we make use of the family $\tilde{K}(\varrho)$ from
Lemma~\ref{kappaSchlangelift} and the representation \eqref{TheResolvent}. We may write
\begin{align*}
(A_{\Dom} - \lambda)^{-1} &= {\mathcal B}_T(\lambda) + (1 - {\mathcal B}(\lambda){\mathcal A}(\lambda)) 
\tilde{K}(|\lambda|^{1/m})\tilde{\kappa}^{-1}_{|\lambda|^{1/m}}
F(\lambda)^{-1}{\mathcal T}_0(\lambda) \\
&={\mathcal B}_T(\lambda) + \tilde{K}(|\lambda|^{1/m})\tilde{\kappa}^{-1}_{|\lambda|^{1/m}}
F(\lambda)^{-1}{\mathcal T}_0(\lambda) \\
&\hspace*{6em} - {\mathcal B}(\lambda){\mathcal A}(\lambda)\tilde{K}(|\lambda|^{1/m})
\tilde{\kappa}^{-1}_{|\lambda|^{1/m}}F(\lambda)^{-1}{\mathcal T}_0(\lambda).
\end{align*}
By construction of the parametrix we have 
$\|{\mathcal B}_T(\lambda)\|_{\L(x^{-m/2}L^2_b)} = O(|\lambda|^{-1})$ as 
$|\lambda| \to \infty$. In view of $\|{\mathcal T}_0(\lambda)\|_{\L(x^{-m/2}L_b^2,\C^{d''})} 
= O(|\lambda|^{-1})$ and \eqref{Festimate} we further obtain
$$
\|\tilde{\kappa}^{-1}_{|\lambda|^{1/m}}F(\lambda)^{-1}
{\mathcal T}_0(\lambda)\|_{\L(x^{-m/2}L_b^2,\tilde\Sing_{\max})}= O(|\lambda|^{-1})
\;\text{ as } |\lambda| \to \infty,
$$
and consequently, using \eqref{Kestimate1}, we get
$$
\|\tilde{K}(|\lambda|^{1/m})\tilde{\kappa}^{-1}_{|\lambda|^{1/m}}
F(\lambda)^{-1}{\mathcal T}_0(\lambda)\|_{\L(x^{-m/2}L_b^2)}= O(|\lambda|^{-1})
\;\text{ as } |\lambda| \to \infty.
$$
On the other hand, by \eqref{Festimate} and the estimates \eqref{Kestimate1}--\eqref{Kestimate3} we have
\begin{align}
\label{Ab1}\|\kappa_{|\lambda|^{1/m}}^{-1}(A-\lambda)\tilde{K}(|\lambda|^{1/m})\tilde{\kappa}^{-1}_{|\lambda|^{1/m}}
F(\lambda)^{-1}\|_{\L(\C^{d''},x^{-m/2}L_b^2)} &= O(|\lambda|), \\
\label{Ab2}\|\kappa_{|\lambda|^{1/m}}^{-1}\gamma B_{j}\tilde{K}(|\lambda|^{1/m})\tilde{\kappa}^{-1}_{|\lambda|^{1/m}}
F(\lambda)^{-1}\|_{\L(\C^{d''},\K^{m-m_j-1/2,m/2-m_j})} &= O(|\lambda|^{m_j/m})
\end{align}
as $|\lambda| \to \infty$ for all $j=1,\ldots,K$.

Let $\tilde{\omega} \in C_0^{\infty}([0,1))$ be any cut-off function near zero. For large $\lambda$
we may write
{\small
$$
{\mathcal B}(\lambda){\mathcal A}(\lambda)\tilde{K}(|\lambda|^{1/m})
\tilde{\kappa}^{-1}_{|\lambda|^{1/m}}F(\lambda)^{-1}{\mathcal T}_0(\lambda) =
{\mathcal B}(\lambda)\tilde{\omega}{\mathcal A}(\lambda)\tilde{K}(|\lambda|^{1/m})
\tilde{\kappa}^{-1}_{|\lambda|^{1/m}}F(\lambda)^{-1}{\mathcal T}_0(\lambda),
$$}
and for the components of the parametrix
$$
{\mathcal B}(\lambda) = \begin{pmatrix} {\mathcal B}_T(\lambda) & {\mathcal B}_1(\lambda) & \cdots & {\mathcal B}_K(\lambda)\end{pmatrix}
$$
we have by construction the norm estimates
\begin{align}
\label{Ab3}\|{\mathcal B}_T(\lambda)\tilde{\omega}\kappa_{|\lambda|^{1/m}}\|_{\L(x^{-m/2}L_b^2(\overline{Y}^{\wedge};E),x^{-m/2}L_b^2(\overline{M};E))} &= O(|\lambda|^{-1}), \\
\label{Ab4}\|{\mathcal B}_j(\lambda)\tilde{\omega}\kappa_{|\lambda|^{1/m}}\|_{\L(\K^{m-m_j-1/2,m/2-m_j},x^{-m/2}L_b^2(\overline{M};E))} &= O(|\lambda|^{-m_j/m})
\end{align}
as $|\lambda| \to \infty$ for $j=1,\ldots,K$ (the operator family $g(\lambda) = \omega'{\mathcal B}_j(\lambda)\tilde{\omega}$
satisfies the estimate \eqref{opsymbabsch} with $\mu = -m_j$ in $\K^{m-m_j-1/2,m/2-m_j} \to
\Dom_{\wedge,\min}$ with respect to the normalized dilation group action $\kappa_{\varrho}$ on both
spaces).

In view of \eqref{Ab1}--\eqref{Ab4} and $\|{\mathcal T}_0(\lambda)\|_{\L(x^{-m/2}L_b^2,\C^{d''})} = O(|\lambda|^{-1})$
we now conclude that, as $|\lambda| \to \infty$,
$$
\|{\mathcal B}(\lambda){\mathcal A}(\lambda)\tilde{K}(|\lambda|^{1/m})
\tilde{\kappa}^{-1}_{|\lambda|^{1/m}}F(\lambda)^{-1}{\mathcal T}_0(\lambda)\|_{
\L(x^{-m/2}L_b^2)} = O(|\lambda|^{-1}).
$$

Summing up, we obtain
$$
\|(A_{\Dom} - \lambda)^{-1}\|_{\L(x^{-m/2}L_b^2)} = O(|\lambda|^{-1})
\;\text{ as } |\lambda| \to \infty
$$
as desired.
\end{proof}

The following lemma completes the proof of Theorem \ref{Maintheorem}.

\begin{lemma}\label{FFdach}
With the notation of the proof of Theorem \ref{Maintheorem}, let
\begin{align*}
F(\lambda) = {\mathcal T}(\lambda){\mathcal A}(\lambda) &: \tilde\Sing_{\max} \to \C^{d''}, \\
F_{\wedge}(\lambda) = {\mathcal T}_{\wedge}(\lambda){\mathcal A}_{\wedge}(\lambda) &: 
\tilde\Sing_{\wedge,\max} \to \C^{d''}.
\end{align*}
Then
\begin{equation}\label{FFdachNull}
\|(F(\lambda) - F_{\wedge}(\lambda)\theta)
\tilde{\kappa}_{|\lambda|^{1/m}}\|_{\L(\tilde\Sing_{\max},\C^{d''})} \to 0
\;\text{ as } |\lambda| \to \infty.
\end{equation}
\end{lemma}
\begin{proof}
For proving \eqref{FFdachNull} it is sufficient to consider the restrictions
$$
(F(\lambda) - F_{\wedge}(\lambda)\theta)\tilde{\kappa}_{|\lambda|^{1/m}} :
\tilde\Sing_{\sigma_0} \to \C^{d''}
$$
for all $\sigma_0 \in \Sigma$. First of all, observe that	
\begin{gather*}
F(\lambda)\tilde{\kappa}_{|\lambda|^{1/m}} 
= {\mathcal T}(\lambda){\mathcal A}(\lambda)\tilde{K}(|\lambda|^{1/m}), \;\text{ and }\\
F_{\wedge}(\lambda)\theta \tilde{\kappa}_{|\lambda|^{1/m}} 
= F_{\wedge}(\lambda)\kappa_{|\lambda|^{1/m}}\theta
= {\mathcal T}_\wedge(\lambda){\mathcal A}_{\wedge}(\lambda)\kappa_{|\lambda|^{1/m}}\omega\theta
\end{gather*}
with the operator family $\tilde{K}(\varrho)=\omega(\varrho x)
\tilde\kappa_{\varrho}$ from Lemma \ref{kappaSchlangelift}. 
If $\omega_0 \in C_0^{\infty}([0,1))$ is a cut-off function near 
zero with $\omega \prec \omega_0$, then
\begin{align*}
(F(\lambda) - F_{\wedge}(\lambda)\theta)\tilde{\kappa}_{|\lambda|^{1/m}}
&= {\mathcal T}(\lambda){\mathcal A}(\lambda)\tilde{K}(|\lambda|^{1/m})
-  {\mathcal T}_\wedge(\lambda){\mathcal A}_\wedge(\lambda)\kappa_{|\lambda|^{1/m}}\omega\theta\\
&= {\mathcal T}(\lambda)\omega_0{\mathcal A}(\lambda)\tilde{K}(|\lambda|^{1/m})
-  {\mathcal T}_\wedge(\lambda)\omega_0{\mathcal A}_\wedge(\lambda)\kappa_{|\lambda|^{1/m}}
   \omega\theta\\
&= {\mathcal T}(\lambda)\omega_0\left({\mathcal A}(\lambda)\tilde K(|\lambda|^{1/m})
-{\mathcal A}_\wedge(\lambda)\kappa_{|\lambda|^{1/m}}\omega\theta\right)\\
&\hspace*{4em} +({\mathcal T}(\lambda)-{\mathcal T}_\wedge(\lambda))\omega_0{\mathcal A}_\wedge(\lambda)
\kappa_{|\lambda|^{1/m}}\omega\theta.
\end{align*}
Now
{\small
$$
({\mathcal T}(\lambda)-{\mathcal T}_\wedge(\lambda))\omega_0{\mathcal A}_\wedge(\lambda)
\kappa_{|\lambda|^{1/m}}\omega\theta =
\Bigl(({\mathcal T}(\lambda)-{\mathcal T}_\wedge(\lambda))\omega_0\kappa_{|\lambda|^{1/m}}\Bigr)
\Bigl(\kappa^{-1}_{|\lambda|^{1/m}}{\mathcal A}_\wedge(\lambda)
\kappa_{|\lambda|^{1/m}}\omega\theta\Bigr),
$$}
and consequently this term is $o(1)$ as $|\lambda| \to \infty$. Recall that
$$
\|({\mathcal T}_k(\lambda)-{\mathcal T}_{k,\wedge}(\lambda))\omega_0\kappa_{|\lambda|^{1/m}}\| =
O(|\lambda|^{-(m_k-1)/m})
$$
for $k = 0,\ldots,K$, where $m_0 = m$. On the other hand, we have
\begin{gather*}
{\mathcal T}(\lambda)\omega_0\Bigl({\mathcal A}(\lambda)\tilde K(|\lambda|^{1/m})
-{\mathcal A}_\wedge(\lambda)\kappa_{|\lambda|^{1/m}}\omega\theta\Bigr) = \\
\Bigl({\mathcal T}(\lambda)\omega_0\kappa_{|\lambda|^{1/m}}\Bigr)\Bigl(\kappa_{|\lambda|^{1/m}}^{-1}
{\mathcal A}(\lambda)\tilde K(|\lambda|^{1/m})
-\kappa_{|\lambda|^{1/m}}^{-1}{\mathcal A}_\wedge(\lambda)\kappa_{|\lambda|^{1/m}}\omega\theta\Bigr).
\end{gather*}
By \eqref{PhiAPhiL} and \eqref{Theta0Estimate} in Lemma~\ref{kappaSchlangelift} (see also Lemma 6.20
in \cite{GKM2}) and Lemma~\ref{Ltheta} it follows that each summand in this matrix multiplication
is $o(1)$ in the norm as $|\lambda| \to \infty$, and so the lemma follows.
\end{proof}

Finally, we want to point out that under the assumptions of Theorem \ref{Maintheorem} we get
the existence of the resolvent of $A$ with polynomial bounds for the norm also for realizations in
Sobolev spaces of arbitrary smoothness $s > -\frac{1}{2}$. The proof follows along the lines
of this section. The advantage in this case is that we need not be as precise with the bounds as for the
case of $x^{-m/2}L^2_b$-realizations.

\begin{theorem}\label{RMG-HigherRegularity}
Let \eqref{BVPCone} be $c$-elliptic with parameter in the closed sector $\Lambda \subset \C$, and consider
the unbounded operator $A$ in $x^{-m/2}H^s_b(\overline{M};E)$ under the boundary condition $Tu = 0$
on some intermediate domain $\Dom^s_{\min}(A_T) \subset \Dom^s(A_T) \subset \Dom^s_{\max}(A_T)$,
where $s > -\frac{1}{2}$.

Let $\Dom_{\wedge}(A_{\wedge,T_{\wedge}}) = \theta(\Dom(A_T))$ be the associated domain for the
model operator $A_{\wedge}$ under the boundary condition $T_{\wedge}u = 0$ in $x^{-m/2}L_b^2(\overline{Y}^{\wedge};E)$
according to Proposition \ref{SingularFunctionManifold}, and assume that $\Lambda$ is a sector
of minimal growth for $A_{\wedge}$ with this domain.

Then
$$
A_{\Dom^s} - \lambda : \Dom^s(A_T) \to x^{-m/2}H_b^s(\overline{M};E)
$$
is invertible for large $\lambda \in \Lambda$, and the resolvent can be written in the form
\eqref{InverseDarstellung}. Moreover, there exists $M(s) \in \R$
such that
$$
\|\bigl(A_{\Dom^s} - \lambda\bigr)^{-1}\|_{\L(x^{-m/2}H_b^s)} = O(|\lambda|^{M(s)})
$$
as $|\lambda| \to \infty$.
\end{theorem}



\begin{thebibliography}{10}

\bibitem{AV63}
M.~Agranovich and M. Vishik, \emph{Elliptic problems with a
parameter and parabolic problems of general type}, Russ. Math. Surveys
\textbf{19} (1963), 53--159.

\bibitem{AmLauNis}
B.~Ammann, R.~Lauter, and V.~Nistor, \emph{Pseudodifferential operators on manifolds
with a Lie structure at infinity}, preprint 2003 (math.AP/0304044 at arXiv.org), to
appear in Ann. of Math.

\bibitem{BdM71}
L.~Boutet de Monvel, \emph{Boundary problems for pseudo-differential operators},
Acta Math. \textbf{126} (1971), 11--51.

\bibitem{CoSeiSchr04}
S.~Coriasco, J.~Seiler, and E.~Schrohe, \emph{Realizations of differential
operators on conic manifolds with boundary}, preprint 2004 (math.AP/0401395 at arXiv.org).

\bibitem{Dauge}
M.~Dauge, \emph{Elliptic boundary value problems on corner domains}, vol.~1341 of \emph{Lecture
Notes in Mathematics}. Springer-Verlag, Berlin, 1988.

\bibitem{Eskin}
G.~Eskin, \emph{Boundary value problems for elliptic pseudodifferential equations}, Translations of Mathematical Monographs, vol.~52. American Mathematical Society, Providence, R.I., 1981.

\bibitem{Gil}
J.~Gil, \emph{Full asymptotic expansion of the heat trace for non-self-adjoint elliptic cone
operators}, Math. Nachr. \textbf{250} (2003), 25--57.

\bibitem{GKM1}
J.~Gil, T.~Krainer, and G.~Mendoza, \emph{Geometry and spectra of 
closed extensions of elliptic cone operators}, preprint 2004 (math.AP/0410178 at arXiv.org),
to appear in Canadian J. Math.

\bibitem{GKM2}
\bysame, \emph{Resolvents of elliptic cone operators}, preprint 2004 (math.AP/0410176 at arXiv.org),
to appear in J. Funct. Anal.

\bibitem{GKM3}
\bysame, \emph{On rays of minimal growth for elliptic cone operators},
preprint 2005 (math.SP/0511727 at arXiv.org).

\bibitem{GiMe01}
J.~Gil and G.~Mendoza, \emph{Adjoints of elliptic cone operators},
Amer. J. Math. \textbf{125} (2003), no.~2, 357--408.

\bibitem{GrubbBuch} 
G.~Grubb, \emph{Functional calculus of pseudodifferential boundary problems},
2nd ed., Progress in Mathematics, vol.~65. Birkh{\"a}user, Basel, 1996.

\bibitem{GruSe95}
G.~Grubb and R.~Seeley, \emph{Weakly parametric pseudodifferential operators
and Atiyah-Patodi-Singer boundary problems}, Invent. Math. \textbf{121} (1995),
no. 3, 481--529.

\bibitem{JerKen}
D.~Jerison and C.~Kenig, \emph{The inhomogeneous Dirichlet problem in Lipschitz
domains}, J. Funct. Anal. \textbf{130} (1995), no.~1, 161--219.

\bibitem{KaSchu03}
D.~Kapanadze and B.-W.~Schulze, \emph{Crack Theory and Edge Singularities},
Mathematics and its Applications, vol.~561, Kluwer Academic Publishers,
Dordrecht-Boston-London, 2003.

\bibitem{K1} 
V.~Kondrat'ev, \emph{Boundary problems for
elliptic equations in domains with conical or angular points},
Trans. Mosc. Math. Soc. \textbf{16} (1967), 227--313.

\bibitem{KozMazRos}
V.~Kozlov, V.~Maz'ya, and J.~Ro{\ss}mann, \emph{Spectral problems associated with corner
singularities of solutions to elliptic equations}, vol.~85 of \emph{Mathematical Surveys and
Monographs}. American Mathematical Society, Providence, RI, 2001.

\bibitem{Le97} 
M.~Lesch, \emph{Operators of {F}uchs type, conical
singularities, and asymptotic methods}, Teubner-Texte zur Math.
vol 136, B.G. Teubner, Stuttgart, Leipzig, 1997.

\bibitem{LioMag}
J.-L.~Lions and E.~Magenes, \emph{Non-homogeneous boundary value problems and applications}, Grundlehren der mathematischen Wissenschaften, vol.~181--183. Springer-Verlag, New York-Heidelberg, 1972, 1973.

\bibitem{LoRes01}
P.~Loya, \emph{On the resolvent of differential operators on conic
manifolds}, Comm. Anal. Geom. \textbf{10} (2002), no.~5, 877--934.

\bibitem{MazNazPlam}
V.~Maz'ya, S.~Nazarov, and B.~Plamenevsky, \emph{Asymptotic theory of elliptic boundary value
problems in singularly perturbed domains {I}, {II}}, vol.~111 and 112 of \emph{Operator Theory: Advances
and Applications}, Birkh{\"a}user Verlag, Basel, 2000.

\bibitem{MimiTay}
D.~Mitrea, M.~Mitrea, and M.~Taylor, \emph{Layer potentials, the Hodge Laplacian, and global boundary
problems in nonsmooth Riemannian manifolds}, vol.~713 of \emph{Memoirs of the American Mathematical
Society}. American Mathematical Society, Providence, RI, 2001.

\bibitem{MiNistor}
M.~Mitrea and V.~Nistor, \emph{Boundary value problems and layer potentials on manifolds with
cylindrical ends}, preprint 2004 (math.AP/0410186 at arXiv.org).

\bibitem{MiTayVas}
M.~Mitrea, M.~Taylor, and A.~Vasy, \emph{Lipschitz domains, domains with corners, and the
Hodge Laplacian}, Comm. Partial Differential Equations \textbf{30} (2005), 1445--1462.

\bibitem{Nistor}
V.~Nistor, \emph{Pseudodifferential operators on non-compact manifolds and analysis on polyhedral
domains}, Contemp. Mathematics \textbf{366} (2005), 307--328.

\bibitem{RBM1} 
R.~Melrose, \emph{Transformation of boundary value
problems}, Acta Math.  \textbf{147} (1981), no. 3-4, 149--236.

\bibitem{RBM2} 
\bysame, \emph{The Atiyah-Patodi-Singer index
theorem}, Research Notes in Mathematics, A~K~Peters, Ltd.,
Wellesley, MA, 1993.

\bibitem{MM} 
R.~Melrose and G.~Mendoza, \emph{Elliptic operators
of totally characteristic type}, MSRI Preprint, 1983.

\bibitem{Schrohe99}
E.~Schrohe, \emph{Fr{\'e}chet algebra techniques for boundary value problems on noncompact
manifolds: Fredholm criteria and functional calculus via spectral invariance},
Math.~Nachr. \textbf{199} (1999), 145--185.

\bibitem{SchroSchu}
E.~Schrohe and B.-W.~Schulze, \emph{Boundary value problems in {B}outet de {M}onvel's algebra for manifolds with conical singularities {I},{II}},
Math. Top., vol.~5 and 8, Akademie Verlag, Berlin, (1994, 1995), pp.~97--209, 70--205.

\bibitem{SzNorthHolland} 
B.-W.~Schulze, \emph{Pseudo-differential operators on manifolds with
singularities}, Studies in Mathematics and its Applications, vol.~24. 
North-Holland Publishing Co., Amsterdam, 1991.

\bibitem{SzWiley98}
\bysame, \emph{Boundary value problems and singular
pseudo-differential operators}, Pure and Applied Mathematics.
John Wiley \& Sons, Ltd., Chichester, 1998.

\bibitem{Seeley}
R.~Seeley, \emph{The resolvent of an elliptic boundary problem},
Amer. J. Math. \textbf{91} (1969), 889--920. 

\bibitem{SeeleyBook}
\bysame, \emph{Topics in Pseudo-Differential Operators},
CIME Conference of Pseudo-Differential Operators 1968,
Edizione Cremonese, Roma, 1969, pp.~169--305.

\end{thebibliography}
\end{document}